\documentclass{irmadegm}
\usepackage[psamsfonts]{amssymb}
\usepackage{amsmath,cmmib57}
\usepackage{amsthm}
\usepackage{latexsym}
\usepackage{xspace}

\makeatletter
\def\thmhead@plain#1#2#3{%
  \thmname{#1}\thmnumber{\@ifnotempty{#1}{ }#2}%
  \thmnote{ \the\thm@notefont(#3)}}
\let\thmhead\thmhead@plain
\def\swappedhead#1#2#3{%
  \thmnumber{#2}\thmname{\@ifnotempty{#2}{. }#1}%
  \thmnote{ \the\thm@notefont(#3)}}
\thm@notefont{\rm} \makeatother

\theoremstyle{definition} 

 \newtheorem{definition}{Definition}[section]
 \newtheorem{remark}[definition]{Remark}
 \newtheorem{example}[definition]{Example}
 
 \newtheorem{terminology}[definition]{Terminology}

\theoremstyle{plain}      

 \newtheorem{proposition}[definition]{Proposition}
 \newtheorem{theorem}[definition]{Theorem}
 \newtheorem{corollary}[definition]{Corollary}
 \newtheorem{lemma}[definition]{Lemma}

\newcommand{\mmod}{\operatorname{mod}}
\newcommand{\ha}[1]{\frac{#1}{2}}
\newcommand{\lexp}{\operatorname{Exp}}
\newcommand{\Ah}{\hat{A}}
\newcommand{\Q}{\mathbb{Q}}
\newcommand{\dei}{\text{d}i}
\newcommand{\dej}{\text{d}j}
\newcommand{\lieh}{\mathfrak{h}}
\newcommand{\liek}{\mathfrak{k}}
\newcommand{\etatil}{\tilde{\eta}}
\newcommand{\rhotil}{\tilde{\rho}}
\newcommand{\Zh}{\hat{Z}}
\newcommand{\teh}{\hat{\theta}}
\newcommand{\muh}{\hat{\mu}}
\newcommand{\Ad}{\operatorname{Ad}}
\newcommand{\liez}{\mathfrak{z}}
\newcommand{\lspan}{\operatorname{span}}
\newcommand{\Xtil}{\tilde{X}}

\newcommand{\lieg}{\mathfrak{g}}
\newcommand{\Sq}{\operatorname{Sq}}
\newcommand{\Sgn}{\operatorname{Sgn}}

\newcommand{\Mhh}{\hat{M}
\hspace{-1.05ex}\hat{\rule{0ex}{2.0ex}}\hspace{1.05ex}}
\newcommand{\dehh}{\hat{\de}
\hspace{-.95ex}\hat{\rule{0ex}{2.05ex}}\hspace{.95ex}}
\newcommand{\vep}{\varepsilon}

\newcommand{\cA}{\mathcal{A}}
\newcommand{\Wtil}{\tilde{W}}

\newcommand{\kruisje}[1]{\, \mbox{$_{#1}$}\hspace{-.2ex}\mbox{$\ltimes$}
\,}
\newcommand{\strong}{\mbox{$\si$-strong$^*$}\xspace}

\newcommand{\alh}{\hat{\alpha}}
\newcommand{\Mh}{\hat{M}}

\newcommand{\cO}{\mathcal{O}}

\newcommand{\te}{\theta}

\newcommand{\Tr}{{\operatorname{Tr}}}
\newcommand{\nsf}{n.s.f.\xspace\!}

\newcommand{\deo}{\de_1}
\renewcommand{\det}{\de_2}
\newcommand{\deth}{\hat{\de}_2}

\newcommand{\deoh}{\hat{\de}_1}

\newcommand{\Wh}{\hat{W}}

\newcommand{\Nfih}{\cN_{\vfih}}
\newcommand{\sluit}{{- \; \strong}}
\newcommand{\lrecht}{\longrightarrow}

\newcommand{\Ytil}{\tilde{Y}}

\newcommand{\betil}{\tilde{\be}}

\newcommand{\vfitil}{\tilde{\vfi}}
\newcommand{\cR}{\mathcal{R}}

\newcommand{\Rh}{\hat{R}}
\newcommand{\Ma}{M_a}
\newcommand{\Mb}{M_b}
\newcommand{\dea}{\de_a}
\newcommand{\deb}{\de_b}
\newcommand{\ala}{\al_a}
\newcommand{\alb}{\al_b}
\newcommand{\bea}{\be_a}
\newcommand{\beb}{\be_b}
\newcommand{\Mah}{\Mh_a}
\newcommand{\Mbh}{\Mh_b}
\newcommand{\deah}{\deh_a}
\newcommand{\debh}{\deh_b}

\newcommand{\cUa}{\cU_a}
\newcommand{\cUb}{\cU_b}
\newcommand{\cVa}{\cV_a}
\newcommand{\cVb}{\cV_b}

\newcommand{\deh}{\hat{\Delta}}

\newcommand{\vfih}{\hat{\vfi}}

\newcommand{\lah}{\hat{\la}}

\newcommand{\pih}{\hat{\pi}}
\newcommand{\sih}{\hat{\si}}

\newcommand{\cU}{{\cal U}}

\newcommand{\Jh}{\hat{J}}

\newcommand{\cI}{{\cal I}}
\newcommand{\cL}{{\cal L}}

\newcommand{\cG}{{\cal G}}

\newcommand{\cN}{{\cal N}}

\newcommand{\cD}{{\cal D}}

\newcommand{\sla}{\lambda}
\newcommand{\la}{\Lambda}
\newcommand{\Om}{\Omega}
\newcommand{\om}{\omega}
\newcommand{\vfi}{\varphi}
\newcommand{\al}{\alpha}
\newcommand{\be}{\beta}
\newcommand{\ga}{\gamma}
\newcommand{\sde}{\delta}
\newcommand{\de}{\Delta}

\newcommand{\si}{\sigma}

\newcommand{\cV}{{\cal V}}
\newcommand{\cW}{{\cal W}}

\newcommand{\cst}{\text{C}$\hspace{0.1mm}^*$}

\newcommand{\deop}{\de \hspace{-.3ex}\raisebox{0.9ex}[0pt][0pt]
{\scriptsize\fontshape{n}\selectfont op}}

\newcommand{\deopt}{\de_2 \hspace{-1ex}\raisebox{1ex}[0pt][0pt]{\scriptsize\fontshape{n}
\selectfont op}}

\newcommand{\dehop}{\deh \hspace{-.3ex}\raisebox{0.9ex}[0pt][0pt]{\scriptsize\fontshape{n}\selectfont op}}
\newcommand{\dehopo}{\deh_1 \hspace{-1ex}\raisebox{1ex}[0pt][0pt]{\scriptsize\fontshape{n}\selectfont op}}
\newcommand{\dehopt}{\deh_2 \hspace{-1ex}\raisebox{1ex}[0pt][0pt]{\scriptsize\fontshape{n}\selectfont op}}

\newcommand{\recht}{\rightarrow}

\newcommand{\R}{\mathbb{R}}
\newcommand{\T}{\mathbb{T}}
\newcommand{\C}{\mathbb{C}}

\newcommand{\tekst}[1]{\quad\mbox{#1}\quad}
\newcommand{\ot}{\otimes}

\newcommand{\Z}{\mathbb{Z}}

\newcommand{\Dim}{\operatorname{dim}}
\newcommand{\Ker}{\operatorname{Ker}}

\newcommand{\Atil}{\tilde{A}}

\newcommand{\Htil}{\tilde{H}}

\newcommand{\io}{\iota}

\newcommand{\su}{\mathfrak{su}}
\newcommand{\sll}{\mathfrak{sl}}
\newcommand{\SU}{\operatorname{SU}}
\newcommand{\SL}{\operatorname{SL}}
\newcommand{\PSL}{\operatorname{PSL}}

\renewcommand{\Im}{\operatorname{Im}}
\renewcommand{\Re}{\operatorname{Re}}
\newcommand{\N}{\mathbb{N}}

\begin{document}

\title{On Low-Dimensional Locally Compact Quantum Groups}

\author{Stefaan Vaes\thanks{Research Assistant of the Fund for Scientific 
Research -- Flanders (Belgium)\ (F.W.O.)}\ \ and Leonid Vainerman}

\address{
Institut de Math{\'e}matiques de Jussieu, Alg{\`e}bres
d'Op{\'e}rateurs,
Plateau 7E \\ 175, rue du Chevaleret, F-75013 Paris, France \\
email:\,\texttt{vaes@math.jussieu.fr}
\\[4pt]
D\'epartement de Math\'ematiques et Mechanique, Universit\'e de Caen, \\
Campus II -- Boulevard de Mar\'echal Juin, \\ B.P. 5186, F-14032 Caen
Cedex, France \\
email:\,\texttt{lvainerman@yahoo.com}}

\markboth{Stefaan Vaes and Leonid Vainerman} {On Low-Dimensional
Locally Compact Quantum Groups}

\maketitle

\begin{abstract} 
Continuing our research on extensions of locally compact quantum groups, 
we give a classification of all cocycle matched pairs of Lie algebras in 
small dimensions and prove that all of them can be exponentiated to 
cocycle matched pairs of Lie groups. Hence, all of them give rise to 
locally compact quantum groups by the cocycle bicrossed product construction. 
We also clarify the notion of an extension of locally compact quantum groups 
by relating it to the concept of a closed normal quantum subgroup and the 
quotient construction. Finally, we describe the infinitesimal objects of
locally compact quantum quantum groups with 2 and 3 generators - Hopf 
$*$-algebras and Lie bialgebras.
\end{abstract}

\section{Introduction}

In this paper we continue the research on extensions of locally
compact (l.c.) quantum groups, initiated in \cite{VV}. The first
wide class of l.c.\ quantum groups, namely G.I. Kac algebras, was
introduced in the early sixties (see \cite{Kac1}) in order to
explain in a symmetric way duality for l.c.\ groups. This class
included besides usual l.c.\ groups and their duals also
nontrivial (i.e., non-commutative and non-cocommutative) objects
\cite{KP1}, \cite {KP2}.
 The general Kac algebra theory was completed independently on the one hand by
 G.I.~Kac and the second author \cite{KVai} and on the other hand by M.~Enock
 and J.-M.~Schwartz (for a survey see \cite{E-S}). However,
this theory was not general enough to cover important new examples
constructed starting from the eighties \cite{B}, \cite{KK},
\cite{KK2}, \cite{Maj2}, \cite{PW}, \cite{PuszWor}, \cite{VS},
\cite{VDW}, \cite{W1} - \cite{WZ}, which motivated essential efforts
to get a generalization that would cover these examples and that would 
be as elegant and symmetric as the theory of Kac algebras. Important 
steps in this direction were made by S.~Baaj and G.~Skandalis \cite{B-S1},
S.L.~Woronowicz \cite{WorComp}, \cite{W},
\cite{W1}, T.~Masuda and Y. Nakagami \cite{Mas-Nak} and A.~Van Daele 
\cite{VanDaele1}. The general theory of l.c.\ quantum groups was
proposed by J.~Kustermans and the first author \cite{KV1}, \cite{KV2} 
(see \cite{KV3} for an overview). Some motivations and applications of this 
theory can be found in the recent lecture notes \cite{KVVVW}.

The mentioned examples of l.c.\ quantum groups are, first of all,
formulated algebraically, in terms of generators of Hopf $^*$-algebras
and commutation relations between them. Then one represents the
generators as (typically, unbounded) operators on a Hilbert space and
tries to give a meaning to the commutation relations as relations between these 
operators. There is no general approach to this nontrivial problem, and one 
elaborates specific methods in each specific case.
Finally, it is necessary to associate an operator algebra with the
above system of operators and commutation relations and to construct 
comultiplication, antipode and invariant weights
as applications related to this algebra. This problem is even more difficult 
than the previous one and again one must consider separately each specific case 
(see the same papers). So, it would be desirable to have some general constructions 
of l.c.\ quantum groups which would allow to construct systematically concrete 
examples in a unified way.

One of such possibilities is offered by the {\it cocycle bicrossed
product construction}. According to G.I. Kac \cite{Kac}, in the
simplest case the needed data for this construction contains:
\begin{enumerate}
\item A pair of finite groups $G_1$ and $G_2$ equipped with
their mutual actions on each other (as on sets) or, equivalently, $G_1$ and $G_2$ 
must be subgroups of a certain group $G$ such that
$G_1\cap G_2=\{e\}$ and any $g\in G$ can be written as $g=g_1g_2\
(g_1\in G_1, g_2\in G_2)$ - we write briefly $G=G_1G_2$. We then say,
that $G_1$ and $G_2$ form {\it a matched pair of groups} \cite{Tak}.
\item A pair of compatible 2-cocycles for these
actions, so $G_1$ and $G_2$ must form {\it a cocycle matched pair} (in what
follows we often write simply "cocycle" rather then "2-cocycle").
\end{enumerate}
Then, due to \cite{Kac}, one can construct a finite-dimensional Kac algebra from 
cocycle crossed products of the algebras of functions on each of the groups $G_1$ 
and $G_2$ with the cocycle action of the other group, and this construction gives 
exactly all extensions of the above groups in the category of finite-dimensional 
Kac algebras.

It is tempting to similarly treat Lie groups instead of finite groups, being 
supported by the theory of cocycle bicrossed products and extensions of l.c.
\ groups developed in \cite{VV} (in fact, in \cite{VV}, the general theory of 
cocycle bicrossed products and extensions of l.c.\ {\em quantum} groups was 
developed). But first of all it turns out that the above definition of a 
matched pair of groups in terms of the equality $G=G_1G_2$ does not cover all 
interesting examples, see \cite{B-S2}. Following S.~Baaj and G.~Skandalis, one 
can just require $G_1G_2$ to be an open subset of $G$ with complement of measure 
zero. Then, the Lie algebras $\lieg_1$ of
$G_1$ and $\lieg_2$ of $G_2$ are Lie subalgebras of the Lie algebra
$\lieg$ of $G$ and $\lieg = \lieg_1 \oplus \lieg_2$ as the direct sum
of vector spaces, i.e., they form a matched pair of Lie algebras
\cite{Majbook}, 8.3. Thus, to get matched pairs of Lie groups one can
start with matched pairs of Lie algebras (which are easier to find)
and then try to exponentiate them.

To construct in this way cocycle matched pairs of Lie groups, one has to 
resolve two problems. First, given a matched pair of Lie algebras
$(\lieg_1,\lieg_2)$ with $\lieg = \lieg_1 \oplus \lieg_2$, one
can always exponentiate $\lieg$ to a connected and
simply connected Lie group $G$ and then
find Lie subgroups $G_1$ and $G_2$ whose Lie algebras are $\lieg_1$
and $\lieg_2$, respectively. However, such a choice of $G$ does not guarantee 
that $G_1 G_2$ is dense in $G$, even if $\dim(G_1)=\dim(G_2)=1$ \cite{Majbook}, 
\cite{Skand}, \cite{VV}, and it also may happen that $G_1\cap G_2 \neq \{e\}$. 
So, it is necessary to pass to some non-connected Lie group $G$ with the same 
Lie algebra $\lieg$ in order to find a matched pair of its subgroups $G_1$ and 
$G_2$ \cite{SV2}, \cite{VV}.
Secondly, given a matched pair of Lie groups, one has to find the corresponding 
cocycles. We give a solution of both these problems for {\em real} Lie groups
$G_1$ and $G_2$ with $\dim(G_1)=1,\ \dim(G_2)
\leq 2$ and construct essentially all possible (up to obvious redundancies) 
matched pairs of such Lie groups having at most 2 connected components. 
Then, using the machinery of cocycle bicrossed products developed in \cite{VV}, 
we construct l.c.\ quantum groups which are extensions of the mentioned Lie 
groups. Our discussion is motivated, apart from the above work by G.I. Kac, 
also by the works by S. Majid \cite{Maj1} - \cite{Majbook}, S. Baaj and 
G. Skandalis \cite{B-S1}, \cite{B-S2}, \cite{Skand}, and by the works on 
extensions of Hopf algebras \cite{A}, \cite{A-D}, \cite{Schn}.

The material is organized as follows. In Section 2, we recall the necessary 
facts of the theory of l.c.\ quantum groups and, following \cite{VV}, the main 
features of the cocycle bicrossed product construction
for l.c.\ groups in connection with the theory of extensions. In the last 
subsection we introduce the notion of a closed normal quantum subgroup of a
l.c.\ quantum group and explain its relation to the theory of extensions. 
As we explained above, the basic notion of
this theory is that of a matched pair of l.c.\ groups. If the groups forming a 
matched pair are Lie groups, we naturally have a matched pair of their Lie 
algebras. But the converse problem, to construct a matched pair of Lie groups
from a given matched pair $(\lieg_1,\lieg_2)$ of Lie algebras, is much more 
subtle. In particular, in Section 3 we show that any matched pair with 
$\lieg_1=\lieg_2=\C$ can be exponentiated to a matched pair of complex Lie 
groups, but there are simple examples of matched pairs of real and complex Lie 
algebras for which the exponentiation is impossible. 

The study of matched pairs of Lie algebras with $\dim \lieg_1 = n,\ 
\dim \lieg_2 = 1$ in Section~4 splits in three cases. In case 1, when 
$\lieg_1$ is an ideal in $\lieg,\ G$ can be constructed as semi-direct product 
of connected and simply connected Lie groups corresponding to $\lieg_1$ and 
$\lieg_2$ (this is possible also for $\dim \lieg_2 > 1$). In case 2, when 
$\lieg_1$ contains an ideal of codimension 1, the 
results of Section 3 show that for complex Lie algebras the exponentiation 
always exists when 
$n=1$ and it does not exist in general if $n\geq 2$. For real Lie algebras 
we show that for $n\leq 4$ there always exists the exponentiation to a matched 
pair of Lie groups with at most two connected components, and for $n\geq 5$ the
exponentiation does not exist in general. In the remaining case 3, for complex 
Lie algebras 
the exponentiation always exists when $n\leq 3$ and it does not exist in general 
if $n\geq 4$. For real Lie algebras we show that the exponentiation always 
exists when $n\leq 4$.

Section~5 is devoted to the complete classification of all matched pairs of 
real Lie algebras $\lieg_1$ and $\lieg_2$ when $\dim(\lieg_1)=1,\ \dim(\lieg_2)
\leq 2$ and to their explicit exponentiation to matched 
pairs of real Lie groups having at most 2 connected components. Here, we also
describe the l.c.\ quantum groups obtained from these matched pairs by the 
bicrossed product construction. In Section~6, we calculate the cocycles for all
the above mentioned matched pairs. Finally, Section~7 is devoted to the 
description of l.c.\ quantum groups with 2 and 3 generators and their 
infinitesimal objects - Hopf $*$-algebras and Lie bialgebras, having the 
structure of a cocycle bicrossed product - equivalently, those that can be 
obtained as extensions (we call them decomposable). At last, to complete the
picture of low-dimensional l.c.\ quantum groups, we review the indecomposable 
ones and their infinitesimal objects.

\subsubsection*{Acknowledgements} The first author would like to
thank the research group Analysis of the Department of Mathematics
of the K.U.Leuven for the nice working atmosphere while this work
was initiated. He also wants to thank the whole Operator Algebra
team of the Institut de Math{\'e}matiques de Jussieu in Paris for
their warm hospitality and the many useful discussions while
this work was finalized. The second author is grateful to the
research group Analysis of the Department of Mathematics of the K.U.
Leuven, to l'Institut de Recherche Math\'ematique Avanc\'ee de
Strasbourg and to Max-Planck-Insitut f\"ur Mathematik in Bonn for
the warm hospitality and financial support during his work on this
paper.

\section{Preliminaries}

\subsubsection*{General notations} Let $B(H)$ denote the algebra of
all bounded linear operators on a Hilbert space $H$, let $\ot$ denote
the tensor product of Hilbert spaces or von Neumann algebras and
$\Sigma$ (resp., $\sigma$) the flip map on it. If $H, K$ and $L$
are Hilbert spaces and $X \in B(H \ot L)$ (resp., $X \in B(H \ot
K), X \in B(K \ot L)$), we denote by $X_{13}$ (resp., $X_{12},\
X_{23}$) the operator $(1 \ot \Sigma^*)(X \ot 1)(1 \ot \Sigma)$
(resp., $X\ot 1,\ 1\ot X$) defined on $H \ot K \ot L$. Sometimes,
when $H = H_1 \ot H_2$ itself is a tensor product of two Hilbert
spaces, we switch from the above leg-numbering notation with
respect to $H \ot K \ot L$ to the one with respect to the finer
tensor product $H_1 \ot H_2\ot K \ot L$, for example, from
$X_{13}$ to $X_{124}$. There is no confusion here, because the
number of legs changes.

Given a comultiplication $\de$, denote by $\deop$ the opposite 
comultiplication $\sigma \de$. Our general reference to the modular 
theory of {\it normal semi-finite faithful} (n.s.f.) weights on von 
Neumann algebras is \cite{Stra}. For any  weight $\theta$ on a von 
Neumann algebra $N$, we use the notations
\begin{align*}
{\cal M}^+_\theta &= \{ x \in N^+ \mid \theta(x) < + \infty \},
\qquad {\cal N}_\theta = \{ x \in N \mid x^*x \in {\cal
M}^+_\theta \} \quad\text{and}\\ {\cal M}_\theta &=
\operatorname{span} {\cal M}^+_\theta \; .
\end{align*}

\subsubsection*{L.c.\ quantum groups}
A pair $(M,\de)$ is called a (von Neumann algebraic) l.c.\ quantum group \cite{KV2} when
\begin{itemize}
\item $M$ is a von Neumann algebra and $\de : M \recht M \ot M$ is
a normal and unital $*$-homomorphism satisfying the coassociativity relation : 
$(\de \ot \io)\de = (\io \ot \de)\de$.
\item There exist \nsf weights $\varphi$ and $\psi$ on $M$ such that
\begin{itemize}
\item $\varphi$ is left invariant in the sense that $\varphi \bigl( (\om \ot
\io)\de(x) \bigr) = \varphi(x) \om(1)$ for all $x \in {\cal M}_{\varphi}^+$ and $\om \in M_*^+$,
\item $\psi$ is right invariant in the sense that $\psi \bigl( (\io \ot
\om)\de(x) \bigr) = \psi(x) \om(1)$ for all $x \in {\cal M}_{\psi}^+$ and $\om \in M_*^+$.
\end{itemize}
\end{itemize}
Left and  right invariant weights are unique up
to a positive scalar \cite{KV1}, Theorem~7.14.

Represent $M$ on the Hilbert space of a GNS-construction $(H,\io,\Lambda)$ for 
the left invariant \nsf weight $\varphi$ and define a unitary $W$ on $H \ot H$ by
$$
W^* (\Lambda(a) \ot \Lambda(b)) = (\Lambda \ot \Lambda)(\de(b)(a \ot 1)) 
\quad\text{for all}\; a,b \in N_{\phi}\; .
$$
Here, $\Lambda \ot \Lambda$ denotes the canonical GNS-map for the tensor 
product weight $\varphi \ot \varphi$.
One proves that $W$ satisfies the pentagonal equation: $W_{12} W_{13} W_{23} = W_{23} W_{12}$, and we say that $W$ is a multiplicative unitary.
The von Neumann algebra $M$ and the comultiplication on it can be given in terms of $W$ respectively as
$$M = \{ (\io \ot \om)(W) \mid \om \in B(H)_* \}^\sluit \; $$
and $\de(x) = W^* (1 \ot x) W$, for all $x \in M$. Next, the l.c.\ quantum 
group $(M,\de)$ has an antipode $S$, which is the unique \strong closed
linear map from $M$ to $M$ satisfying $(\io \ot \om)(W) \in \cD(S)$ for all $\om \in B(H)_*$
and $S(\io \ot \om)(W) = (\io \ot \om)(W^*)$ and such that the elements $(\io \ot \om)(W)$
form a \strong core for $S$. $S$ has a polar decomposition $S = R \tau_{-i/2}$ where $R$ is an
anti-automorphism of $M$ and $(\tau_t)$ is a strongly continuous one-parameter group of
automorphisms of $M$. We call $R$ the unitary antipode and $(\tau_t)$ the scaling group of
$(M,\de)$. From \cite{KV1}, Proposition~5.26 we know that $\sigma (R \ot R) \de = \de R$. So
$\varphi R$ is a right invariant weight on $(M,\de)$ and we take $\psi:= \varphi R$.

Let us denote by $(\sigma_t)$ the modular automorphism group of $\varphi$. From \cite{KV1}, Proposition~6.8 we get the existence of a number $\nu > 0$, called the scaling constant, 
such that $\psi \, \sigma_t = \nu^{-t} \, \psi$ for all
$t \in \R$. Hence, we get the existence of a unique positive, self-adjoint operator $\sde_M$
affiliated to $M$, such that $\sigma_t(\sde_M) = \nu^t \, \sde_M$ for all $t \in \R$ and $\psi =
\varphi_{\sde_M}$, see \cite{KV1}, Definition~7.1. Formally this means that $\psi(x) =
\varphi(\sde_M^{1/2} x \sde_M^{1/2})$, and for a precise definition of $\varphi_{\sde_M}$ we refer
to \cite{SV3}. The operator $\sde_M$ is called the modular element of $(M,\de)$. If 
$\sde_M=1$ we call $(M,\de)$ unimodular. The scaling constant can be characterized as 
well by the relative invariance $\varphi \, \tau_t = \nu^{-t} \, \varphi$.

The dual l.c.\ quantum group $(\Mh,\deh)$ is defined in \cite{KV1}, Section~8. Its von 
Neumann algebra $\Mh$ is $$\Mh = \{(\om \ot \io)(W) \mid \om \in B(H)_* \}^\sluit$$ and the
comultiplication $\deh(x) = \Sigma W (x \ot 1) W^* \Sigma$ for all $x \in \Mh$. If we turn the
predual $M_*$ into a Banach algebra with product $\om \, \mu = (\om \ot \mu)\de$ and define
$$\lambda: M_* \recht \Mh : \lambda(\om) = (\om \ot \io)(W),$$ then $\lambda$ is a
homomorphism and $\lambda(M_*)$ is a \strong dense subalgebra of
$\Mh$. To construct explicitly a left invariant \nsf weight $\hat
\varphi$ with GNS-construction $(H,\io,\lah)$, first introduce the
space $$\cI = \{ \om \in M_* \mid \; \text{there exists} \;
\xi(\om) \in H \;\text{s.t.}\; \om(x^*) = \langle \xi(\om),
\Lambda(x) \rangle \;\text{when}\; x \in {\cal N}_{\varphi} \}.$$
If $\om \in \cI$, then such a vector $\xi(\om)$ clearly is
uniquely determined. Next, one proves that there exists a unique
\nsf weight $\hat\varphi$ on $\Mh$ with GNS-construction
$(H,\io,\lah)$ such that $\lambda(\cI)$ is a core
for $\lah$ (when we equip $\Mh$ with the \strong topology and $H$ with
the norm topology) and such that
$$\lah(\lambda(\om)) = \xi(\om) \tekst{for all} \om \in \cI \; .$$ One proves that the weight
$\hat\varphi$ is left invariant, and the associated multiplicative unitary is denoted by $\hat{W}$.
From \cite{KV1}, Proposition~8.16 it follows that $\hat{W} = \Sigma W^* \Sigma$.

Since $(\Mh,\deh)$ is again a l.c.\ quantum group, we can introduce the antipode $\hat{S}$,
the unitary antipode $\hat{R}$ and the scaling group $(\hat{\tau}_t)$ exactly as we did it 
for $(M,\de)$. Also, we can again construct the dual of $(\Mh,\deh)$, starting from the left
invariant weight $\vfih$ with GNS-construction $(H,\io,\lah)$. From \cite{KV1}, Theorem~8.29
we have that the bidual l.c.\ quantum group $(\Mhh,\dehh)$ is  isomorphic to $(M,\de)$.

We denote by $(\sih_t)$ the modular automorphism groups of the weight $\vfih$. The modular
conjugations of the weights $\varphi$ and $\hat\varphi$ will be denoted by $J$ and $\Jh$ respectively.
Then it is worthwhile to mention that $$R(x) = \Jh x^* \Jh \quad\text{for all} \; x \in M
\qquad\text{and}\qquad \Rh(y) = J y^* J \quad\text{for all}\; y \in \Mh \; .$$
Let us mention important special cases of l.c.\ quantum groups.

a) {\it Kac algebras} \cite{E-S}. From \cite{E-S}, we know that $(M,\de)$ is a Kac 
algebra if and only if $(\tau_t)$ is trivial
and $\sigma_t \, R = R \, \sigma_{-t}$ for all $t \in \R$. Now, denote by $(\sigma'_t)$ the modular
automorphism group of $\psi$. Because $\psi = \varphi R$ we get that $\sigma'_t \, R = R \,
\sigma_{-t}$ for all $t \in \R$. Hence $(M,\de)$ is a Kac algebra if and only if $(\tau_t)$ 
is trivial and $\sigma'=\sigma$. From \cite{SV3}, we know that $\sigma'_t(x) = \sde_M^{it} \sigma_t(x) \sde_M^{-it}$ for all $x \in M$ and $t \in \R$. Hence $\sigma'=\sigma$ if and 
only if $\sde_M$ is affiliated to the center of $M$.

In particular, $(M,\de)$ is a Kac algebra if $M$ is commutative.
Then $(M,\de)$ is generated by a usual l.c.\ group $G:\
M=L^{\infty}(G),\ (\de f)(g,h) = f(gh)$ \linebreak $(Sf)(g) =
f(g^{-1}),\ \varphi(f)=\int f(g)\; dg$, where $f\in
L^{\infty}(G),\ g,h\in G$ and we integrate with respect to the
left Haar measure $dg$ on $G$. The right invariant weight $\psi$
is given by $\psi(f) = \int f(g^{-1}) \; dg$. The modular element
$\sde_M$ is given by the strictly positive function $g \mapsto
\sde_G(g)^{-1}$.

The von Neumann algebra $M=L^\infty(G)$ acts on $H=L^2(G)$ by multiplication and $$(W_G\xi)(g,h)=\xi(g,g^{-1}h)$$ for all
$\xi\in H\ot H=L^2(G\times G)$. Then $\Mh=\cL(G)$ is the group von Neumann algebra generated
by the operators $(\lambda_g)_{g\in G}$ of the left regular representation of $G$ and
$\deh(\lambda_g)=\lambda_g\ot\lambda_g$. Clearly, $\dehop:=\sigma\deh=\deh$, so $\deh$ is
cocommutative.

b) A l.c.\ quantum group is called {\it compact} if its Haar measure is finite: $\varphi(1)<+\infty$,
which is equivalent to the fact that the norm closure of \linebreak $\{(\io\ot\om)(W)\vert\om\in
B(H)_*\}$ is a unital $C^*$-algebra. A l.c.\ quantum group $(M,\de)$ is called {\it discrete} if
$(\Mh,\deh)$ is compact.

\subsubsection*{Crossed and bicrossed products}
An {\it action} of a l.c.\ quantum group $(M,\de)$ on a von
Neumann algebra $N$ is a normal, injective and unital
$*$-homomorphism $\al: N \recht M \ot N$ such that $(\io \ot
\al)\al(x) =  (\de \ot \io)\al(x)$ for all $x \in N$. This
generalizes the definition of an action of a (separable) l.c.\
group $G$ on a ($\sigma$-finite) von Neumann algebra $N$, as a
continuous map $G \recht\operatorname{Aut} N : s \mapsto \al_s$
such that $\al_{st}=\al_s \al_t$ for all $s,t\in G$. Indeed,
putting $M=L^\infty(G)$, one can identify $M \ot N$ with
$L^\infty(G,N)$ and $M \ot M \ot N$ with $L^\infty(G \times G,N)$
and define the above homomorphism $\al$ by $(\al(x))(s) =
\al_{s^{-1}}(x)$. The fixed point algebra of an action $\al$ is
defined by $N^\alpha=\{x\in N\mid \al(x)=1\ot x\}$.

A {\it cocycle} for an action of a l.c.\ group $G$ on a
commutative von Neumann algebra $N$ is a Borel map $u:G \times G
\recht N$ such that $\al_r(u(s,t)) \; u(r,st) = u(r,s) \; u(rs,t)$
nearly everywhere. Then, putting $M=L^\infty(G)$, one can define a
unitary $\cU \in M \ot M \ot N$ by $\cU(s,t)=u(t^{-1},s^{-1})$
satisfying $$(\io \ot \io \ot \al)(\cU) (\de \ot \io \ot \io)(\cU)
= (1 \ot \cU) (\io \ot \de \ot \io)(\cU) \; .$$ For the general
definition of a cocycle action of a l.c.\ quantum group on an
arbitrary von Neumann algebra, we refer to Definition~1.1 in
\cite{VV}.

The {\it cocycle crossed product} $G \kruisje{\al,\cU} N$ is the von Neumann subalgebra of $B(L^2(G)) \ot N$
generated by $$\al(N) \tekst{and} \{ (\om \ot \io \ot \io)(\Wtil) \mid \om \in L^1(G) \} \; ,$$
where $\Wtil = (W_G \ot 1)\cU^*$. This is a von Neumann algebraic version of the twisted $C^*$-algebraic crossed product \cite{Pac-Rae}.
There exists a unique action $\alh$ of $(\cL(G),\deh)$ on $G \kruisje{\al,\cU} N$ such that
\begin{align*}
\alh(\al(x)) & = 1 \ot \al(x) \tekst{for all} x \in N, \\
(\io \ot \alh)(\Wtil) &= W_{G,12} \Wtil_{134} \; ,
\end{align*}
and for any \nsf weight $\theta$ on $N$, we can
define the {\it dual} \nsf weight $\tilde\theta$ on $G \kruisje{\al,\cU} N$ by the formula
$$\tilde\theta = \theta \al^{-1} \; (\hat\varphi \ot \io \ot \io) \alh.$$

\begin{definition}\label{mpg} (see \cite{B-S-V}) Let $G,G_1$ and $G_2$
be (separable) l.c.\ groups and let 
a homomorphism $i: G_1 \recht G$ and an anti-homomorphism $j:G_2
\recht G$ have closed images and be homeomorphisms onto these images.
Suppose that $i(G_1) \cap j(G_2) = \{e\}$ and that the complement of
$i(G_1)j(G_2)$ in $G$ has measure zero.
Then we call $G_1$ and $G_2$ a
matched pair of l.c.\ groups.
\end{definition}

Observe that this definition of a matched pair of l.c.\ groups, due to
Baaj, Skandalis and the first author, is more
general than the one studied in \cite{B-S2} and \cite{VV}. Indeed, in
\cite{B-S-V}, there is given an example of a matched pair in the sense of the
definition above, which does not fit in the definition of
\cite{B-S2}. More specifically, consider the map $$\te:G_1 \times G_2
\recht G: (g,s) \mapsto i(g) j(s) \; ,$$ which is clearly injective. In \cite{B-S2} and \cite{VV},
the map $\te$ is supposed to have a range $\Om$ which is open in $G$,
with complement of measure zero and such that $\te$ is a homeomorphism
of $G_1 \times G_2$ onto $\Om$. In the example of \cite{B-S-V}, the
range of $\te$ has an empty interior. However, the following proposition holds:

\begin{proposition} \label{LieOK}
If, in Definition~\ref{mpg}, $G$ is a Lie group, then the map $$\te:G_1 \times G_2
\recht G: (g,s) \mapsto i(g) j(s)$$ has an open range $\Om$ and is a
diffeomorphism of $G_1 \times G_2$ onto $\Om$, where $G_1$ and $G_2$
are Lie groups under the identification with closed subgroups of $G$.
\end{proposition}
\begin{proof}
Denote by $\lieg,\lieg_1,\lieg_2$ the Lie algebras of $G,G_1,G_2$,
respectively. Then, we have an injective homomorphism and
anti-homomorphism $$\dei : \lieg_1 \recht \lieg \quad\text{and}\quad \dej
: \lieg_2 \recht \lieg \; .$$ Because $i(G_1) \cap j(G_2) = \{e\}$, we
get $\dei(\lieg_1) \cap \dej(\lieg_2) = \{0\}$ (otherwise, the
exponential mapping produces elements in $i(G_1) \cap j(G_2)$). Hence,
we can take a linear subspace $\liek$ of $\lieg$ (not necessarily a
Lie subalgebra) such that $\lieg = \dei(\lieg_1) \oplus \dej(\lieg_2)
\oplus \liek$ as vector spaces. We first prove that $\liek=\{0\}$.

Denote by $\exp_{\lieg}$ the exponential mapping of $G$ and
analogously for $\exp_{\lieg_{1,2}}$. Take open subsets $\cU_i \subset
\lieg_i$, $\cV \subset \liek$ containing $0$ such that $\exp_{\lieg}$
is a diffeomorphism of $\dei(\cU_1) \times \dej(\cU_2) \times \cV$ onto an open
subset of $G$. Define
\begin{align*}
\rho: \cU_1 \times \cU_2 \times \cV \recht G :
\rho(v,w,z) &= i(\exp_{\lieg_1}(v)) \; j(\exp_{\lieg_2}(w)) \;
\exp_{\lieg}(z) \\&= \exp_{\lieg}(\dei(v)) \; \exp_{\lieg}(\dej(w)) \;
\exp_{\lieg}(z) \; .
\end{align*}
Because $\lieg = \dei(\lieg_1) \oplus \dej(\lieg_2)
\oplus \liek$, we find that $\text{d}\rho(0,0,0)$ is bijective. So,
for $\cU_1,\cU_2,\cV$ small enough, $\rho$ is a diffeomorphism onto an
open subset of $\cW$ of $G$ containing $e$ and $\exp_{\lieg_i}$ will
be a diffeomorphism of $\cU_i$ onto an open subset $\cW_i$ of $G_i$.

It is clear that $\te(\cW_1 \times \cW_2) \subset \cW$ and
$\rho^{-1}(\te(\cW_1 \times \cW_2)) = \cU_1 \times \cU_2 \times
\{0\}$. As a diffeomorphism, $\rho$ is a Borel isomorphism and so, if
$\liek \neq \{0\}$, $\te(\cW_1 \times \cW_2)$ has measure zero in
$G$. This contradicts the result of \cite{B-S-V}, saying that $\te$ is
automatically a Borel isomorphism. Hence, $\liek=\{0\}$.

But then, $\te : \cW_1 \times \cW_2 \recht \cW$ is a
diffeomorphism. In particular, $\cW \subset \Om$. If now $i(g_0)j(s_0)
\in \Om$, it follows that $i(g_0)j(s_0)
\in i(g_0) \cW j(s_0) \subset \Om$. Hence, $\Om$ is open in $G$ and
$\te$ is a diffeomorphism of $G_1 \times G_2$ onto $\Om$, because we
know that $\te$ is injective.
\end{proof}

In \cite{B-S-V}, it is proved that $\te$ is automatically a Borel isomorphism, 
i.e. it induces an isomorphism between $L^\infty(G_1 \times G_2)$ and
$L^\infty(G)$. Hence, 
this data allows to construct as follows two actions: $\al$ of $G_1$ on $M_2=L^\infty(G_2)$
and $\be$ of $G_2^[\text{\rm op}]$ on $M_1=L^\infty(G_1)$ verifying certain compatibility relations.

Define $\Om$ to be the image of $\te$ and define the Borel isomorphism
$$
\rho : G_1 \times G_2 \recht \Om^{-1} : (g,s) \mapsto j(s)i(g) \; .
$$
So $\cO=\theta^{-1}(\Om \cap \Om^{-1})$ and $\cO'=\rho^{-1}(\Om \cap
\Om^{-1})$ are Borel subsets of $G_1 \times G_2$, with complement of
measure zero, and $\rho^{-1} \theta$ is a Borel isomorphism
of $\cO$ onto $\cO'$. For all $(g,s) \in \cO$ define $\be_s(g) \in G_1$ 
and $\al_g(s) \in G_2$ such that
$$
\rho^{-1}(\theta(g,s)) = (\be_s(g),\al_g(s)) \; .
$$
Hence we get $j \bigl( \al_g(s) \bigr) \; i\bigl( \be_s(g) \bigr) = 
i(g)j(s)$ for all $(g,s) \in \cO$.

\begin{lemma}\label{42} (\cite{VV}, Lemma 4.8)

Let $(g,s) \in \cO$ and $h \in G_1$. Then $(hg,s) \in \cO$ if and only if $(h,\al_g(s)) \in
\cO$, and in that case
$$\al_{hg}(s) = \al_h \bigl( \al_g(s) \bigr) \qquad\text{and}\qquad \be_s(hg) =
\be_{\al_g(s)}(h) \; \be_s(g) \; .$$
Let $(g,s) \in \cO$ and $t \in G_2$. Then $(g,ts) \in \cO$ if and only if $(\be_s(g),t) \in
\cO$ and in that case
$$\be_{ts}(g) = \be_t \bigl( \be_s(g) \bigr) \qquad\text{and}\qquad \al_g(ts) =
\al_{\be_s(g)}(t) \; \al_g(s) \; .$$
Finally, for all $g \in G_1$ and $s \in G_2$ we have $(g,e) \in \cO$, $(e,s) \in \cO$, and
$$\al_g(e) = e, \quad \al_e(s) = s, \quad \be_s(e)=e \tekst{and} \be_e(g)=g \; .$$
\end{lemma}
This can be viewed as a definition of a matched pair of l.c.\ groups in terms of mutual actions.

The cocycles for the above actions can be introduced as measurable maps
$\cU  : G_1 \times G_1 \times G_2 \recht U(1)$ and
$\cV  : G_1 \times G_2 \times G_2 \recht U(1)$, where
$U(1)$ is the unit circle in $\C$, satisfying
\begin{align}
\cU(g,h,\al_k(s)) \; \cU(gh, k, s) &= \cU(h,k,s) \; \cU(g, hk, s), \notag\\
\cV(\be_s(g),t,r) \; \cV(g,s, rt) &= \cV(g,s,t) \; \cV(g, ts, r) \label{cocU} ,\\
\cV(gh, s, t) \; \bar{\cU}(g,h,ts) & = \bar{\cU}(g,h,s) \;
\bar{\cU}(\be_{\al_h(s)}(g), \be_s(h) , t) \notag \\ &
\qquad\qquad\qquad \cV(g, \al_h(s), \al_{\be_s(h)}(t)) \;
\cV(h,s,t) \notag
\end{align}
nearly everywhere. Then we have a definition of a cocycle matched pair of l.c.\ groups.

Fixing a cocycle matched pair of l.c.\ groups $G_1$ and $G_2$, denoting $H_i=L^2(G_i)$ $(i=1,2)$, $H=H_1\ot H_2$ and identifying $\cU$ and $\cV$ with unitaries in $M_1\ot M_1\ot M_2$ and
in $M_1\ot M_2\ot M_2$ respectively, define unitaries $W$ and $\Wh$ on $H \ot H$ by
\begin{align*}
\Wh &= (\be \ot \io \ot \io) \bigl( (W_{G_1} \ot 1) \cU^* \bigr)
\; (\io \ot \io \ot \al) \bigl( \cV (1 \ot \hat W_{G_2})
\bigr)\;\text{and}\; W = \Sigma \Wh^* \Sigma \; .
\end{align*}
On the von Neumann algebra $M = G_1 \kruisje{\al,\cU} L^{\infty}(G_2)$, let us define a faithful $*$-homomorphism
$$\de : M \recht B(H \ot H) : \de(z) = W^*(1 \ot z) W \ (\forall z\in M)$$
and denote by $\varphi$ the dual weight of the canonical left invariant trace $\varphi_2$ on $L^{\infty}(G_2)$.
Then, Theorem 2.13 of \cite{VV} shows that $(M,\de)$ is a l.c.\
quantum group with $\varphi$ as a left invariant weight, which we call the {\em cocycle bicrossed product} of $G_1$ and $G_2$. One can also show that its scaling constant is 1.
The dual l.c.\ quantum group is $(\Mh,\deh)$, where $\Mh= G_2 \kruisje{\al,\cU} L^{\infty}(G_1)$and $\deh(z) = \Wh^* (1 \ot z) \Wh$ for all $z \in \Mh$.

One can get explicit formulas for the modular operators, modular
conjugations of the left invariant weights, unitary antipodes, scaling
groups and modular elements of both $(M,\de)$ and its dual in terms of
the above mutual actions, the cocycles and the modular functions
$\delta,\delta_1$ and $\delta_2$ of the l.c.\ groups $G,G_1$ and $G_2$. In particular, one can characterize all cocycle bicrossed products of l.c.\ groups which are Kac algebras.
\begin{proposition} \label{charKac}
The l.c.\ quantum group $(M,\de)$ is a Kac algebra if and only if
\begin{align*}
& \sde\bigl( i(g \be_s(g)^{-1}) \bigr) \; \sde_1\bigl(g^{-1}
\be_s(g) \bigr) \; \sde_2\bigl( \al_g(s) s^{-1} \bigr) \ = 1
\tekst{and}\\ & \frac{\sde_1 \bigl( \be_s(g) \bigr)}{\sde_1(g)} =
\frac{\sde_2 \bigl( \al_g(s) \bigr)}{\sde_2(s)} \; .\end{align*}
\end{proposition}
This proposition implies three helpful corollaries.
\begin{corollary}\label{410}
If $\al$ or $\be$ is trivial, $(M,\de)$ and
$(\Mh,\deh)$ are Kac algebras.
\end{corollary}
\begin{corollary}\label{411}
If both $\al$ and $\be$ preserve modular functions and Haar measures, then $(M,\de)$ and $(\Mh,\deh)$ are Kac algebras.
\end{corollary}
Remark that the conditions of this corollary are fulfilled if both
groups are discrete. Indeed, any discrete group is unimodular and the
Haar measure is constant at an arbitrary point of such a group.
\begin{corollary}\label{hom}
If $(G_1,G_2)$ is a fixed matched pair of l.c.\ groups and cocycles $\cU$ and $\cV$ satisfy (\ref{cocU}), we get a cocycle
bicrossed product $(M,\de)$.
If one of these cocycle bicrossed products is a Kac algebra, then all of them are Kac algebras.
\end{corollary}
\begin{proof}
The necessary and sufficient conditions for $(M,\de)$ to be a Kac
algebra in Proposition~\ref{charKac} are independent of $\cU$
and $\cV$.
\end{proof}

It is easy to check that the above measurable mutual actions $\al_g$ and $\be_s$ of $G_1$ and $G_2$ are in fact the restrictions of the
canonical continuous actions $\tilde\al_g$ of $G_1$ on
$G/G_1$ and $\tilde\be_s$ of $G_2$ on $G_2 \backslash G$ (topologies on $G_1$ and $G_2 \backslash G$ and, respectively, on $G_2$ and
$G/G_1$, are in general different). This allows, in particular, to express the $C^*$-algebras of the $C^*$-algebraic
versions of the split extension (i.e.\ with trivial cocycles) and its
dual respectively as $G_1 \kruisje{\tilde{\al}} C_0(G/G_1)$ and
$C_0(G_2 \backslash G) \rtimes_{\tilde{\be}} G_2$, see
 \cite{B-S-V}.

\subsubsection*{Extensions of l.c.\ groups}
To clarify the following definition, recall that any normal
$*$-homo\-mor\-phism $\be : M_1 \recht \Mh$ of l.c.\ quantum
groups satisfying $\deh \be = (\be \ot \be) \deo$ generates two
canonical actions: $\mu$ of $(\Mh_1,\deoh)$ on $M$ and $\theta$ of
$(\Mh_1,\dehopo)$ on $M$ (\cite{VV}, Proposition 3.1). On a formal
level, this can be understood easily: the morphism $\be$ gives
rise to a dual morphism $\betil: M \recht \Mh_1$ and $\mu$ should
be thought of as $\mu = (\betil \ot \io)\de$, while $\te$ should
be thought of as $\te = (\betil \ot \io)\deop$.
\begin{definition} \label{32}
Let $G_i\ (i=1,2)$ be l.c.\ groups and let $(M,\de)$ be a l.c.\ quantum
group. We call $$(L^\infty(G_2),\det) \overset{\al}{\lrecht}
(M,\de) \overset{\be}{\lrecht} (\cL(G_1),\deoh)$$ a short exact
sequence, if $$\al: L^\infty(G_2) \recht M \quad\text{and}\quad
\be: L^\infty(G_1) \recht \Mh$$ are normal, faithful
$*$-homomorphisms satisfying $$\de \al = (\al \ot \al) \det
\quad\text{and}\quad \deh \be = (\be \ot \be) \deo$$ and if
$\al(L^\infty(G_2))  = M^\theta$, where $\theta$ is the canonical
action of $(\cL(G_1),\dehopo)$ on $M$  generated by the morphism
$\be$. In this situation, we call $(M,\de)$ an extension of $G_2$
by $\hat{G}_1$.
\end{definition}
The faithfulness of the morphisms $\al$ and $\be$ reflects the exactness of the sequence in
the first and third place. The formula $\al(L^\infty(G_2))=M^\theta$ reflects its exactness in the second
place.
Given a short exact sequence as above, one can check that the dual sequence
$$(L^\infty(G_1),\deo)
\overset{\be}{\lrecht} (\Mh,\deh) \overset{\al}{\lrecht} (\cL(G_2),\deth)$$
is exact as well.

Given a cocycle matched pair of l.c.\ groups, one can check that
their cocycle bicrossed product is an extension in the sense of
Definition~\ref{32}. Moreover, it belongs to a special class of
extensions, called cleft extensions (\cite{VV}, Theorem~2.8). This
theorem also shows that, conversely, all cleft extensions of l.c.\
groups (and of l.c.\ quantum groups) are given by the cocycle
bicrossed products. This means that, whenever $(M,\de)$ is a cleft
extension of $G_2$ by $\hat{G}_1$, the pair consisting of
$(L^\infty(G_1),\deo)$ and $(L^\infty(G_2),\det)$ is a cocycle
matched pair in the sense of \cite{VV}, Definition~2.1 and
$(M,\de)$ is isomorphic to their cocycle bicrossed product. From the
results of \cite{B-S-V}, it follows that this precisely means that
$(G_1,G_2)$ is a matched pair in the sense of Definition~\ref{mpg}
with cocycles as in Equation~\eqref{cocU}.

By definition, two extensions
\begin{align*}
& (L^\infty(G_2),\det) \overset{\al_a}{\lrecht} (M_a,\de_a)
\overset{\be_a}{\lrecht} (\cL(G_1),\deoh)\qquad \text{and}\\
& (L^\infty(G_2),\det) \overset{\al_b}{\lrecht} (M_b,\de_b)
\overset{\be_b}{\lrecht} (\cL(G_1),\deoh)
\end{align*}
are called isomorphic, if there is an isomorphism $\pi :
(\Ma,\dea) \recht (\Mb,\deb)$ of l.c.\ quantum groups satisfying
$\pi \ala = \alb$ and $\pih \bea = \beb$, where $\pih$ is the
canonical isomorphism of $(\Mah,\deah)$ onto $(\Mbh,\debh)$
associated with $\pi$.

Given a matched pair $(G_1,G_2)$ of l.c.\ groups, any couple of cocycles $(\cU,\cV)$ satisfying (\ref{cocU}) generates as above a cleft extension $$(L^\infty(G_2),\det) \overset{\al}{\lrecht} (M,\de)
\overset{\be}{\lrecht} (\cL(G_1),\deoh).$$  The extensions given by two pairs of cocycles $(\cUa,\cVa)$ and $(\cUb,\cVb)$, are isomorphic if and only if there exists a measurable map $\cR$ from
$G_1 \times G_2$ to $U(1)$, satisfying
\begin{align*}
\cUb(g,h,s) &= \cUa(g,h,s) \; \cR(h,s) \; \cR(g,\al_h(s)) \; \bar{\cR}(gh,s) \\
\cVb(g,s,t) &= \cVa(g,s,t) \; \cR(g,s) \; \cR(\be_s(g),t) \; \bar{\cR}(g,ts)
\end{align*}
almost everywhere. If this is the case, the pairs $(\cUa,\cVa)$ and $(\cUb,\cVb)$ will be
called cohomologous.
Then the set of equivalence classes of cohomologous pairs of cocycles $(\cU,\cV)$ satisfying (2.2), exactly corresponds to the set $\Gamma$ of classes of isomorphic extensions associated with $(G_1,G_2)$.

The set $\Gamma$ can be given the structure of an abelian group by defining
$$\pi(\cUa,\cVa) \; \cdot \; \pi(\cUb,\cVb) = \pi(\cUa \cUb,\cVa \cVb)$$
where $\pi(\cU,\cV)$ denotes the equivalence class containing the pair $(\cU,\cV)$.
The group $\Gamma$ is called the group of extensions of $(L^\infty(G_2),\det)$ by
$(\cL(G_1),\deoh)$ associated with the matched pair of l.c.\ groups $(G_1,G_2)$. The unit of
this group corresponds to the class of cocycles cohomologous to trivial. The
corresponding extension is called {\it split extension}; all other
extensions are called {\it non-trivial
extensions}.

\subsubsection*{Closed normal quantum subgroups}

Definition~\ref{32} is the partial case of the general definition of a short 
exact sequence
$$(M_2,\det) \overset{\al}{\lrecht} (M,\de) \overset{\be}{\lrecht} (\Mh_1,\deoh),$$
where $(M_1,\Delta_1),\ (M_2,\det)$ and $(M,\de)$ are l.c.\ quantum groups, see 
Definition~3.2 in \cite{VV}.
We explain the relation between this notion and the following notion of a 
closed normal quantum subgroup. 

\begin{definition} \label{closed}
A l.c.\ quantum group $(M_2,\det)$ is called a closed quantum
subgroup of $(M,\de)$ if there exists a normal, faithful
$^*$-homomorphism \linebreak $\al : M_2 \recht M$ such that $\de \al = (\al \ot
\al)\det$.
\end{definition}

This definition might need some justification: in \cite{Kus},
J.~Kustermans defines morphisms between l.c.\ quantum groups on the
(natural) level of universal C$^*$-algebraic quantum groups. So, it
might seem strange to require the existence of a normal morphism on
the von Neumann algebra level. We
claim, however, that this precisely characterizes the closedness (or
properness of the injective embedding). Let
us illustrate this with an example. Consider the identity map from
$\R_d$ with the discrete topology to $\R$ with its usual
topology. Dualizing, we get a morphism $\al : C_0(\R) \recht
M(C_0(\R_d)) = C_b(\R_d)$ which is injective. It is clear that we want
to exclude this type of morphisms. This is precisely achieved by
requiring the normality (weak continuity) of the morphism. To
conclude, we mention that in the case where $M_2 = \cL(G_2)$ and
$M=\cL(G)$, we precisely are in the situation of an identification
$\pi : G_2 \recht G$ of $G_2$ with a closed subgroup of $G$ and
$\al(\lambda_g)=\lambda_{\pi(g)}$, see Theorem~6 in \cite{TT}.

Next, we define normality of a closed quantum subgroup. Recall that
when $\cA_1$ is a Hopf subalgebra of a Hopf algebra $\cA$, $\cA_1$ is
called normal if $\cA_1$ is invariant under the adjoint action. Using
Sweedler notation, this means
$$\sum a_{(1)} x S(a_{(2)}) \in \cA_1 \quad\text{for all}\;\; x \in
\cA_1, a \in \cA \; .$$
Recalling that $S((\io \ot \om)(W)) = (\io \ot \om)(W^*)$ and that
$$\de((\io \ot \om)(W)) = (\io \ot \io \ot \om)(W_{13} W_{23})$$ because
of the pentagon equation, it is easy to verify that the operator
algebraic version of normality is given as follows.

\begin{definition} \label{normal}
If $\al: M_2 \recht M$ turns $(M_2,\det)$ into a closed quantum
subgroup of the l.c.\ quantum group $(M,\de)$, we say that
$(M_2,\det)$ is normal if $$W (\al(M_2) \ot 1) W^* \subset \al(M_2)
\ot B(H) \; .$$
\end{definition}

As could be expected, we now prove the bijective correspondence
between closed normal quantum subgroups and short exact sequences.

\begin{theorem}
Suppose that $\al: M_2 \recht M$ turns $(M_2,\det)$ into a closed
normal quantum
subgroup of $(M,\de)$. Then, there exists a
unique (up to isomorphism) l.c.\ quantum group $(M_1,\de_1)$ and a unique
$\be : M_1 \recht \Mh$ such that
$$(M_2,\det) \overset{\al}{\lrecht} (M,\de) \overset{\be}{\lrecht}
(\Mh_1,\deoh)$$ is a short exact sequence.

If, conversely, we have a short exact sequence, then $\al: M_2 \recht M$
turns $(M_2,\det)$ into a closed
normal quantum
subgroup of the l.c.\ quantum group $(M,\de)$.
\end{theorem}
\begin{proof}
Suppose first that we have a short exact sequence. Consider the
coaction $\te$ of $(\Mh_1,\dehopo)$ on $M$ associated with $\be$. By
definition of exactness, we have $\al(M_2)=M^\te$. Let $x \in M_2$. It
suffices to prove that $(\te \ot \io)(W(\al(x) \ot 1)W^*) = W_{23} (1 \ot \al(x)
\ot 1) W_{23}^*$. From Proposition~3.1 of \cite{VV} and with the
notations introduced over there, it follows that
it is sufficient to prove that $1 \ot \al(x)$ commutes with $Z_1$, or
equivalently, $\mu(\al(x))=1 \ot \al(x)$. But, $$\mu(\al(x)) = (\Rh_1
\ot R)\te(R(\al(x))) = (\Rh_1 \ot R)\te(\al(R_2(x))) = 1 \ot \al(x) \;
.$$
This proves the most easy, second part of the theorem.

Next, suppose that we have a closed normal quantum subgroup
$(M_2,\det)$ of $(M,\de)$. Using Proposition~3.1 from \cite{VV}, the
morphism $\al$ generates
two actions: $\muh$ is an action of $(\Mh_2,\deh_2)$ on
$\Mh$ and $\teh$ is an action of $(\Mh_2,\dehopt)$ on $\Mh$ and
they are determined by
$$\muh(x) = \Zh_1 (1 \ot x) \Zh_1^* \quad\text{and}\quad \teh(x)
= \Zh_2 (1 \ot x) \Zh_2^* \quad\text{for all}\;\; x \in \Mh \; ,$$
where $$\Zh_1 = (\io \ot \al)(\Wh_2^*) \quad\text{and}\quad \Zh_2
= (J_2 \ot J) \Zh_1 (J_2 \ot J) \; .$$ The actions $\muh$ and
$\teh$ are related by the formula $\teh(x) = (\Rh_2 \ot
\Rh)\muh(\Rh(x))$ and satisfy
$$(\muh \ot \io)(\Wh) = (\io \ot \al)(\Wh_2)_{13} \Wh_{23}
\quad\text{and}\quad (\teh \ot \io)(\Wh) = \Wh_{23} (\io \ot
\al)(\Wh_2)_{13} \; .$$ Using the definition of the left invariant
weight $\vfih$ on the dual $(\Mh,\deh)$ of $(M,\de)$, we easily
conclude that $\vfih$ is invariant under the action $\muh$ and
moreover, for all $x \in \Nfih$ and $\om \in \Mh_{2,*}$, we have
$(\om \ot \io)\muh(x) \in \Nfih$ and
$$\lah((\om \ot \io)\muh(x)) = (\om \ot \io)(\Zh_1) \lah(x) \; .$$
From Proposition~4.3 of \cite{SV1}, it then follows that $\Zh_1$
is the canonical implementation of the action $\muh$ (in the sense
of Definition~3.6 of \cite{SV1}). We want to prove that $\muh$ is
integrable (see Definition~1.4 in \cite{SV1}) and we will use
Theorem~5.3 of \cite{SV1} to do this. So, we have to construct a
normal $^*$-homomorphism $\rho: \Mh_2 \kruisje{\muh} \Mh \recht
B(H)$ such that $$\rho(\muh(x)) = x \quad\text{for all}\; x \in \Mh
\quad\text{and}\quad (\io \ot \rho)(\Wh_2 \ot 1) = \Zh_1^* \; .$$
We first define
$$\rhotil : \Mh_2 \kruisje{\muh} \Mh \recht B(H \ot H) :
\rhotil(z) = \cV (\al \ot \io)(\Zh_1^* z \Zh) \cV^* \; ,$$ where
$\cV = (\Jh \ot \Jh) W (\Jh \ot \Jh)$ has the properties $\cV \in
M \ot \Mh'$ and $\deop(y) = \cV^* (1 \ot y) \cV$ for all $y \in
M$. This map $\rhotil$ is well-defined for the following reasons.
For $x \in \Mh$, we have $\Zh_1^* \muh(x) \Zh_1 = 1 \ot x$. We can
apply $\al \ot \io$ and because $\cV \in M \ot \Mh'$, we find that
$\rhotil(\muh(x)) = 1 \ot x$. Next, for $\om \in \Mh_{2,*}$, we
find
\begin{align*}
\Zh_1^* ((\om \ot \io)(\Wh_2) \ot 1) \Zh_1 &= (\om \ot \io \ot
\io)(\io \ot \io \ot \al)(\Wh_{2,23} \Wh_{2,12} \Wh^*_{2,23}) \\
&= (\om \ot \io \ot \io)(\Wh_{2,12} (\io \ot \al)(\Wh_2)_{13}) \;
.
\end{align*}
Again, it is possible to apply $\al \ot \io$ and we find
$$(\io \ot \rhotil)(\Wh_2 \ot 1) = \cV_{23} (\io \ot \al \ot
\al)(\io \ot \deopt)(\Wh_2) \cV_{23}^* = \Zh^*_{1,13} \; ,$$
because $\al$ is a morphism and $\cV^*$ implements $\deop$. The
$\rho$ that we were looking for, is then obtained as $\rhotil(z) =
1 \ot \rho(z)$ for all $z \in \Mh_2 \kruisje{\muh} \Mh$. Hence,
$\muh$ is integrable.

Define the von Neumann algebra $M_1:= \Mh^{\muh}$, the fixed point
algebra of $\muh$. Because $(\io \ot \deh)\muh = (\muh \ot
\io)\deh$, it is clear that $\deh(M_1) \subset M_1 \ot \Mh$. We
claim that also $\deh(M_1) \subset \Mh \ot M_1$. For this, we will
need the normality. Observe that the right leg of $\Zh_1$
generates $\al(M_2)$. Hence, by definition, $M_1 = \Mh \cap
\al(M_2)'$. Because $R \al = \al R_2$ and $R(x) = \Jh x^* \Jh$, we
conclude that $\Jh M_1' \Jh =$ \linebreak $(\Mh \cup \al(M_2))^{\prime\prime}$.
By normality, we know that $W(\al(M_2) \ot 1) W^* \subset \al(M_2)
\ot \Mh$. Because $W(\Mh \ot 1) W^* = \dehop(\Mh) \subset \Mh \ot
\Mh$, we get $W(\Jh M_1' \Jh \ot 1) W^* \subset \Jh M_1' \Jh \ot
\Mh$. Writing $\Jh \ot J$ around this equation, we find $W^* (M_1'
\ot 1) W \subset$ \linebreak $M_1' \ot B(H)$. Because $W \in M \ot \Mh$, it
follows that $W^* (M_1' \ot \Mh') W \subset$ \linebreak $M_1' \ot B(H)$. Taking
commutants, we conclude that $M_1 \ot 1 \subset W^*(M_1 \ot
\Mh)W$. Bringing the $W$ to the other side, we have proven our
claim that $\dehop(M_1) \subset M_1 \ot \Mh$.

Defining $\de_1$ to be the restriction of $\deh$ to $M_1$, we have
found a von Neumann algebra with comultiplication $(M_1,\de_1)$.
In order to produce invariant weights, we first prove that
$\Rh(M_1) \subset M_1$. Verifying the following equality on a
slice of $\Wh$, we easily arrive at the formula
$$(\io \ot \muh)\deh(x) = ((\teh \ot \io)\deh(x))_{213}
\quad\text{for all}\;\;x \in \Mh \; .$$ Let now $x \in M_1$. Then
$\deh(x) \in M_1 \ot M_1$, so that
$$\deh(x)_{13} = (\io \ot \muh)\deh(x) = ((\teh \ot
\io)\deh(x))_{213} \; .$$ So,
\begin{equation} \label{amaidasgoe}
\deh(x) \in \Mh^{\teh} \ot M_1 = \Rh(M_1) \ot M_1 \; ,
\end{equation}
because of the relation between $\muh$ and $\teh$. If we regard
the restriction of $\dehop$ as a map from $M_1$ to $\Mh \ot M_1$,
then it will be an action of $(\Mh,\dehop)$ on $M_1$. But then we
know that the $\si$-strong$^*$ closure of
$$\{ (\om \ot \io)\dehop(x) \mid x \in M_1, \om \in \Mh_* \}$$
equals $M_1$. Combining this with Equation~\eqref{amaidasgoe}, we
find that $M_1 \subset \Rh(M_1)$. Applying $\Rh$, we get the
equality $M_1 = \Rh(M_1)$. In particular, we also have $M_1 =
\Mh^{\teh}$.

Because the restriction of $\Rh$ to $M_1$ will be an
anti-automorphism of $M_1$ anti-commuting with the
comultiplication $\de_1$, it now suffices to produce a left invariant
weight on $(M_1,\de_1)$, in order to get that $(M_1,\de_1)$ is a
l.c.\ quantum group. Choose an arbitrary n.s.f.\ weight $\eta$ on
$M_1$. Because $\muh$ is integrable, also $\teh$ is integrable and
we can define an n.s.f.\ operator valued weight $T$ from $\Mh$ to
$M_1=\Mh^{\teh}$ by the formula $T(z) = (\vfih_2 \ot \io)\teh(z)$
for all $z \in \Mh^+$. Defining $\etatil = \eta T$, we get an
n.s.f.\ weight $\etatil$ on $\Mh$. We claim that the weight
$\etatil$ is invariant under the action $\muh$. In fact, by
verifying the next formula on slices of $\Wh$, we easily get that
$$(\io \ot \teh)\muh(x) = ((\io \ot \muh)\teh(x))_{213}
\quad\text{for all}\;\;x \in \Mh \; .$$ So, for all $z \in \Mh^+$,
$\om \in \Mh_{2,*}^+$ and $\om' \in \Mh_*^+$, we get
$$\om'(T((\om \ot \io)\muh(z))) = (\om \ot \om')\muh(T(z)) =
\om(1) \; \om'(T(z)) \; ,$$ because $T(z)$ belongs to the extended
positive part of $M_1$. Hence, \linebreak $\etatil((\om \ot \io)\muh(z)) =
\om(1)\; \etatil(z)$, proving our claim.

Above, we already observed that $\vfih$ is invariant under $\muh$.
It then follows from Lemma~3.9 in \cite{SV1} that the Connes
cocycle $u_t = [D\vfih : D \etatil]_t$ belongs to $M_1$ for all $t
\in \R$. From the theory of operator valued weights, we know that
$\si^{\etatil}_t (x) = \si^\eta_t(x)$ for all $x \in M_1$. Hence,
$(u_t)$ is a cocycle with respect to the modular group
$(\si_t^\eta)$ on $M_1$. So, there exists a unique n.s.f.\ weight
$\vfi_1$ on $M_1$ such that $[D\vfi_1 : D\eta]_t = u_t$. Define
$\vfitil_1 = \vfi_1 T$. From operator valued weight theory, we get
$$[D \vfitil_1 : D \etatil]_t = [D\vfi_1 : D\eta]_t = u_t = [D\vfih : D
\etatil]_t\; ,$$ which yields $\vfitil_1 = \vfih$. Let now $x \in
\Mh^+$. Because $$(\io \ot \teh)\deh(x) = ((\io \ot
\deh)\teh(x))_{213} \; ,$$ we find, for all $\om,\om' \in
\Mh_*^+$,
$$\om'(T((\om \ot \io)\deh(x))) = (\om \ot \om') \de_1(T(x)) \;
.$$ When $x \in \Mh^+$ is such that $T(x)$ is bounded, we conclude
that
$$\vfitil_1((\om \ot \io)\deh(x)) = \vfi_1((\om \ot
\io)\de_1(T(x))) \; .$$ Because $\vfitil_1=\vfih$, the left hand
side equals $\om(1)\; \vfitil_1(x) = \om(1) \; \vfi_1(T(x))$.
Hence, for all $x \in \Mh^+$ such that $T(x)$ is bounded and for
all $\om \in M_{1,*}^+$, we find
$$\om(1) \; \vfi_1(T(x)) = \vfi_1((\om \ot \io)\de_1(T(x))) \; .$$
Take an increasing net $(u_i)$ in $\Mh^+$ such that $T(u_i)$
converges increasingly to $1$. Take $y \in M_1$. By lower
semi-continuity, we get
\begin{align*}
\om(1) \; \vfi_1(y^* y) = \sup_i \om(1) \; \vfi_1(T(y^* u_i y)) &=
\sup_i \vfi_1((\om \ot \io)\de_1(T(y^* u_i y))) \\ &= \vfi_1((\om
\ot \io)\de_1(y^*y)) \; .
\end{align*}
Hence, $\vfi_1$ is an n.s.f.\ left invariant weight on
$(M_1,\de_1)$ and the latter is a l.c.\ quantum group.

Define $\be$ to be the identity map, embedding $M_1$ into $\Mh$.
In order to obtain that
$$(M_2,\det) \overset{\al}{\lrecht} (M,\de) \overset{\be}{\lrecht}
(\Mh_1,\deoh)$$ is a short exact sequence, it remains to show that
$\al(M_2) = M^\te$, where $\te$ is the canonical action of
$(\Mh_1,\dehopo)$ on $M$, associated to $\be$ (see Proposition~3.1
in \cite{VV}). Because $\te = (\Rh_1 \ot R)\mu R$ and
$R(\al(M_2))=\al(M_2)$, it suffices to show that $\al(M_2)=M^\mu$.
From Proposition~3.1 in \cite{VV}, we immediately deduce that
$M^\mu = M \cap \be(M_1)'$. Above we already saw that $M_1 = \Mh
\cap \al(M_2)'$. Because $\al(M_2)$ is a two-sided coideal of
$(M,\de)$, it follows from Th{\'e}or{\`e}me~3.3 in \cite{E1} that
$$M \cap (\Mh \cap \al(M_2)')' = \al(M_2) \; .$$ So, we have a
short exact sequence of l.c.\ quantum groups.

Finally, we should prove the uniqueness of this short exact
sequence up to isomorphism. Suppose that we have another short
exact sequence
$$(M_2,\det) \overset{\al}{\lrecht} (M,\de) \overset{\gamma}{\lrecht}
(\Mh_3,\deh_3) \; .$$ We still have the same action $\muh$ of
$(\Mh_2,\deh_2)$ on $\Mh$ and a reasoning as in the previous
paragraph yields that $\gamma(M_3) = \Mh^{\muh}$. Hence, it
follows that $\gamma$ gives an isomorphism of l.c.\ quantum groups
between $(M_3,\de_3)$ and the l.c.\ quantum group $(M_1,\de_1)$
constructed above.
\end{proof}

\section{Matched pairs of Lie groups and Lie algebras}

In what follows we consider Lie groups and Lie algebras over the field $k=\C$
or $\R$.

\begin{definition} \label{def2.1}
We call $(G_1,G_2)$ a matched pair of Lie groups if, in
Definition~\ref{mpg}, $G$ is a Lie group.
\end{definition}

Observe that it follows from Proposition~\ref{LieOK} that $\te$ is a
diffeomorphism of $G_1 \times G_2$ onto the open subset $\Om$ of $G$.

The infinitesimal form of this definition is as follows (see \cite{Majbook}).

\begin{definition} \label{def2.2}
We call $(\lieg_1,\lieg_2)$ a matched pair of Lie algebras, if there
exists a Lie algebra $\lieg$ with Lie subalgebras $\lieg_1$ and
$\lieg_2$ such that $\lieg = \lieg_1 \oplus \lieg_2$ as vector spaces.
\end{definition}

These conditions are equivalent to the existence of a left action $\triangleright:
\lieg_2\otimes \lieg_1\to \lieg_1$ and a right action $\triangleleft:
\lieg_2\otimes \lieg_1\to \lieg_2$, so that $\lieg_1$ is a left $\lieg_2$-module and
$\lieg_2$ is a right $\lieg_1$-module and
\begin{enumerate}
\item $x\triangleright[a,b]=[x\triangleright a,b]+[a,x\triangleright b]+
(x\triangleleft a)\triangleright b-(x\triangleleft b)\triangleright a,
$
\item $[x,y]\triangleleft a=[x,y\triangleleft a]+[x\triangleleft a,y]+
x\triangleleft (y\triangleright a)-y\triangleleft (x\triangleright a),
$
\end{enumerate}
for all $a,b\in \lieg_1,\ x,y\in \lieg_2$. Then, for the decomposition of vector
spaces above we have
$$
[a\oplus x,b\oplus y]=([a,b]+x\triangleright b - y\triangleright a)\oplus
([x,y]+x\triangleleft b-y\triangleleft a)
$$
(see \cite{Majbook}, Proposition 8.3.2).

Two matched pairs of Lie algebras, $(\lieg_1,\lieg_2)$ and
$(\lieg'_1,\lieg'_2)$, are called isomorphic if there is an isomorphism of the
corresponding Lie algebras $\lieg$ and $\lieg'$ sending $\lieg_i$ onto
$\lieg'_i\ (i=1,2)$.

Let us explain the relation between the two notions of a matched pair.

\begin{proposition} \label{prop2.3}
Let $(G_1,G_2)$ be a matched pair of Lie groups in the sense of
Definition~\ref{def2.1}. If $\lieg$ denotes the Lie algebra of $G$,
and if $\lieg_1$, resp.\ $\lieg_2$, are the Lie subalgebras corresponding to
the closed subgroups $i(G_1)$, resp.\ $j(G_2)$, then
$(\lieg_1,\lieg_2)$ is a matched pair of Lie algebras.
\end{proposition}
\begin{proof}
The fact that $\lieg = \lieg_1 \oplus \lieg_2$ as vector spaces
follows from the fact that $\te$ is a diffeomorphism in the neighbourhood of
the unit
element.
\end{proof}

The converse problem, to construct a matched pair of Lie groups
from a given matched pair $(\lieg_1,\lieg_2)$ of Lie algebras, is
much more subtle. Indeed, one can take, of course, the connected,
simply connected Lie group $G$ of the corresponding $\lieg$ and
find unique connected, closed subgroups $G_1$ and $G_2$ of $G$
whose tangent Lie algebras are $\lieg_1$ and $\lieg_2$,
respectively. However, in the proof of the following proposition, we
see that $(G_1,G_2)$
is not necessarily a matched pair of Lie groups even if
$\dim \lieg_1=\dim \lieg_2=1$.

\begin{proposition} \label{propnieuw}
Every matched pair of complex Lie algebras $\lieg_1=\lieg_2=\C$ can be
exponentiated to a matched pair of Lie groups $(G_1,G_2)$ where
$G_1,G_2$ are either $(\C,+)$ or $(\C \setminus \{0\},\cdot)$.

Every matched pair of real Lie algebras $\lieg_1=\lieg_2=\R$ can be
exponentiated to a matched pair of Lie groups $(G_1,G_2)$ where
$G_1,G_2$ are either $(\R,+)$ or $(\R \setminus \{0\},\cdot)$.
\end{proposition}
\begin{proof}
Consider first the complex case. The only two-dimensional complex Lie algebras are the abelian one and
the one with generators $X,Y$ and relation $[X,Y] =Y$. If $\lieg$ is
abelian, the mutual actions of $\lieg_1$ and $\lieg_2$ on each other
are trivial and exponentiation is obviously a direct sum.

If $\lieg$ is generated by $[X,Y]=Y$, we either have that $\lieg_1$ or
$\lieg_2$ is equal to $\C Y$, in which case one of the actions is
trivial and $G$ can be constructed as semi-direct product of the connected, 
simply connected Lie groups of $\lieg_1$ and $\lieg_2$, or we have 
that both $\lieg_1$ and $\lieg_2$ differ from
$\C Y$. In the latter case, there is, up to isomorphism, only one
possibility, namely $\lieg_1 = \C X$, $\lieg_2 = \C (X+Y)$. Define on
$\C \setminus \{0\} \times \C$ the Lie group with product
$$
(t,s)(t',s')=(tt',s+ts') \; .
$$
Define $G_1=G_2=\C \setminus \{0\}$ with embeddings $i(g) = (g,0)$ and
$j(s) = (s,s-1)$, we indeed get a matched pair of complex Lie
groups with mutual actions
\begin{equation} \label{exp11}
\al_g(s) = g(s-1) + 1 \; , \quad \be_s(g) = \frac{sg}{g(s-1) +
1} \; .
\end{equation}
The real case is completely analogous.
\end{proof}

\begin{remark} \label{remnieuw}
The connected simply
connected complex Lie group $G$ of $\lieg$ consists of all pairs 
$(t,s)$ with $t,s\in \C$ and the product
$$
(t,s)(t',s')=(t+t',s+\exp(t)s')
$$
(see, for example, \cite{FH}, {\S}10.1), and its closed subgroups $G_1$ 
and $G_2$ corresponding to the decomposition $\lieg = \C X \oplus \C (X+Y)$ 
above consist respectively of all pairs of the form
$(g,0)$ and $(s,\exp(s)-1)$ with $g,s\in \C$. These groups do not form a matched
pair because $G_1\cap G_2=\{(2\pi in,0)\vert n\in \Z\}=Z(G)$. So, it
is crucial not to take $G$ simply connected above.

Allowing $g,t,s$ above to be only real, we come to the example of a matched pair
of real Lie groups from \cite{VV}, Section~5.3. Here $\lieg$ is a
real Lie algebra generated by $X$ and $Y$ subject to the relation $[X,Y]=Y$ and
one considers the decomposition $\lieg=\R X\oplus\R(X+Y)$.
Then, to get a matched pair of Lie groups, we consider $G$ as the variety
$\R \setminus \{0\} \times \R$ with the product
$$
(s,x)(t,y) = (st,x+sy)
$$
and embed $G_1=G_2=\R \setminus \{0\}$ by the formulas $i(g) = (g,0)$ and
$j(s) = (s,s-1)$.
Remark that here, it is impossible to take the connected component of the
unity of the group of affine transformations of the real line as $G$, because
it is easy to see that for its closed subgroups $G_1$ and $G_2$ corresponding to
the above mentioned subalgebras, the set $G_1G_2$ is not dense in $G$.
\end{remark}

The next example shows that in general, for a given matched pair of Lie 
algebras, it is even possible that $G_1 \cap G_2 \neq \{e\}$
for {\em any} corresponding pair of Lie groups, which means that
such a matched pair of Lie algebras cannot be exponentiated to a matched pair of
Lie groups in the sense of Definition~\ref{def2.1}.

\begin{example} \label{ex2.4}
Consider a family of complex Lie algebras $\lieg=\text{span}\{X,Y,Z\}$ with $[X,Y]=Y,\
[X,Z]=\al Z,\ [Y,Z]=0$, where $\al\in \C\setminus \{0\}$, and the decomposition
$\lieg=\text{span}\{X,Y\}\oplus\C (X+\al Z)$. The
corresponding connected simply connected complex Lie group $H$ consists of all triples
$(t,u,v)$ with $t,u,v\in \C$ and the product
$$(t,u,v)(t',u',v')=(t+t',u+\exp(t)u',v+\exp(\al t)v')$$
(see, for example, \cite{FH}, {\S}10.3), and its closed subgroups $H_1$ and $H_2$
corresponding to the decomposition above consist respectively of all triples of the form
$(t,u,0)$ and $(s,0,\exp(\al s)-1)$ with $t,u,s\in \C$. These groups do not form a matched
pair because $H_1\cap H_2=\{({\frac{2\pi in}{\al}},0,0)\vert n\in
\Z\}$.

We claim that, if $1/\al \not\in \Z$ and if $G$ is any complex Lie group with Lie
algebra $\lieg$, such that $G_1,G_2$ are closed subgroups of $G$ with
tangent Lie algebras $\lieg_1$, resp.\ $\lieg_2$, then $G_1 \cap G_2 \neq
\{e\}$. Indeed, since the Lie group $H$ is connected and simply
connected, the connected component $G^{(e)}$ of $e$ in $G$ can be identified
with the quotient of $H$ by a
discrete central subgroup. If $\al \not\in \Q$, the center of $H$ is
trivial, so that
we can identify $G^{(e)}$ and $H$. Under this identification, the
connected components of $e$ in $G_1,G_2$ agree with
$H_1,H_2$. Because $H_1 \cap H_2 \neq \{e\}$, our claim follows. If
$\al =\frac{m}{n}$ for $m,n \in \Z \setminus \{0\}$ mutually prime, the center of $H$
consists of the elements $\{(2\pi n N,0,0) \mid N \in \Z \}$. Hence,
the different possible quotients of $H$ are labeled by $N \in \Z$ and
are given by the triples $(a,u,v) \in \C^3$, $a \neq 0$ and the
product
\begin{equation} \label{group-law}
(a,u,v) (a',u',v') = (aa',u+a^{nN}u',v+a^{mN}v') \; .
\end{equation}
The
closed subgroups corresponding to $\lieg_1$ and $\lieg_2$ are given by
$(a,u,0)$ and $(b,0,b^{mN} - 1)$ with $a,b,u \in \C$ and $a,b \neq
0$. The intersection of both subgroups is non-trivial whenever $mN
\neq \pm 1$. This proves our claim.

Considering now the complex Lie algebras above as real Lie algebras with
generators $X,iX,Y,iY,Z,iZ$ and the decomposition above as a decomposition of
real Lie algebras, we get a matched pair of real Lie algebras which
cannot be exponentiated to a matched pair of real Lie groups.

In the remaining case $\al=1/n$ with $n \in \Z \setminus \{0\}$, we
can consider the Lie group $G$ defined by Equation~\eqref{group-law}
with $m=N=1$. Consider $G_1 = \C \setminus \{0\} \times \C$ with
$(a,u)(a',u') = (aa',u+a^n u')$ and $G_2 = \C \setminus
\{0\}$. Writing $i(a,u) = (a,u,0)$ and $j(v) = (v,0,v-1)$, we get a
matched pair of Lie groups with mutual actions
$$\al_{(a,u)}(v) = a(v-1)+1 \quad\text{and}\quad \be_v(a,u) = \bigl(
\frac{va}{a(v-1)+1}, \frac{u}{(a(v-1)+1)^n} \bigr)\; .$$
\end{example}

In the next section, we study more closely the exponentiation of a
matched pair of Lie algebras when one of the Lie algebras has
dimension 1.

\section{Matched pairs of Lie groups and Lie algebras in
  dimension $n+1$}

We use systematically the following terminology.

\begin{terminology} \label{term3.1}
A matched pair $(\lieg_1,\lieg_2)$, resp.\ $(G_1,G_2)$, is said to be
of dimension $n_1+n_2$ if the dimension of $\lieg_i$, resp.\ $G_i$, is
$n_i$.
\end{terminology}

Suppose that $\lieg = \lieg_1 \oplus \lieg_2$ is a matched pair of
Lie algebras with $\Dim \lieg_2 = 1$. Put $\lieg_2 = k A$. For
all $X \in \lieg_1$, we define $\be(X) \in \lieg_1$ and $\chi(X)
\in k$ such that $$[X,A] = \be(X) + \chi(X) A \; .$$ Then, $\be$
and $\chi$ are linear, and, for all $X,Y \in \lieg_1$, the Jacobi
identity for $\lieg$ gives:
\begin{align}
\chi([X,Y]) &= 0 \; , \label{eq1} \\
\be([X,Y]) &= [X,\be(Y)] + [\be(X),Y] + \be(X) \chi(Y) - \be(Y)
\chi(X) \;  . \notag
\end{align}
By induction, one verifies that
\begin{align*}
\be^n([X,Y]) &= \sum_{k=0}^n
\begin{pmatrix} n \\ k \end{pmatrix} [\be^k(X) , \be^{n-k}(Y)] \\ &\quad +
\sum_{k=1}^n \begin{pmatrix} n \\ k-1 \end{pmatrix} \bigl( \be^k(X)
\chi(\be^{n-k}(Y)) - \be^k(Y) \chi(\be^{n-k}(X))  \bigr) \; .
\end{align*}
Hence,
\begin{equation} \label{eq2}
\chi(\be^n([X,Y])) = n \bigl( \chi(\be^n(X)) \chi(Y) -
\chi(\be^n(Y)) \chi(X) \bigr) \; .
\end{equation}
Then, we claim that the linear forms $\chi, \chi \be$ and $\chi \be^2$ are
linearly dependent. If not, we find $X_0,X_1,X_2 \in \lieg_1$, such
that $\chi(\be^i(X_j)) = \delta_{ij}$ for $i,j \in \{0,1,2\}$, where
$\delta_{ij}$ is the Kronecker symbol. Because $\chi(X_1) =
\chi(X_2)=0$, we get $\chi(\be^n([X_1,X_2])) = 0$ for all $n$. Define
$$\lieg_0 = \bigcap_{i=0}^2 \Ker \chi \be^i \; .$$ Then,
$\be([X_1,X_2]) \in \lieg_0$. Using Equation~\eqref{eq1}, we get
\begin{equation} \label{eq3}
\be([X_1,X_2]) = [X_1,\be(X_2)] + [\be(X_1),X_2] \; .
\end{equation}
On the other hand, using Equation~\eqref{eq2}, it follows that
$[X_1,\be(X_2)] \in \lieg_0$, because $\chi(\be(X_2)) =
\chi(X_1)=0$. Combining this with Equation~\eqref{eq3}, we get that
$[\be(X_1),X_2] \in \lieg_0$. Nevertheless, using once again
Equation~\eqref{eq2}, we get that $\chi(\be^2([\be(X_1),X_2])) = -2$,
contradicting the fact that $[\be(X_1),X_2] \in \lieg_0$. So, we have
proved that $\chi, \chi \be$ and $\chi \be^2$ are linearly dependent.

Hence, we can separate three different possibilities.
\begin{list}{}{\setlength{\labelwidth}{1.2cm}\setlength{\leftmargin}{1.5cm}
\setlength{\labelsep}{.5em}\setlength{\itemindent}{0cm}}
\item[Case 1.] $\chi = 0$.
\item[Case 2.] $\chi \neq 0$ and $\be(\Ker \chi) \subset \Ker \chi$. \label{casetwo}
\item[Case 3.] $\chi$ and $\chi \be$ linearly independent, and $\be(\Ker \chi
  \cap \Ker \chi \be)  \subset \Ker \chi \cap \Ker \chi \be \; .$
\end{list}

\centerline{\bf Case~1} 

The action of $\lieg_1$ on $\lieg_2$ is trivial and
$\lieg_2$ acts on $\lieg_1$ by automorphisms. To exponentiate such a
matched pair it suffices to use a semi-direct product of $k$ and the
connected, simply connected Lie group $G_1$ of $\lieg_1$.

\centerline{\bf Case~2} 

In this case $\lieg_0:=\Ker(\chi)$ is an ideal of $\lieg$, on
which $\lieg_2$ acts as an automorphism group. There exists an $a \in
k$ such that $\chi \be = a \chi$, by assumption. Take $X_0 \in
\lieg_1$ such that $\chi(X_0) = 1$. Then, we get
$$
[X_0,A]= A + aX_0 + Y_0\quad (Y_0\in \lieg_0) \; .
$$
Putting $\Atil = A+aX_0$, we get
$[X_0,\Atil]= \Atil +Y_0$. This suggests how to
exponentiate $\lieg$. We start with the connected,
simply connected Lie group $G_0$ of $\lieg_0$ and observe that,
due to the action of $\Atil$, $k$ acts by automorphisms
$(\rho_x)_{x \in k}$ on $G_0$. We exponentiate $k
\Atil + \lieg_0$ on the space $k \times G_0$ with product $$(x,g) (y,h) =
(x+y, g \rho_x(h)) \; .$$ Next, we observe that, due to the action
of $X_0$, $k$ is acting by automorphisms on $k \ltimes G_0$, and
one would suggest to make another semi-direct product.
But in this case the subgroups corresponding to $\lieg_1$ and $k A$
do not form a matched pair of Lie
groups if $a \neq 0$.  The subgroups corresponding to $\lieg_1$ and
$k A$ do form a matched pair of Lie groups when $a =0$. So, in what
follows, we treat the case $a \neq 0$. Example \ref{ex2.4} shows that 
if $k=\C$, the exponentiation is in general impossible if 
$\dim \lieg_1 \geq 2$. So, we restrict to the case $k=\R$.

We first prove some general results which allow to obtain the 
exponentiation whenever $n \leq 4$. Afterwards, we will give an
example showing that if dimension $n\geq 5$, the exponentiation is 
in general impossible.

Writing $\mu = \Ad X_0$ and $\rho = \Ad \Atil$ on
$\lieg_0$, we are given $\mu$ and $\rho$, derivations of $\lieg_0$ and an 
element $Y_0 \in \lieg_0$ such that $[\mu,\rho]=\rho + \Ad Y_0$. Further, we 
have $[X_0,\Atil]=\Atil + Y_0$,
$\lieg_1 = \R X_0 + \lieg_0$ and $\lieg_2 = \R A$, where $A = \Atil -
a X_0$.

We introduce the notation $\R^*$ for the Lie group $\R \setminus
\{0\}$ with multiplication and $\R^*_+$ for the subgroup of elements
$s > 0$.

\begin{proposition} \label{exponentiation_general}
The matched pair $(\lieg_1,\lieg_2)$ has an exponentiation with at most
two connected components in the following cases:
\begin{enumerate}
\item $\rho$ is inner and the center of $\lieg_0$ is trivial.
\item $\lieg_0$ is abelian.
\item $\lieg_0 = \langle X,Y \rangle \oplus \liez_0$, with $[X,Y]=Y$
  and $\liez_0$ central in $\lieg_0$.
\item $\lieg_0=\langle X,Y,Z \rangle$ with $[X,Y]=Z$ and $Z$ central.
\end{enumerate}
\end{proposition}

In the proof of this proposition we use systematically the following lemma.

\begin{lemma} \label{closedsubgroup}
Let $G_0$ be a Lie group with Lie algebra $\lieg_0$, with center
$\liez_0$. Suppose that $[\lieg_0,\lieg_0] \subset
\liez_0$. Let $(\mu_s)$ be an action of $\R^*$ by automorphisms of
$G_0$. Denote by $\mu$ the derivation of $\lieg_0$ corresponding to
$(\mu_s)_{s > 0}$ and denote by $\te$ the involutive automorphism
$\te:=d \mu_{-1}$, giving rise to a decomposition $\lieg_0 = \lieg_0^+
\oplus \lieg_0^-$.

So, $\mu$ leaves $\lieg_0^-$ invariant and we suppose that $\mu$ is
invertible on $\lieg_0^-$. Further, $\mu$ leaves $[\lieg_0,\lieg_0]$ 
invariant and we suppose that $\mu$ is invertible on $[\lieg_0,\lieg_0]$ 
as well.

Define $G:= \R^* \ltimes G_0$ with product $(s,g)(t,h) = (st,g
\mu_s(h))$ and Lie algebra $\R X_0 + \lieg_0$, where $X_0$ denotes the
canonical generator of $\lieg$ corresponding to $\R^*$. Let $C \in
\lieg_0$.  Then, there exists a closed
subgroup
$$K= \{(s,v(s)) \mid s \in \R^* \}$$
of $G$, with tangent Lie algebra $\R (X_0 + C)$, where $v : \R^*
\recht G_0$ is a smooth function.
\end{lemma}
\begin{proof}
Denote by $\lexp_\lieg$ the exponential mapping $\lieg \recht
G$. Then, we get a smooth function $v:\R^*_+ \recht G_0$ such that
$\lexp_\lieg((\log s)(X_0 + C)) = (s,v(s))$ for $s > 0$. Then, $v(1)
= e$ and $v'(1) = C$. Further, $v(st) = v(s) \mu_s(v(t))$ for all
$s,t > 0$. We say that $v$ is a $(\mu_s)$-cocycle. So, we are looking 
for an element $g_0 \in \lieg_0$, so
that we can define $v(-1) = g_0$. Then, $(-1,g_0) \in K$ and, for $s >
0$,
$$(-1,g_0) (s,v(s)) = (-s,g_0 \mu_{-1}(v(s))) \; , \quad
(s,v(s))(-1,g_0) = (-s,v(s) \mu_s(g_0)) \; .$$
Hence, we are done if we can find an element $g_0 \in G_0$ such that
$g_0 \mu_{-1}(v(s))= v(s) \mu_s(g_0)$, because than we can define
$v(-s)$ to be this expression.

Denote by $\lexp_{\lieg_0}$ the exponential mapping $\lieg_0 \recht
G_0$. Then, we look for $D \in \lieg_0$ such that
\begin{equation} \label{eeneerstegelijkheid}
\mu_{-1}(v(s)) = \lexp_{\lieg_0}(-D) v(s)
\lexp_{\lieg_0}(\exp((\log s) \mu)(D)) \; .
\end{equation}
We want to derive at $s = 1$. For this, observe that
\begin{align*}
\lexp_{\lieg_0}(-D) & \lexp_{\lieg_0}(\exp( (\log s) \mu)(D)) \\ &=
\lexp_{\lieg_0}\bigl(-D + \exp( (\log s) \mu)(D) - \frac{1}{2}
[D,\exp( (\log s) \mu)(D)] \bigr) \; ,
\end{align*}
where we have used that
$[\lieg_0,\lieg_0] \subset \liez_0$. Taking the derivative at $s=1$ of
Equation~\eqref{eeneerstegelijkheid}, we look for $D \in \lieg_0$
such that
$$\te(C) = (\Ad \lexp_{\lieg_0}(-D))(C) + \mu(D) -
\frac{1}{2} [D,\mu(D)] \; .$$
The equation becomes
$$\te(C) = \exp(- \Ad D)(C) + \mu(D) -
\frac{1}{2} [D,\mu(D)] \; .$$
Using once again that $[\lieg_0,\lieg_0] \subset \liez_0$, we get the
equation
$$\te(C) - C = \mu(D) - [D,C] -
\frac{1}{2} [D,\mu(D)] \; .$$
Define $Y_2 = \mu^{-1}(\te(C) - C)$, which is possible because
$\mu$ is invertible on $\lieg_0^-$. Next, define $Z_2 :=
\frac{1}{2}\mu^{-1}([Y_2,\te(C)+C])$, which is possible because
$\mu$ is invertible on $[\lieg_0,\lieg_0]$. Then, $Z_2 \in \liez_0$
and we define $D= Y_2 + Z_2$. Then,
\begin{align*}
\mu(D) & - [D,C] - \frac{1}{2} [D,\mu(D)] \\ &= \te(C) - C
 + \frac{1}{2}[Y_2,\te(C)+C] - [Y_2,C] - \frac{1}{2}
 [Y_2,\te(C) - C] \\ &= \te(C) - C  \; .
\end{align*}
So, writing $g_0 = \lexp_{\lieg_0}(D)$, we may conclude that the
left and right hand side of the to be proven equation
$$\mu_{-1}(v(s))= g_0^{-1} v(s) \mu_s(g_0)$$
have the same derivative at $s=1$. But, both sides of the equations
are $(\mu_s)_{s>0}$-cocycles and hence, both sides are equal for all
$s>0$.
\end{proof}

\begin{proof}[Proof of Proposition~\ref{exponentiation_general}]
{\it Part 1 : $\rho$ is inner and the center of $\lieg_0$ is trivial.}

Take the unique $B \in \lieg_0$ such that $\rho = \Ad B$ on
$\lieg_0$. Take a new generator $\Ah = \Atil - B$ in $\lieg$. Because
the center of $\lieg_0$ is trivial and
$$0 = [\mu, \rho - \Ad B] = \rho + \Ad (Y_0 - \mu(B)) = \Ad (Y_0 + B -
\mu(B)) \; ,$$
we get $Y_0 + B - \mu(B) = 0$ and hence $[X_0,\Ah] = \Ah$. Because
$[\Ah, \lieg_0]=\{0\}$, the connected, simply connected Lie group of
$\lieh:=\R \Ah + \lieg_0$ is given by $H:=\R \oplus G_0$, where $G_0$ is the
connected, simply connected Lie group of $\lieg_0$. Using the
derivation $\Ad X_0$ on $\lieh$, we get an action $(\mu_s)$ of
$\R^*_+$ on $H$ of the form
$$\mu_s(x,g)=(sx,\eta_s(g)) \; , \quad\text{where}\;\; s > 0, x \in
\R, g \in G_0 \; , $$ and where $(\eta_s)$ is an action of $\R^*_+$ on $G_0$. We
easily extend $(\mu_s)$ to an action of $\R^*$, defining
$\mu_s(x,g)=(sx,\eta_{|s|}(g))$. So, we can define a Lie group $G:=
\R^* \ltimes H$, with Lie algebra $\lieg$, on the space $\R^* \times
\R \times G_0$ with product
$$(s,x,g)(t,y,h) = (st,x+sy,g\eta_{|s|}(h)) \; .$$
Define the closed subgroup $G_1$ consisting of the elements $(s,0,g)$,
where $s \in \R^*$ and $g \in G_0$. The tangent Lie algebra of $G_1$
is precisely $\lieg_1$. To define the closed subgroup $G_2$, recall
that $A = \Atil - a X_0 = \Ah - a X_0 + B$. Using the exponential
mapping of $G$, we find a smooth function $v : \R^*_+ \recht G_0$ such
that
$$\{(s, \frac{1}{a} (1-s) , v(s) ) \mid s \in \R^*_+ \}$$ is a closed
subgroup of $G$ with tangent Lie algebra $\R A$. We define
$$G_2 := \{(s, \frac{1}{a} (1-s) , v(|s|) ) \mid s \in \R^* \}$$ and
then, one verifies that $G_2$ is a closed subgroup of $G$ with tangent
Lie algebra $\R A$. It is easy to see that $(G_1,G_2)$ is a matched
pair of Lie groups.

{\it Part 2 : $\lieg_0$ is abelian.}

In this case, $\mu$ and $\rho$ are linear transformations of
$\lieg_0$ satisfying $\mu \rho - \rho \mu = \rho$. Then, for any polynomial $P$, we get
\begin{equation}\label{wateenmooitje}
P(\mu) \rho =\rho P (\mu + \io) \; .
\end{equation}
Let $V$ be the complexified vector space of the real vector space
$\lieg_0$ and consider $\mu$ as a linear operator on $V$, which we still denote by
$\mu$. Then, we have a
direct sum decomposition $$V = \bigoplus_{\lambda \in \C} E_\lambda \;
,$$ where $E_\lambda$ is the generalized eigenspace corresponding to
$\lambda \in \C$. A vector $X \in V$ belongs to $E_\lambda$ if and
only if
$(\mu - \lambda \io)^n X = 0$ for $n$ big enough. We also get a direct
sum decomposition $$\lieg_0 = \bigoplus_{r \in \R} F_r \; ,$$ where $F_r$ is
such that $$\C F_r = \bigoplus_{\lambda, \Re(\lambda)=r} E_\lambda \; .$$
The subspaces $F_r$ are invariant under $\mu$.
Also $\rho$ extends to $V$ and using
Equation~\eqref{wateenmooitje}, it is clear that $\rho(E_\lambda)
\subset E_{\lambda+1}$. Hence, $\rho(F_r)=F_{r+1}$. Denote by
$\vep(r)$ the entire part of $r \in \R$, such that $\vep(r) \in \Z$
and $\vep(r) \leq r < \vep(r)+1$. Then, we define
$$\lieg_0^+ = \bigoplus_{\vep(r) \;\text{is even}} F_r \quad\text{and}\quad
\lieg_0^-=\bigoplus_{\vep(r) \;\text{is odd}} F_r \; .$$
Defining $\te(X) =X$ for $X \in \lieg_0^+$ and $\te(X) = -X$ for $X
\in \lieg_0^-$, we obtain an involution $\te$ of $\lieg_0$ satisfying
$\te\mu = \mu\te$ and $\te \rho = - \rho \te$. Observe that $\mu$
leaves the subspaces $\lieg_0^+$ and $\lieg_0^-$ globally invariant
and that $\mu - \lambda \io$ is invertible on $\lieg_0^+$ when
$\vep(\lambda)$ is odd, while $\mu - \lambda \io$ is invertible on
$\lieg_0^-$ when $\vep(\lambda)$ is even.

Define the Lie algebra $\lieh := \R \Atil + \lieg_0$. Its connected,
simply connected Lie group $H$ lives on the space $\R \times \lieg_0$,
with product $$(x,X)(y,Y) = (x+y,X + \exp(x \rho)(Y)) \; .$$ We extend
the derivation $\mu = \Ad X_0$ to $\lieh$ and observe that
$\mu(\Atil)=\Atil + Y_0$. Defining $\te(\Atil) = - \Atil + Z_0$, where
$$Z_0 = (\mu - \io)^{-1}(\te(Y_0) + Y_0) \in \lieg_0^+ \; ,$$
one verifies that $\te$ is an involutive automorphism of $\lieh$
commuting with the derivation $\mu$. Putting together $\mu$ and
$\te$, we obtain an action $(\mu_s)$ of $\R^*$ on $H$ such that
$\mu_s(x,X) = (sx, \eta_s(X) u(s,x))$, where $(\eta_s)$ is an action
of $\R^*$ on $\lieg_0$. We define the Lie group $G := \R^*
\ltimes H$, which lives on the space $\R^* \times \R \times \lieg_0$, with
product $$(s,x,X)(t,y,Y) = (st, x+sy, X+\exp(x \rho)(\eta_s(Y)u(s,y)))
\; .$$
Define the closed subgroup $G_1$ consisting of the elements $(s,0,X)$,
where $s \in \R^*$ and $X \in \lieg_0$. The tangent Lie algebra of $G_1$
is precisely $\lieg_1$. Finally, we have to find a closed subgroup
$G_2$ with tangent Lie algebra $\R A = \R(\Atil - a X_0)$, which
consists of the elements $(s,\frac{1}{a} (1-s),v(s))$, $s \in \R^*$
and $v : \R^* \recht \lieg_0$ a smooth function.
Conjugating with the element $(1,-\frac{1}{a},0)$, an equivalent
question is to find a closed subgroup with tangent Lie algebra $\R(X_0
+ Z_0)$ (for a certain $Z_0 \in \lieg_0$), consisting of the elements
$(s,0,w(s))$, $s \in \R^*$ and $w : \R^* \recht G_0$ a smooth
function. This is possible applying
Lemma~\ref{closedsubgroup} to the action $(\eta_s)$ of $\R^*$ on $\lieg_0$. 
Then, it is clear that $(G_1,G_2)$ form a matched pair of Lie groups.

{\it Part 3 : $\lieg_0 = \langle X,Y \rangle \oplus \liez_0$, with $[X,Y]=Y$ 
and $\liez_0$ central.}

Every derivation of $\lieg_0$ leaves the center $\liez_0$ invariant. On the 
quotient Lie algebra
$\lieg_0 / \liez_0$ every derivation is inner. So, changing the generators
$X_0$ and $\Atil$ of $\lieg$, we may suppose that we are in the
following situation:
$$[X_0,\Atil] = \Atil + Y_0 \; , \quad [X_0,\lieg_0] \subset \liez_0
\; , \quad [\Atil,\lieg_0] \subset \liez_0 \; ,$$
and $\lieg_1 = \R X_0 + \lieg_0$, $\lieg_2 = \R (\Atil - a X_0 + B_1)$
for a certain element $B_1 \in \lieg_0$.

Write $\mu = \Ad X_0$ and $\rho = \Ad \Atil$ as derivations on
$\lieg_0$. Because $\mu$ and $\rho$ preserve $[\lieg_0,\lieg_0]$, we
get $\mu(Y_1) = \rho(Y_1) = 0$. Further,
$[Y_0,\lieg_0]=([\mu,\rho]-\rho)(\lieg_0) \subset \liez_0$. Because
$[\lieg_0,\lieg_0]=\R Y_1$, we get $[Y_0,\lieg_0]=\{0\}$,
which
gives $Y_0 \in \liez_0$ and so, $[\mu,\rho]=\rho$. Suppose $\mu(X_1) = Z_1$ 
and $\rho(X_1) = Z_2$. Because $[\mu,\rho]=\rho$, we get
\begin{equation}\label{onthouden}
\mu(Z_2) - Z_2 = \rho(Z_1) \; .
\end{equation}
As in part~2, we can find an involutive automorphism $\te$ of
$\liez_0$, such that $\te$ commutes with $\mu$ and anti-commutes with
$\rho$ and such that $\mu-\lambda \io$ is invertible on $\liez_0^+$
when $\vep(\lambda)$ is odd and invertible on $\liez_0^-$ when
$\vep(\lambda)$ is even. Defining $Z_3 = \mu^{-1}(\te(Z_1) - Z_1) \in
\liez_0^-$, we can extend $\te$ to an involutive automorphism of $\lieg_0$, by
putting $\te(X_1) = X_1+Z_3$ and $\te(Y_1) = Y_1$. Then, $\te$ commutes
with $\mu$, by definition of $Z_3$, and anti-commutes with
$\rho$, because $\rho(Z_3) = - Z_2 - \te(Z_2)$. This last equality can
be deduced as follows: because $(\mu - \io)\rho = \rho \mu$, we get
$\rho \mu^{-1} = (\mu - \io)^{-1} \rho$ on $\liez_0^-$; in particular,
$$\rho(Z_3) = (\mu - \io)^{-1}(-\te(\rho(Z_1)) - \rho(Z_1)) =
-Z_2-\te(Z_2) \; ,$$ where we used Equation~\eqref{onthouden}.

Next, we define the Lie algebra $\lieh:=\R \Atil + \lieg_0$. We extend
$\te$ to an involutive automorphism of $\lieh$ by defining $\te(\Atil)
= -\Atil + Z_4$, where $Z_4 : = (\mu - \io)^{-1}(\te(Y_0) + Y_0) \in
\liez_0^+$. Extending also $\mu = \Ad X_0$ to $\lieh$ (recall that $\mu(\Atil) =
\Atil + Y_0$), we observe that $\te$ and $\mu$ commute.

As above, the derivation $\rho$ gives rise to an action of $\R$ on
$G_0$, where $G_0$ lives on the space $\R^2 \times \liez_0$ with
product $(x,y,Z)(x',y',Z') = (x+x',y+\exp(x)y',Z+Z')$. Then, $H := \R
\kruisje{\rho} G_0$ gives an exponentiation of $\lieh$, on which the
derivation $\mu$ and the involutive automorphism $\te$ are combined to
produce an action $(\mu_s)$ of $\R^*$ by automorphisms of $H$. We
define $G=\R^* \kruisje{\mu} H$ and the product is given by $(s,x,g)(t,y,h)
= (st,x+sy, \ldots)$. The precise form of the product can be written
as above. The closed subgroup $G_1$ consists again of the elements
$(s,0,g)$ for $s \in \R^*$ and $g \in G_0$. To find $G_2$, we have to
construct, again as above, the closed subgroup $G_2$ of $G$ with tangent Lie
algebra $\R(\Atil - a X_0 + B_1)$, where $B_1 \in
\lieg_0$, and such that $G_2$ consists of the elements
$(s,\frac{1}{a}(1-s),v(s))$, where $s \in \R^*$. Conjugating, the
problem is reduced again to finding a closed subgroup of $G$ with
tangent Lie algebra $\R(X_0 + B_2)$, for some arbitrary $B_2 \in
\lieg_0$, consisting of the elements $(s,0,w(s))$, $s \in
\R^*$. We solve this problem in the closed subgroup $K:=\R^*
\ltimes G_0$ of $G$, whose Lie algebra is $\liek:=\R X_0 + \lieg_0$.
Suppose that $B_2 = x_1 X_1 + y_1 Y_1 \mod \liez_0$. If $b \in
\R$, $$\exp(\Ad(bY_1))(X_0+B_2) = x_1 X_1 + (y_1 - b x_1) Y_1 \mod
\liez_0 \; .$$ If $\lexp_\liek$ denotes the exponential mapping from
$\liek$ to $K$, we observe that conjugation by $\lexp_\liek(bY_1)$ for
a well chosen $b\in \R$, reduces the problem to either $x_1 = 0$ or
$y_1 = 0$. Both cases are solved by Lemma~\ref{closedsubgroup},
because $\R X_1 + \lieg_0$ and $\R Y_1 + \lieg_0$ are abelian. So,
the proof of part~3 is done.

{\it Part 4 : $\lieg_0=\langle X,Y,Z \rangle$ with $[X,Y]=Z$ and $Z$ central.}

A general derivation $\mu$ of $\lieg_0$ has the form
$$\mu(X) = x_1 X + y_1 Y + z_1 Z \; , \quad \mu(Y) = x_2 X + y_2 Y +
z_2 Z \; , \quad \mu(Z) = (x_1 + y_2)Z \; .$$ Perturbing $\mu$ with
$\Ad (-z_2 X + z_1 Y)$, $\mu$ is of the form
$$\mu \begin{pmatrix} X \\ Y \end{pmatrix} =  P \begin{pmatrix} X \\ Y 
\end{pmatrix} \; , \quad \mu(Z) = \Tr(P) \; Z \; .$$
Because $[\mu,\rho] = \rho + \Ad Y_0$, we only have two possibilities. Either 
$\rho=0$, or $\mu$ and $\rho$ have after inner perturbation and a change 
of basis in $\lieg_0$ (respecting the relations of $\lieg_0$), the form
\begin{align*}
& \mu(X) = \al X \; , \quad \mu(Y) = (\al-1) Y \; , \quad \mu(Z) = (2 \al - 1) Z \; , \\
& \rho(X) = 0 \; , \quad \rho(Y) = X \; , \quad \rho(Z) = 0 \; .
\end{align*}
Then, necessarily, $Y_0 = \lambda Z$ for some $\lambda \in \R$. We have $[X_0,\Atil] 
= \Atil + \lambda Z$, $\lieg_1=\R X_0 + \lieg_0$ and $\lieg_2 = \R(\Atil - a X_0 + B)$
for some $B \in \lieg_0$. To prove the existence of an exponentiation,
we apply $\exp(- \frac{1}{a} \Ad \Atil)$ and observe that it is an 
equivalent question to exponentiate
$\lieg_1 = \R(X_0 + \frac{1}{a} \Atil) + \lieg_0$ and $\lieg_2 = \R(X_0 + C)$ 
for some $C \in \lieg_1$. Write $C=eX+fY+gZ$. If we now replace $X_0$ by 
$X_0 + (g+2ef)Z$, we see that none of the relations above change, because $\rho(Z) = 0$, 
but only $C$ changes to $eX+fY-2ef Z$. So, we may suppose that $C$ has this last form.

Define $\lieh:=\R \Atil + \lieg_0$ and exponentiate as $H:= \R \kruisje{\rho} G_0$. Define $\te(\Atil) = - \Atil$, $\te(X)=pX$, $\te(Y)=-pY$ and $\te(Z)=-Z$, where $p=\pm 1$ will
be determined later. Then, $\te$ is an involutive automorphism of $\lieh$ that commutes 
with the derivation $\mu = \Ad X_0$ of $\lieh$. Both combine to an action $\mu$ of $\R^*$
on $H$ and we define $G=\R^* \kruisje{\mu} H$. The only problem left is the definition of 
the good closed subgroup of $G$ with tangent Lie algebra $\R(X_0 + C)$. Suppose first
that $\al \neq 1$ and $\al \neq \frac{1}{2}$. Then, we take $p=1$ and observe that $\mu$ 
is invertible on $\lieg_0^- = \langle Y,Z \rangle$ and on $[\lieg_0,\lieg_0]=\R Z$.
So, Lemma~\ref{closedsubgroup} provides us with the needed
subgroup. If $\al = 1$, we take $p=-1$ and we are done as
well. Finally, take $\al=\frac{1}{2}$ and $p=1$. Define $D=4fY$ and
observe that
$$\mu(D) - [D,C]- \frac{1}{2}[D,\mu(D)] = -2f Y + 4ef Z = \te(C) - C \; .$$
Hence, it follows from the proof of Lemma~\ref{closedsubgroup} that we can define 
the right closed subgroup of $G$.

Finally, we have to consider the case where $\rho$ is trivial. This is very much 
similar to part~1 of this proof, but simpler.
\end{proof}

Now we prove that at least one of the conditions of
Proposition~\ref{exponentiation_general} is fulfilled when
$\dim \lieg_0 \leq 3$, up to one exceptional case, that we
exponentiate explicitly \lq by hands\rq.

\begin{corollary} \label{exponentiation_casetwo}
In case~2 every real matched pair of dimension $n+1$ with $n \leq 4$ 
can be exponentiated to a matched pair of real Lie groups with at most two 
connected components.
\end{corollary}
\begin{proof}
If $\dim \lieg_0 = 1$ or $2$, then either $\lieg_0$ is abelian, or $\lieg_0$ 
has trivial center and all derivations are inner. If $\dim \lieg_0 = 3$ and 
$\lieg_0$ is of rank~3, then $\lieg_0 = \mathfrak{sl}_2$ or $\mathfrak{su}_2$. 
In both cases, every derivation is inner and the center is trivial.

When $\lieg_0$ has rank~1 or rank~0, one can always apply Proposition~\ref{exponentiation_general}.

Finally, there are only three real non-isomorphic 3-dimensional $\lieg_0$ 
of rank~2 defined respectively by : $ a)\ [H,X]=X,\ [H,Y] = \al Y,
\ [X,Y]=0\ (\al\in \R);\quad b)\ [H,X]=X+Y,\ [H,Y]=Y,\ [X,Y]=0;\quad 
c)\ [H,X]=rX + Y,\ [H,Y]=-X + r Y,\ [X,Y]=0\ (r\in \R)$ (see \cite{FH}).

All these cases (except a), $\al=1$, which we will study separately), can be 
treated in a similar way. Namely, a general derivation $\mu$ of $\lieg_0$ 
has the following form: 

$a),\ \al \neq 1:\ \mu(H) = aX+bY \; , 
\quad \mu(X) = cX \; , \quad \mu(Y) = dY,$ and it is inner if $d = \al c$.

$b)\ \mu(H) = xX + yY \; , \quad \mu(X) = aX+bY \; , \quad \mu(Y) = aY,$
and it is inner if $a=b$.

$c)\ \mu(H) = xX+yY \; , \quad \mu(X) = aX+bY \; , \quad \mu(Y) = -b X + aY,$ 
and it is inner if $a = rb$

Then, since $\mu$ and $\rho$ are derivations of $\lieg_0$, we observe that in
all cases above $[\mu,\rho]$ is inner. Hence, $\rho = [\mu,\rho] - \Ad Y_0$ is 
inner. Also, the center of $\lieg_0$ is always trivial.

At last, we study separately $\lieg_0$ defined by $[H,X]=X$, $[H,Y] = Y$ 
and $[X,Y]=0$. A general derivation $\mu$ of $\lieg_0$ has the form
$$
\mu(H) = xX + y Y \; , \quad \mu(X) = aX + bY \; , \quad \mu(Y) = cX
+ dY \; ,
$$ 
which is inner if $b=c=0$ and $a=d$. Because $[\mu,\rho] = \rho + \Ad Y_0$, 
we conclude that there exists
a $\lambda \in \R$ such that, on $\langle X,Y \rangle$, $\mu \rho -
\rho \mu = \rho + \lambda \io$. It is easy to see that either $\rho =
- \lambda \io$, in which case $\rho$ is inner and the exponentiation exists, 
or that we can change the
basis $X,Y$ and perturb $\mu$ and $\rho$ innerly such that
\begin{align*}
& \mu(H) = 0 \; , \quad \mu(X) = X \; , \quad \mu(Y) = 0, \\
& \rho(H) = 0 \; , \quad \rho(X) = 0 \; , \quad \rho(Y) = X \; .
\end{align*}
Our matched pair lives now in the Lie algebra $\lieg$ with generators
$X_0,\Atil,X,Y,Z$ with relations $\Ad X_0 = \mu$, $\Ad \Atil = \rho$
on $\lieg_0$, $[X_0,\Atil]=\Atil$. Observe that $Y_0$ disappeared
because, after the necessary inner perturbations of $\mu$ and $\rho$,
we arrived at $[\mu,\rho]=\rho$. The Lie subalgebras $\lieg_1$ and
$\lieg_2$ are $\langle X_0,H,X,Y \rangle$ and $\R (\Atil - a X_0 +
B)$, respectively, where $B \in \lieg_0$ is arbitrary. We exponentiate
$\R \Atil + \lieg_0$ on the space $\R^4$ with product
$$
(a,h,x,y)(a',h',x',y') = (a+a',h+h',x + \exp(h) x' + a \exp(h) y', y
+ \exp(h)y')
$$ and we denote this Lie group by $H$. Write $B=\al H + \be X + 
\ga Y$. Suppose first that $\al \neq 0$. Then, we make $\R^*$ act
on $H$ by the automorphisms $\mu_s(a,h,x,y)  = (sa,h,|s|x,\Sgn(s)y)$. Take $G =
\R^* \ltimes H$ in which we consider the closed subgroup $G_1$ of
elements $(s,0,h,x,y)$, $s \in \R^*, h,x,y \in \R$. Next, we have to find a 
good (as in the proof of Proposition~\ref{exponentiation_general}) closed 
subgroup $G_2$ with tangent Lie algebra $\R(\Atil - a X_0 +
B)$. Conjugating with the element $(1,-\frac{1}{a},0,0,0)$, we have to find 
a good closed subgroup of $\R^* \ltimes G_0$ with tangent Lie algebra 
$\R(X_0 + \al'H + \be'X + \ga'Y)$,
where $\al'=-\frac{\al}{a} \neq 0$. If $\al' \neq -1$, such a subgroup is given by
$$
\{ (s,\al' \log |s| , \frac{\be'}{\al'+1}(|s|^{\al'+1} - 1) ,
\frac{\ga'}{\al'}(\Sgn(s) |s|^{\al'} - 1)) \mid s \in \R^* \} \; .
$$
If $\al'=-1$, such a subgroup is given by
$$
\{ (s,-\log |s|, \be' \log |s| , \ga'(1- \frac{1}{s})) \mid s \in \R^* \} \; .
$$ 
One can easily check that $(G_1,G_2)$ is indeed a matched pair of Lie groups.

If next, $\al = 0$, we consider the action of $\R^*$ by automorphisms of $H$ given by
$\mu_s(a,h,x,y)  = (sa,h,sx,y)$. Take $G =
\R^* \ltimes H$ in which we consider the closed subgroup $G_1$ of
elements $(s,0,h,x,y)$, $s \in \R^*, h,x,y \in \R$.
Conjugating with the element $(1,-\frac{1}{a},0,0,0)$, we have to find a good 
closed subgroup of $\R^* \ltimes G_0$ with tangent
Lie algebra $\R(X_0 + \be'X + \ga' Y)$ and this is given by
$$
\{ (s,0,\be' (s-1) , \ga' \log |s|) \mid s \in \R^* \} \; .
$$
Again, we arrive at a matched pair of Lie groups.
\end{proof}

As it was promised, we now present a matched pair of dimension $5+1$ which has no exponentiation.
\begin{example}
Consider the $6$-dimensional Lie algebra $\lieg$ with generators 
\newline $A, X_0, X_1, X_2, Y, Z$ and relations
\begin{align*}
& [X_0,A]=A+Y \; , \quad [X_0,X_1] = \frac{1}{2} X_1 + X_2 \; , \quad [X_0,X_2] = \frac{1}{2} X_2 \; , \\ & [X_0,Y] = Y \; , \quad [X_0,Z] = \frac{3}{2} Z \; , \\
& [A,X_1] = [A,Y] = [A,Z] = 0 \; , \quad [A,X_2] = Z \; , \quad [X_1,X_2]=Y \; , \\ & [X_1,Y] = Z \; , \quad [X_1,Z]=[X_2,Y]=[X_2,Z] = [Y,Z] = 0 \; .
\end{align*}
Defining $\lieg_1 = \langle X_0,X_1,X_2,Y,Z \rangle$ and $\lieg_2 = \R(A + X_0)$, we have a matched 
pair of Lie algebras. There does not exist a matched pair of Lie groups
which exponentiates the given matched pair of Lie algebras.
\end{example}
\begin{proof}
It is easy to check that $\lieg$ is indeed a Lie algebra: it is clear that $\lieg_0:=\langle X_1,X_2,Y,Z 
\rangle$ is a Lie algebra. Next, we write $\rho = \Ad A$ and $\mu = \Ad X_0$
on $\lieg_0$ and it is clear that $\mu$ and $\rho$ are derivations of $\lieg_0$ satisfying $[\mu,\rho]=
\rho + \Ad Y$. The connected, simply connected Lie group $G_0$
of $\lieg_0$ has $\R^4$ as an underlying space, with product
$$
(x_1,x_2,y,z) (x_1',x_2',y',z') = (x_1+x_1',x_2+x_2',y+y'+x_1x_2',z+z'+\frac{x_1^2}{2} x_2' + x_1y') \; .
$$
The derivation $\rho$ gives rise to an action $(\rho_a)$ of $\R$ on $G_0$ given by $\rho_a(x_1,x_2,y,z) = (x_1,x_2,y+ax_2,z)$ and we can define the Lie group
$H:= \R \kruisje{\rho} G_0$ on the space $\R^5$. Finally, $\mu$ gives rise to an action $(\mu_x)$ of $\R$ on $H$ given by
\begin{align*}
\mu_x(a,x_1,x_2,y,z) = &(\exp(x) a, \exp(\frac{1}{2}x) x_1, \exp(\ha{1}x)(x_2+xx_1), \\ & \exp(x)(y+xa + \ha{1} x x_1^2), \exp(\ha{3} x)(z+xax_1 + \ha{6} x x_1^3)) \; .
\end{align*}
Defining $G=\R \kruisje{\mu} H$ on the space $\R^6$, we obtain the connected, simply connected 
Lie group of $\lieg$. One observes easily that the center of $G$ is trivial.
Hence, $G$ is the only connected Lie group with Lie algebra $\lieg$.

We claim the following: any automorphism $\al$ of $G$ leaves the closed normal subgroup $G_0$  
invariant and modulo an inner perturbation by $\Ad g$ $(g \in G)$, we have $\al = \io$ on $G/G_0$. Indeed, let $\al$ be an automorphism of $G$. Denote $\te = d\al$, the corresponding 
automorphism of $\lieg$. It is straightforward to check that $\lieg$ has only
two ideals of dimension $4$: $\lieg_0$ and $\langle A,X_2,Y,Z \rangle$. As a Lie algebra, 
the first one has rank~2 and the second one rank~1. Because $\te(\lieg_0)$ is
an ideal of dimension $4$ of $\lieg$ isomorphic with $\lieg_0$, we get $\te(\lieg_0)=\lieg_0$. Then also 
$\al(G_0)=G_0$. Because $[X_0,A]= A \mod \lieg_0$, a perturbation of $\al$ by $\Ad (x,a,0,0,0,0)$
for certain $x,a \in \R$, leaves two possibilities for $\te$: $\te = \io \mod \lieg_0$ or 
$\te(X_0) = X_0 \mod \lieg_0$ and $\te(A) = - A \mod \lieg_0$. We have to prove that the
second option is impossible. Hence, suppose that $\te(X_0) = X_0 + C$ and $\te(A) = -A + D$ 
with $C,D \in \lieg_0$. If we arrive at a contradiction, then our claim is proved.
Because $\te$ is an automorphism of $\lieg_0$ it has the following form
\begin{align*}
& \te(X_1) = b X_1 + g X_2 + h Y + k Z \; , \quad \te(X_2) = d X_2 + e Y + f Z \; , \\
& \te(Y) = bd Y + be Z \; , \quad \te(Z) = b^2 d Z \; ,
\end{align*}
for $e,f,g,h,k \in \R$ and $b,d \in \R^*$. Because $\te$ is an automorphism of $\lieg$,
we get $\te \rho=(- \rho + \Ad D) \te$. Suppose that $D=x_1X_1 + x_2 X_2 +y Y +z Z$.
Verifying the equality on $Y$ gives $x_1 b d Z = 0$ and, hence,
$x_1=0$. Next, we verify on $X_2$ and conclude that $b^2 d Z= - dZ$,
which yields the required contradiction $b^2=-1$.

Suppose now that $\cG$ is a Lie group with Lie algebra $\lieg$ and with two closed subgroups 
$\cG_1,\cG_2$ whose tangent Lie algebras are $\lieg_1,\lieg_2$, respectively
and such that $(\cG_1,\cG_2)$ is a matched pair of Lie groups. Because $G$ is the only 
connected Lie group with Lie algebra $\lieg$, we can identify $G$ and $\cG^{(e)}$, the
connected component of $e$ in $\cG$.

Because any automorphism of $G$ leaves $G_0$ invariant, we find that $G_0$ is a normal 
subgroup of $\cG$. We can naturally identify $\bigl( \cG / G_0 \bigr)^{(e)}$
and $G/G_0$. Define $H = \cG / G_0$. We claim that for any $\eta \in H / H^{(e)}$, 
there exists 
a unique representative $u(\eta) \in H$ such that $\Ad u(\eta) = \io$ on $H^{(e)}$.
Denote by $\pi_1$ the quotient map from $\cG$ to $H$ and by $\pi_2$ the quotient map 
from $H$ to $H/H^{(e)}$. Let $\eta \in H / H^{(e)}$. Take $g \in G$ such that 
$\pi_2(\pi_1(g)) = \eta$.
Then, $\Ad g$ defines an automorphism of $\cG^{(e)}=G$. Using our claim, we can take an $h \in G$, 
such that $\Ad (hg)$ is trivial on $G/G_0$. This means that $\Ad \pi_1(hg)$ is
trivial on $H^{(e)}$. But $\pi_2(\pi_1(hg)) = \eta$. So, we have proven the existence of the 
required representative. The uniqueness is trivial because $H^{(e)} = G/G_0$ has
trivial center. Then, the map $\Psi:H/H^{(e)} \oplus H^{(e)} \recht H : (\eta,g) \mapsto u(\eta) g$ 
is an isomorphism of Lie groups. (Recall that $H/H^{(e)}$ has dimension zero and
the discrete topology.)

The only closed subgroup of $G$ with tangent Lie algebra $\lieg_1$ is
the group consisting of the elements $(x,0,x_1,x_2,y,z)$. Hence, this
last subgroup coincides with $\cG_1 \cap G$. On the other hand, the
only closed subgroup of $G$ with tangent Lie algebra $\lieg_2$
consists of the elements $(x,\exp(x)-1,0,0,x \exp(x),0)$. Hence, this
last subgroup coincides with $\cG_2 \cap G$. Identifying $H^{(e)}$
with the group $K:=\R^2$ with product $(x,a)(x',a') = (x+x',a +
\exp(x) a')$, we get that $K_1:=\pi_1(\cG_1) \cap H^{(e)} = \{(x,0) \mid x
\in \R \}$ and $K_2:=\pi_1(\cG_2) \cap H^{(e)} = \{(x,\exp(x)-1) \mid x \in
\R \}$. Combining this with the fact that $H$ is isomorphic with
$H/H^{(e)} \oplus H^{(e)}$ and with the fact that the only elements of
$K$ normalizing $K_1$, resp.\ $K_2$, belong to $K_1$, resp.\ $K_2$, we
conclude that $\pi_1(\cG_i) = Q_i \oplus K_i$ for some subgroups $Q_i
\subset H/H^{(e)}$ and $i=1,2$. Because $\cG_1 \cG_2$ is dense in
$\cG$, it follows that $K_1 K_2$ should be dense in $K$. This is
clearly not the case.
\end{proof}

\centerline{\bf Case~3} 

Now we have $\lieg_0 := \Ker \chi \cap \Ker \chi \be$ and this is still an
ideal of $\lieg$. We can take $X,Y \in \lieg_1$ such that
$$\chi(X) = 1 \; , \quad \chi(Y) = 0 \; , \quad \chi(\be(X)) = 0 \; ,
\quad \chi(\be(Y))= 1 \; .$$
Using Equation~\eqref{eq2}, we get $\chi(\be([X,Y]))=-1$ and
$\chi([X,Y]) = 0$. Hence, $[X,Y]= - Y  \mod \lieg_0$. Because
$\chi(\be(X))=0$ and $\chi(\be(Y))=1$, we get
$$\be(X) = aY \mod \lieg_0\; , \quad \be(Y) = X + bY \mod
\lieg_0\; , \quad a,b \in k \; .$$
Checking Equation~\eqref{eq1}, we arrive at $b=0$. Since the quotient
Lie algebra $\lieg / \lieg_0$ is 3-dimensional and of rank 3, it is
isomorphic to $\mathfrak{sl}_2(k)$ \cite{FH} (if $k=\R$, one must analyse
also $\mathfrak{su}_2(\R)$, but it has no 2-dimensional Lie subalgebras).
Then, by the
Levy-Maltsev theorem \cite{NS}, Chapter X, we can find a Lie subalgebra
$\tilde\lieg$ of $\lieg$
isomorphic to $\mathfrak{sl}_2(k)$ and such that $\lieg = \tilde\lieg
\oplus \lieg_0$ as vector spaces. This means that we are always in the
following situation: $\sll_2(k)$, with the generators $A,X,Y$
satisfying
$$[X,A] = aY + A \; , \quad [Y,A]= X \; , \quad [X,Y]=-Y \; ,$$
is
represented by derivations of $\lieg_0$ and in the semi-direct product
$\lieg:= \sll_2(k) \ltimes \lieg_0$ we have the matched pair $\lieg_1
= \langle X,Y \rangle + \lieg_0$, $\lieg_2 = k(A + Z)$ for some
element $Z \in \lieg_0$.

We first analyse the case $k=\R$. If we replace $Y$ by $rY$ and $A$ by 
$(1/r) A$ for $r \neq 0$, then $a$
changes to $a/r^2$. So, we only have to consider three cases: one with
$a>0$, one with $a=0$ and one with $a>0$.

In terms of the standard generators $H=e_{11}-e_{22}$, $K=e_{12}$ and
$L=e_{21}$, we realize the relations above by putting $X=\frac{1}{2} H$,
$Y=L$ and $A=-2a L - \frac{1}{2} K$. This means that we consider the
matched pair $\lieg_1 = \langle H,L \rangle + \lieg_0$, $\lieg_2 = \R (K + 4aL +
Z)$ for some $Z \in \lieg_0$.

\begin{proposition} \label{exponentiation_casethree}
If $k = \R$ and $a \leq 0$, the matched pair has an
exponentiation. When $a < 0$, $G,G_1,G_2$ can be taken connected. When
$a=0$, we need two connected components for $G_1$.

If $k = \R$ and $a > 0$, the matched pair has an exponentiation in at
least the following cases:
\begin{enumerate}
\item $\lieg_0$ is abelian.
\item $\dim \lieg_0 = 2$.
\end{enumerate}
We can take a matched pair of Lie groups where $G,G_1$ have at most two 
connected componenents and $G_2$ has at most four connected components.

In particular, for $k = \R$ 
the exponentiation exists for all matched pairs of dimension $n+1$, $n \leq 4$.
\end{proposition}
\begin{proof}
We start off with the case $a \leq 0$ and we
first suppose that $\lieg_0=\{0\}$.
An exponentiation of the case $a=0$ is given by $F=\SL_2(\R)$ with subgroups
\begin{equation} \label{bbbbb}
F_1 : = \{ \begin{pmatrix} a & 0 \\ x & \frac{1}{a}  \end{pmatrix}
\mid a \neq 0, x \in \R \} \; , \quad F_2 : = \{
\begin{pmatrix} 1 & b \\ 0 & 1 \end{pmatrix} \mid b \in \R \}
\; .
\end{equation}
An exponentiation of the case $a=-\frac{1}{4}$ (as we saw above, it is
sufficient to consider one value of $a < 0$) is given by $F=\SL_2(\R)$
with subgroups
$$F_1 : = \{ \begin{pmatrix} a & 0 \\ x & \frac{1}{a}  \end{pmatrix}
\mid a > 0, x \in \R \} \; , \quad F_2 : = \{
\begin{pmatrix} \cos t & \sin t \\ - \sin t  & \cos t \end{pmatrix}
\mid t \in \R \}
\; .$$
For this last case, there is another exponentiation which is at least
as important. We observe that $F=F_1F_2$ and the multiplication map is
a diffeomorphism of $F_1 \times F_2$ onto $F$. This means that the
mutual actions $\al_g(s)$ and $\be_s(g)$, $g \in F_1, s  \in F_2$ are
everywhere defined and smooth. Identifying $F_2$ with $\T$, we observe
that for all $g \in F_1$, $\al_g$ is a diffeomorphism of $\T$
satisfying $\al_g(1) = 1$. Hence, there exists, for every $g \in F_1$ a 
unique diffeomorphism $\tilde{\al}_g$ of $\R$, such that $\tilde{\al}_g(0) = 
0$ and $p(\tilde{\al}_g(t)) = \al_g(p(t))$ for all $t \in
\R$, where $p(t) = \cos t + i \sin t$. Defining $\tilde{\be}_t(g) =
\be_{p(t)}(g)$, we obtain a matched pair $(F_1,\R)$, in which both
actions are everywhere defined and smooth. Its corresponding big Lie
group $F_{\text{\rm sc}}$ is the connected, simply connected Lie group 
of $\sll_2(\R)$.

Suppose now that $\lieg_0$ is arbitrary and $Z \in \lieg_0$. First,
take $a = 0$. Because any finite-dimensional representation of
$\sll_2(\R)$ can be exponentiated to $\SL_2(\R)$ (although $\SL_2(\R)$
is not simply connected, see e.g.\ \cite{bourbaki}, Chapitre VIII,
par.\ 1, Th\'eor\`eme 2), we get an action $\mu$ of $F:=\SL_2(\R)$ with
automorphisms of $G_0$, the connected, simply connected Lie group of
$\lieg_0$. We define $G:=F \kruisje{\mu} G_0$ and we denote by
$\lexp_\lieg$ its exponential mapping. Because
$$\lexp_\lieg(b(K+Z)) = \Bigl( \begin{pmatrix} 1 & b \\ 0 & 1
\end{pmatrix}, \ldots \Bigr) \; ,$$ it is clear that, defining
$$G_1 := \{(g,k) \mid g \in F_1, k \in G_0 \} \; , \quad G_2 : =
\{\lexp_\lieg(b(K+Z)) \mid b \in \R \} \; ,$$
where $F_1$ is as in Equation~\eqref{bbbbb}, we get the required
matched pair of Lie groups.

When $a=-\frac{1}{4}$, we proceed similarly, but now with the
simply connected Lie group $F_{\text{\rm sc}}$ of $\sll_2(\R)$. We get an
action $\mu$ of $F_{\text{\rm sc}}$ on $G_0$ by automorphisms and we define
$G:= F_{\text{\rm sc}} \kruisje{\mu} G_0$. As we explained above, we can
find in $F_{\text{\rm sc}}$ a matched pair $(F_1,F_2)$ with tangent Lie
algebras $\langle X,Y \rangle$ and $\R (K-L)$, such that $F_2$ can be
identified with $\R$. If $\lexp_\lieg$ denotes the exponential mapping
of $\lieg$, it follows that $$\lexp_\lieg(t(K-L+Z)) = (t,\ldots) \in F_2
\times G_0 \subset F_{\text{\rm sc}} \kruisje{\mu} G_0 \; .$$ So, we can
define in the same way as above the required matched pair of Lie
groups. Observe that this argument would not work with $\SL_2(\R)$
instead of $F_{\text{\rm sc}}$, because then $\lexp_\lieg(2\pi n
(K-L+Z)) \in G_0 \subset G_1$ when $n \in \Z$ and hence, we no longer
have $G_1 \cap G_2 = \{e\}$.

Next, we turn to $a > 0$ and we choose the value $a =
\frac{1}{4}$. Then, $\lieg_2 = \R (K+L+Z)$ for some $Z \in \lieg_0$.
When $\lieg_0=\{0\}$, we can exponentiate as follows.
Write $F= \{ T \in M_2(\R) \mid \operatorname{det} T = \pm 1 \}$. Define
\begin{equation} \label{ccc}
F_1 : = \{ \begin{pmatrix} a & 0 \\ x & \frac{s}{a}  \end{pmatrix}
\mid a > 0, s = \pm 1 , x \in \R \} \; , \quad F_2 : = \{
\begin{pmatrix} b & c \\ c & b \end{pmatrix} \mid b^2 - c^2 = \pm 1 \}
\; .
\end{equation}
It is an easy exercise to check that we indeed get a matched pair of
Lie groups. For general $\lieg_0$, we want to proceed as in the case
$a=0$. We first get an action $\mu$ of $\SL_2(\R)$ on $G_0$. We now run into
the same kind of problems as in the proof of
Proposition~\ref{exponentiation_general}. We should first extend the
action $\mu$ to an action of $F$, by adding an involutive automorphism
of $G_0$ corresponding to the action of $\left(\begin{smallmatrix} 1 & 0 \\ 0 & -1
\end{smallmatrix}\right)$ on $G_0$ and next, we should find the good closed
subgroup with tangent Lie algebra $\R (K+L+Z)$ and with elements whose first 
components are precisely the matrices $\left(\begin{smallmatrix} b & c \\ c & b
\end{smallmatrix}\right)$ with $b^2 - c^2 = \pm 1$.

First, take $\lieg_0$ abelian. Write $\mu_H, \mu_K$ and $\mu_L$ for
the derivations of $\lieg_0$ corresponding to the generators $H,K$ and
$L$ of $\mathfrak{sl}_2(\R)$. From \cite{bourbaki} (Chapitre VIII,
par.\ 1, no.\ 2, Corollaire), we can write
$$\lieg_0 = \bigoplus_{-N \leq n \leq N} E_n \; ,$$
where $N \in \N$, $n$ only takes values in $\Z$ and where
$$\mu_H(X) = n X \;\; \text{for} \;\; X \in E_n \; , \quad \mu_K(E_n)
\subset E_{n+2} \; , \quad \mu_L(E_n) \subset E_{n-2}$$
with $E_m = \{0\}$ if $m < -N$ or $m > N$. We exponentiate $\mu$ to a
homomorphism $\mu : \SL_2(\R) \recht \operatorname{GL}(\lieg_0)$.
If we define the involution $\te$ of $\lieg_0$ by putting $\te(X) = X$
if $X \in E_n$, $n = 0,1 \mmod 4$ and $\te(X) = -X$ if $X \in E_n$ and
$n = 2,3 \mmod 4$, then, $\te$ commutes with $\mu_H$ and anti-commutes
with $\mu_K$ and $\mu_L$. So, we extend $\mu$ to the group $F$ of
matrices with determinant $\pm 1$ by writing $\mu
\left(\begin{smallmatrix} 1 & 0 \\ 0 & -1 \end{smallmatrix}\right) =
\te$. We define, on the space $F \times \lieg_0$, the Lie group $G$
with product $(P,X)(Q,Y) = (PQ,X+\mu(P)Y)$ for $P,Q \in F$ and $X,Y
\in \lieg_0$. Define $G_1$ consisting of the pairs $(P,X)$ for $P \in
F_1$ and $X \in \lieg_0$, with $F_1$ as in
Equation~\eqref{ccc}. Finally, we have to find the good closed
subgroup with tangent Lie algebra $\R (K + L + Z)$ for $Z \in
\lieg_0$. This procedure is described in great detail in the proof of
Proposition~\ref{exponentiation_general}.
After a well chosen conjugation, we have to find a
closed subgroup of $G$ with tangent Lie algebra $\R (H + C)$, for some
$C \in \lieg_0$, consisting of the elements
$$\{ \Bigl( \begin{pmatrix} a & 0 \\ 0 & b
  \end{pmatrix}, v(a,b) \Bigr) \mid ab=\pm 1 \} \; .$$
As we see from the proof of Lemma~\ref{closedsubgroup}, the only point
is to define $v(-1,-1)$ and $v(1,-1)$. It follows from
\cite{bourbaki}, Chapitre VIII, par.\ 1, no.\ 5, Corollaire, that
$\mu\left(\begin{smallmatrix} -1 & 0 \\ 0 & -1
  \end{smallmatrix}\right) X = (-1)^n X$ for $X \in E_n$. Denote by
$\rho$ this involution of $\lieg_0$. The proof of
Lemma~\ref{closedsubgroup} suggests us to define $v(-1,-1) =
\mu_H^{-1}(\rho(C) - C)$ and $v(1,-1) = \mu_H^{-1}(\te(C) - C)$. 
Both are well-defined, because
$$\rho(C) - C \in \bigoplus_{n = 1 \operatorname{mod} 2} E_n \; , 
\quad \te(C) - C \in \bigoplus_{n = 2,3 \operatorname{mod} 4} E_n$$ 
and $\mu_H$ is invertible on both subspaces.
Also, $v(-1,-1)$ and $v(1,-1)$ should be compatible, in
the sense that
$$v(-1,-1) + \rho(v(1,-1)) = v(1,-1) + \te(v(-1,-1)) \; ,$$
or equivalently
$$(\rho-\io) \mu_H^{-1}(\te - \io)(C) = (\te - \io) \mu_H^{-1} (\rho -
\io)(C) \; ,$$
which is the case because the factors commute.

To finish the proof of the proposition, it remains to consider non-abelian
2- dimensional $\lieg_0$. So,
$\lieg_0$ has generators $X,Y$ with relation $[X,Y]=Y$. The Lie
algebra of derivations of $\lieg_0$ is now isomorphic to
$\lieg_0$. Any homomorphism of $\sll_2(\R)$ into $\lieg_0$ must be
trivial, because its kernel is a non-zero ideal of $\sll_2(\R)$. So,
any action of $\sll_2(\R)$ on $\lieg_0$ is necessarily trivial and we
take $G:= F \oplus G_0$, where $G_0$ is the $ax+b$-group. We take
$G_1=F_1 \oplus G_0$ and for any $Z \in \lieg_0$ we can define
$$G_2 = \{ \Bigl( \begin{pmatrix} b & c \\ c & b
  \end{pmatrix} , \lexp_{\lieg_0}((\log |b+c|) Z) \Bigr)
\mid b^2 - c^2 = \pm 1 \} \; .$$
This again provides us with a matched pair of Lie groups.
\end{proof}

\begin{remark}
For the case $\lieg_0 = \{0\}$, we will give more connected
exponentiations below. In the proof of the previous proposition they
are not so interesting, because we will then rather have
$G=\PSL_2(\R)$ and not every representation of $\SL_2(\R)$ factors
through $\PSL_2(\R)$. Hence, we cannot make the right semi-direct
products with $\PSL_2(\R)$ acting.
\end{remark}

\begin{remark}
In case~3, for $k=\R$ and $n \geq 5$, there are indications that there
again exist matched pairs of Lie algebras that cannot be
exponentiated. Their explicit description remains however open.
\end{remark}

Next, we analyze the case $k=\C$. First, let us note, that now there 
are only two non-isomorphic cases: with $a=0$ and with $a\neq 0$.

\begin{proposition}
In case~3, any matched pair of complex Lie algebras can be 
exponentiated to a matched pair of connected complex Lie 
groups if $n \leq 3$.
\end{proposition}
\begin{proof}
First, consider the case $\lieg_0 = \{0\}$. If $a=0$, we proceed
exactly as in the case of $k=\R$ and just replace $\R$
by $\C$. If $a \neq 0$, we consider $\lieg_1 = \langle H,L \rangle$
and $\lieg_2 = \C (K-L)$. Define $F = \PSL_2(\C)$ and define
\begin{align}
F_1 : &= \{ \begin{pmatrix} a & 0 \\ x & \frac{1}{a}  \end{pmatrix}
\mod \{\pm 1 \} \mid a \neq 0, x \in \C \} \; , \label{fff} \\ F_2 : &= \{
\begin{pmatrix} \cos z & \sin z \\ - \sin z  & \cos z \end{pmatrix}
\mod \{\pm 1\} \mid z \in \C \} \; . \notag
\end{align}
Some care is needed in checking that we do get a matched pair of Lie
groups. Writing the product of an element in $F_1$ and an element in
$F_2$, we have to find a unique solution in $F_1,F_2$ of the equation
$$\begin{pmatrix} a \cos z & a \sin z \\ x \cos z - \frac{1}{a} \sin z
  & x \sin z + \frac{1}{a} \cos z \end{pmatrix} = \begin{pmatrix} u &
  v \\ w & r \end{pmatrix} \mod \{\pm 1\}$$
whenever $ur - vw = 1$. Given $u,v,w,r \in \C$ with $ur-vw = 1$, we
proceed as follows: choose $a \in \C$ such that $a^2 = u^2 + v^2$ and
define $\cos z = \frac{u}{a}$, $\sin z = \frac{v}{a}$ and $x =
\frac{uw+vr}{a}$. Then, the required equation holds. If we choose the
other square root of $u^2 + v^2$, then $a,x,\cos z$ and $\sin z$ change
sign and hence, their projections $\mmod \{\pm 1\}$ do not
change. Because clearly $F_1 \cap F_2 = \{e\}$, we have a matched pair
of Lie groups.

If next, $\lieg_0 = \C$, the action of $\sll_2(\C)$ on $\lieg_0$ is
necessarily trivial. Our matched pair has
the form $\lieg = \sll_2(\C) \oplus \C$, $\lieg_1 = \langle H,L,Z
\rangle$, $\lieg_2 = \C (K+4aL+ \lambda Z)$, where $Z$ is the generator
of $\lieg_0=\C$ and $\lambda \in \C$. If $\lambda = 0$, it is clear
how to exponentiate, just adding a copy of $\C$ to $F$ and $F_1$
above. If $\lambda \neq 0$, we change the generator $Z$ and we may
suppose that $\lambda = 1$. If $a=0$, exponentiation is again easy. If
we take $a=-\frac{1}{4}$, we denote by $\T$ the complex torus,
consisting of the pairs $(\cos z ,\sin z) \in \C^2$, we define $G =
\PSL_2(\C) \oplus \T$ with subgroups $G_1:=F_1 \oplus \T$ (with $F_1$
as in Equation~\eqref{fff}) and
$$G_2 = \{ \Bigl(\begin{pmatrix} \cos z & \sin z \\ - \sin z  & \cos z
\end{pmatrix}, (\cos z,\sin z) \Bigr) \mid z \in \C \} \; .$$
We indeed have a matched pair of Lie groups.
\end{proof}

We conclude with the promised counterexample in dimension $4+1$.

\begin{example}
Define $\lieg$ to be the complex Lie algebra $\lieg:=\sll_2(\C) \oplus
\lieg_0$, where $\lieg_0$ has the generators $X,Y$ satisfying
$[X,Y]=Y$. Use the canonical generators $H,K,L$ of $\sll_2(\C)$ and define the matched pair
$$\lieg_1:=\langle H,L,X,Y \rangle \; , \quad \lieg_2 := \C (K - L +
Y) \; .$$
There does not exist an exponentiation of this matched pair of Lie
algebras.
\end{example}
\begin{proof}
The connected, simply connected Lie group of $\lieg$ is given by
$\SL_2(\C) \oplus G_0$, where $G_0$ lives on the space $\C^2$ with
product $(x,y)(x',y') = (x+x',y + \exp(x)y')$. Its center consists of the elements
$(\pm 1,2\pi n,0)$, where $n \in \Z$. Hence, the
only connected Lie groups with Lie algebra $\lieg$ are $G:=\SL_2(\C)
\oplus G_0/H_N$ and $G':=\PSL_2(\C) \oplus G_0/H_N$, where $H_N$ consists of the elements $(2\pi nN,0)$, $n \in \Z$.
If we take $G:=\SL_2(\C) \oplus G_0$,
the connected closed subgroup of $G$ with tangent Lie algebra $\lieg_1$
consists of the elements
$$\Bigl( \begin{pmatrix} a & 0 \\ z & \frac{1}{a}  \end{pmatrix}, x,y
\Bigr) \; , \quad a \neq 0, z,x,y \in \C \; .$$
The connected closed subgroup of $G$ with tangent Lie algebra
$\lieg_2$ consists of the elements
$$\Bigl(\begin{pmatrix} \cos z & \sin z \\ - \sin z  & \cos z
\end{pmatrix}, 0,z \Bigr) \mid z \in \C \; .$$
The intersection of both subgroups is non-trivial, because it contains
the elements $(1,0,2 \pi n)$ with $n \in \Z$. This intersection is not annihilated by any of the central subgroups of $G$.
So, with the same kind of
reasoning as in Example~\ref{ex2.4}, we conclude that the matched pair
cannot be exponentiated.
\end{proof}

\section{Matched pairs of real Lie algebras of dimension $1+1$ and
$2+1$ and their exponentiation}

Next we classify, up to isomorphism, all matched pairs of real Lie algebras of
dimension $1+1$ and $2+1$ and compute explicitly their exponentiation.

From now on, all Lie algebras and Lie groups are understood to be real.

\begin{theorem} \label{prop3.2}
In dimension $1+1$, there exist, up to isomorphism, the following
non-isomorphic matched pairs of Lie algebras. We choose a generator $X$ for $\lieg_1$.
\begin{enumerate}
\item $\chi=0$ and $\be=0$.
\item $\chi=0$ and $\be(X)=X$.
\item $\chi(X)=1$ and $\be=0$.
\item $\chi(X)=1$ and $\be(X)=X$.
\end{enumerate}
In dimension $2+1$, there exist, up to isomorphism, the following
non-isomorphic matched pairs of Lie algebras. We choose generators
$X,Y$ for $\lieg_1$. \renewcommand{\theenumii}{\arabic{enumi}.\arabic{enumii}}\renewcommand{\labelenumii}{\theenumii.}
\begin{enumerate}
\item $\chi=0$ and $\be=0$.
\begin{enumerate}
\item $[X,Y]=0$.
\item $[X,Y]=Y$.
\end{enumerate}
\item $\chi=0$ and $\be \neq 0$:
\begin{enumerate}
\item $[X,Y]=0$, $\be(X) = X$, $\be(Y) = rY$, $-1 \leq r \leq 1$.
\item $[X,Y]=0$, $\be(X) = X+Y$, $\be(Y) = Y$.
\item $[X,Y]=0$, $\be(X) = Y$, $\be(Y)=0$.
\item $[X,Y]=Y$, $\be(X) = Y$, $\be(Y)=0$.
\item $[X,Y]=Y$, $\be(X) = 0$, $\be(Y) = Y$.
\end{enumerate}
\item $\chi \neq 0$ and $\be = 0$: $[X,Y]= aY$, $\chi(X)=1$,
  $\chi(Y)=0$, $a \in \R$.
\item $\chi \neq 0$ and $\be \neq 0$: $\chi(X)=1$, $\chi(Y)=0$ and
\begin{enumerate}
\item $[X,Y]=dY$, $\be(X) = X+ bY$, $\be(Y) =d Y$, either $d=1$ and
  $b \in \R$, or $d \neq 1$ and $b=0$.
\item $[X,Y]=dY$, $\be(X) = Y$, $\be(Y) = 0$, $d \in \R$.
\item $[X,Y]=-Y$, $\be(X) = aY$, $\be(Y) = X$, $a=1,0,-1$.
\end{enumerate}
\end{enumerate}
Every matched pair above can be
exponentiated to a matched pair of Lie groups, having at most 2 connected components.
\end{theorem}
\begin{proof}
In dimension 1+1 the classification is obvious. If either $\be=0$ or $\chi=0$,
an exponentiation can be given using the semi-direct product
of the corresponding connected simply connected Lie groups. In the remaining case,
the exponentiation
was explicitly described in Remark~\ref{remnieuw}.

In dimension $2+1$, it is again natural to separate the cases 1, 2, 3 and 4.
In case~1, the classification follows from the classification of 2-dimensional Lie groups.
In case~2, we observe that $\be([X,Y]) =[\be(X), Y] + [X,\be(Y)]$, i.e., $\be$ is an action.
If $[X,Y]=0$, any linear map $\be$ defines an
action. We have either $\be$ diagonalizable (case 2.1), either $\be$ not
diagonalizable and not nilpotent (case 2.2), or $\be$ nilpotent (case
2.3). Multiplying $A$ by a scalar, one can scale $\be$. Hence, cases 2.2 and 2.3 cover all the
non-diagonalizable $\be$. In case 2.1, we not only scale $\be$,
but also interchange $X$ and $Y$, so that we can limit ourselves to $-1 \leq r \leq
1$.

If $[X,Y]=Y$, we see that $\be(Y) \in \R Y$, and then
$[\be(X),Y]=0$. Hence, $\be(X) \in \R Y$. Replacing, if necessary, $X$
by $X-rY$ and rescaling $A$, we obtain two different cases: 2.4 and 2.5.

Next, suppose that $\chi \neq 0$ and $\be=0$. Then $\chi$ is a
character. Take $X,Y$ such that $\chi(X)=1$ and $\chi(Y)=0$. Hence
$[X,Y]=aY$ for some $a \in \R$. These are all non-isomorphic: to pass
from one $a \in \R$ to another, we have to multiply $X$ by a scalar,
but then $\chi(X)=1$ is violated.

An exponentiation of all the above matched pairs of Lie algebras again can be
given using semi-direct products of the corresponding connected simply connected Lie groups.

Finally, the general discussion before this theorem shows that in case~4, there
are two special situations.
First, if $\chi \be$ is a multiple of $\chi$, we separate
the cases $\chi \be \neq 0$ and $\chi \be =0$. If $\chi \be \neq 0$, we
rescale $A$, so that $\chi \be = \chi$. We can take generators $X,Y$ for $\lieg_1$
such that $\chi(X)=1$ and $\chi(Y)=0$. Then, there exists $b,c
\in \R$, such that
$$\be(X) = X + bY \; , \quad \be(Y) = c Y
\; .$$ Because $\chi([X,Y]) = 0$, we get
$[X,Y]= dY$ for $d \in \R$. Checking Equation~\eqref{eq1}, we get
$c=d$.
If we replace $X$ by $X+rY$, then $b$ changes to $b + r(d-1)$. Hence, we
get two non-isomorphic families: $d=1$, $b \in \R$ and $d \neq 1$, $b
=0$, as stated in case~4.1.

If $\chi \be =0$, an analogous reasoning gives generators $X,Y$ for
$\lieg_1$ such that $[X,Y]=dY$, $\chi(X)=1$, $\chi(Y)=0$,
$\be(X)=bY$, $\be(Y) = 0$, for $b,d \in \R$. Because $\be \neq 0$, we
get $b \neq 0$. Rescaling $Y$, we can assume that $b=1$. This gives case~4.2.

The existence of the exponentiation of these matched pairs has been
proven in Corollary~\ref{exponentiation_casetwo}, it will be described explicitly below.

The classification in case~4.3 as well as the exponentiation follows
from Proposition~\ref{exponentiation_casethree}. To reduce the number
of connected components, we modify the exponentiation slightly.

If $a>0$, in order to exponentiate, we define
\begin{align*}
G &= \frac{\{X \in M_2(\R) \mid \operatorname{det} (X) = \pm
  1\}}{\{\pm 1\}} \; , \\
G_1 = \{(a,x) \mid a \neq 0, x \in \R\}\;&  , \; (a,x) (b,y) =
(ab,x+ay) \; , \; G_2 = (\R \setminus \{0\}, \cdot) \; .
\end{align*}
Making use of the function $\Sq(a) := \Sgn(a) \sqrt{|a|}$, we define
\begin{align}
i(a,x) &= \begin{pmatrix} \frac{1}{\sqrt{|a|}} & 0 \\ 2
\frac{x}{\sqrt{|a|}} & \Sq(a)
\end{pmatrix} \mod \{\pm 1\}\; , \label{eq3.2} \\
j(s) &=  \begin{pmatrix} \sqrt{|s|} & \frac{1}{2} \bigl(
\sqrt{|s|} - \frac{1}{\Sq(s)} \bigr)
\\ 0 & \frac{1}{\Sq(s)} \end{pmatrix} \mod \{\pm 1\} \; . \notag
\end{align}
It is clear that the tangent Lie algebras of $i(G_1)$ and $j(G_2)$ are
given by $\lspan\{\Htil,\Ytil\}$ and $\R (\Htil + \Xtil)$,
respectively. It is not hard to check that we indeed get a matched
pair of Lie groups.

If $a<0$, the easiest way to
exponentiate this matched pair goes as follows (as we saw in the proof
of Proposition~\ref{exponentiation_casethree}, there is also another way of doing so).

Define
\begin{align*}
G &= \operatorname{PSL}_2(\R) =
\frac{\operatorname{SL}_2(\R)}{\{\pm 1\}} \; , \quad G_2 =
(\mathbb{T} = \{z \in
\mathbb{C} \mid |z| = 1 \}, \cdot) \\
G_1 &= \{(a,x) \mid a > 0, x \in \R \} \;  .
\end{align*}
Then define
\begin{align}
i(a,x) &= \begin{pmatrix} \frac{1}{\sqrt{a}} & 0 \\
\frac{x}{\sqrt{a}} & \sqrt{a} \end{pmatrix} \mod \{ \pm 1 \} \; ,
\label{eq3.3} \\ j(\cos t, \sin t) &= \begin{pmatrix} \cos
\frac{t}{2} & \sin \frac{t}{2} \\ - \sin \frac{t}{2} & \cos
\frac{t}{2} \end{pmatrix} \mod \{ \pm 1 \}  \; . \notag
\end{align}
One can check that the tangent Lie subalgebras of $i(G_1)$
and $j(G_2)$ agree with $\lspan\{\Htil,\Ytil\}$ and $\R(\Xtil -
\Ytil)$, respectively and that we get a matched pair of Lie groups.

The final case $a =0$ has been exponentiated in \cite{VV},
Section~5.4, but we recall it for completeness. We take again $G=\operatorname{PSL}_2(\R)$.
$G_1$ consists of pairs $(a,x)$ with $a > 0$ and $x \in \R$ with
product $(a,x)(b,y) = (ab, ay + \frac{x}{b})$. Putting
\begin{equation} \label{weetikveel}
i(a,x) = \begin{pmatrix} a & x \\
0 & \frac{1}{a} \end{pmatrix} \mod \{ \pm 1 \} \; , \qquad j(s) =
\begin{pmatrix} 1 & 0
\\ s & 1 \end{pmatrix} \mod \{ \pm
1 \},  \;
\end{equation}
we get the required exponentiation to a matched pair of Lie groups.
\end{proof}

For any of the obtained matched pairs of Lie groups, we can now perform the bicrossed
product construction in order to get a l.c.\ quantum
group. Whenever one of the corresponding actions is trivial, we obtain a Kac algebra (see
Corollary \ref{410}). When both actions are non-trivial, we find a lot
of l.c.
quantum groups which are not Kac algebras. To take a closer look at
them, we need explicit forms for the corresponding mutual actions, and we use
the formulas
$$
\chi(X) = X_e[ g \mapsto \frac{d}{ds} \al_g(s) |_{s=0} ] \; , \quad
\be(X) = \frac{d}{ds}( (d\be_s)(X) ) |_{s=0},
$$
where $X_e$ is the partial derivative in $e$ in the direction of an arbitrary generator
$X\in \lieg_1$ and $d\be_s$ is the canonical action of $G_2$ on $\lieg_1$ coming from
$\be_s$.

The only case of dimension $1+1$ with both non-trivial actions has already been presented
in Remark~\ref{remnieuw}. It is easy to check that we do not get a Kac algebra,
that $\sde_M \neq 1$, and that the corresponding l.c.\ quantum group is self-dual. For the
details see \cite{VV}, 5.3.

In dimension $2+1$, we analyze cases~4.1, 4.2 and 4.3.

In case~4.1, following the approach of Proposition~\ref{exponentiation_general}, we define the Lie group $G$ on the space $\R \setminus \{0\}
\times \R^2$ with multiplication
$$(s,x,y) (s',x',y') = (ss',x+sx', y + b u_d(s) x' + s^d y') \; ,$$
where $$u_d(s) = \begin{cases} \frac{s^d - s}{d-1} \; ,
  \quad\text{if}\; d \neq 1 \\ s \log |s| \quad\text{if}\; d = 1
\end{cases} \; , \quad \text{where} \; s^d = \Sgn(s) |s|^d \; ,$$
and $G_1$ on the space $\R \setminus \{0\} \times \R$ with multiplication
$(a,x) (a',x') = (aa',x+a^dx')$ and $i(a,x) = (a, a-1, x+b
u_d(a))$. Further, we put $G_2 = \R \setminus \{0\}$  and $j(s) = (s,0,0)$. Then, the mutual actions are
given by
\begin{align}
\al_{(a,x)}(s) &= a(s-1)+1 \; , \label{exp41} \\ \be_s(a,x) &=
\Bigl(\frac{sa}{a(s-1)+1}, \frac{x+b(u_d(a) + u_d(a(s-1)+1) -
  u_d(as))}{(a(s-1)+1)^d} \Bigr)\; . \notag
\end{align}
One can check that the corresponding matched
pair of Lie algebras is isomorphic to the initial one. Indeed, using the obvious
generators, we have:
$$
[X,Y]=dY, \; \chi(X)=1, \; \chi(Y)=0, \; \be(X) = -X - bd Y, \;\be(Y) = -dY \; .
$$
The needed isomorphism is given by $A\mapsto -A,\ Y\mapsto dY\ (\text{if}\ d\neq 0)$
;
if $d=0,\ A\mapsto -A$ establishes an isomorphism with the special case of the initial
matched pair: $b=d=0$.

Because the modular functions of the groups $G_1$ and $G$ are
given by $\sde_1(a,x)$ $=|a|^{-d}$ and $\sde(s,x,y) = |s|^{-d-1}$,
we compute that the first equality of Proposition~\ref{charKac}
does not hold and
$$\quad \sde_M(a,x,s)
= |a(s-1)+1|^{d+1} \; , \quad \sde_{\hat{M}}(a,x,s) = \Bigl|
\frac{as}{a(s-1)+1} \Bigr|^{d-1} \; .$$
So, both the l.c.\ quantum
group and its dual are not Kac algebras, and are non-unimodular.

In case~4.2, the Lie groups $G$ and $G_1$ are defined on $\R^+ \setminus \{0\}
\times \R^2$ and $\R^+ \setminus \{0\}
\times \R$, respectively, with the same multiplication as in case~4.1,
but with
parameter $b=1$. We
consider $G_2$ to be $\R$ with addition and define $i(a,x) = (a,0,x)$
and $j(s) = (1,s,0)$. Then, the mutual actions are
\begin{equation} \label{exp42}
\al_{(a,x)}(s) = as \; , \qquad \be_s(a,x) = (a, x+ u_d(a) s) \;
.
\end{equation}
One can check that the corresponding matched
pair of Lie algebras coincides with the initial one, that the first equality of
Proposition~\ref{charKac} holds and that $\sde_M =1,\ \sde_{\hat{M}}(a,x,s) =
a^{d-1}$. Hence, $(M,\de)$ is a unimodular Kac algebra, and
$(\hat{M},\hat{\de})$ is unimodular if and only if $d=1$.

Finally, the exponentiations of case~4.3 with
$a>0,<0,=0$, are determined by
Equations~\eqref {eq3.2}, \eqref{eq3.3} and \eqref{weetikveel},
respectively. In the case $a>0$, the
mutual actions are given by
\begin{align}
\al_{(a,x)}(s) &= \frac{(x+1)s + a - x - 1}{xs + a -x} \;
,\label{exp43a4} \\
\be_s(a,x) &= \Bigl( \frac{((x+1)s + a - x - 1) \; (xs + a -x)}{as} \; , \; \frac{x \;
((x+1)s + a - x - 1)}{a} \Bigr) \; . \notag
\end{align}
The corresponding matched pair of Lie algebras
$$[X,Y]=Y, \; \chi(X)=-1, \; \chi(Y)=0, \; \be(X)=-X, \; \be(Y)=2X+Y$$
is isomorphic to the initial one: $X\mapsto -X-\frac{Y}{2},\ Y\mapsto -\frac{Y}{4},\ A\mapsto 2A$.

In the case $a<0$, the mutual actions are
\begin{align}
& \al_{(a,x)}(\cos t, \sin t) = \label{exp43a-4} \\
&\qquad \frac{\Bigl( (a^2 + x^2 - 1) + (a^2-x^2+1) \cos t + 2ax
\sin t, 2x - 2x \cos t + 2a \sin t \Bigr)}{(x^2+a^2+1) + (-x^2 +
a^2 - 1) \cos t
+ 2ax \sin t} , \notag \\
& \be_{(\cos t, \sin t)}(a,x) =
\Bigl(\frac{1}{2a}\bigl((x^2+a^2+1) + (-x^2 + a^2 - 1) \cos t +
2ax \sin t\bigr), \notag \\ & \hspace{5cm} \frac{1}{2a}\bigl((x^2
- a^2 + 1) \sin t + 2ax \cos t\bigr) \Bigr)\; . \notag
\end{align}
Observe that these actions are everywhere defined
and continuous. The corresponding matched pair of Lie algebras
$$[X,Y]=Y, \; \chi(X)=-1, \; \chi(Y)=0, \; \be(X)=-Y, \; \be(Y)=X$$
is isomorphic to the initial one: $X\mapsto -X,\ Y\mapsto \frac{Y}{2},\ A\mapsto -2A$.

In the case $a=0$, the mutual actions are
\begin{equation} \label{exp43a0}
\al_{(a,x)}(s) = \frac{s}{a(a+xs)} \; , \quad \be_s(a,x) = (|a+sx|,
\operatorname{Sgn}(a+sx) x) \; .
\end{equation}
The corresponding matched pair of Lie algebras
$$[X,Y]=2Y , \;  \chi(X)=-2 , \; \chi(Y)=0 , \; \be(X)=0 ,\; \be(Y)=X$$
is isomorphic to the initial one: $X\mapsto -\frac{X}{2},\ Y\mapsto -\frac{Y}{2}$.

In all three cases $a>0,<0,=0$, one verifies that the first equality of
Proposition~\ref{charKac} does not hold, $\sde_M=1$,
while $\sde_{\hat{M}} \neq 1 $. In particular, we do not get Kac
algebras.

\section{Cocycle matched pairs of Lie groups and Lie algebras in low dimensions}

So far, we have explained how to construct l.c.\ quantum groups which are bicrossed products
of low-dimensional Lie groups without 2-cocycles. The usage of 2-cocycles gives much more concrete
examples and,
what is more important, gives a more complete picture of low-dimensional l.c.\ quantum
groups. Again, we first explain the infinitesimal picture, i.e. how
2-cocycles for matched pairs of Lie algebras look like, how they are related to the
problem of extensions and then show how to exponentiate them.

The first thing that we need here is the notion of a Lie bialgebra, due to V.G. Drinfeld
\cite{Drin}.
A Lie bialgebra is a Lie algebra $\lieg$ equipped with a Lie bracket $[\cdot,\cdot ]$ and a Lie
cobracket $\delta$, i.e., a linear map $\delta:\lieg\to\lieg\otimes \lieg$ satisfying
the co-anticommutativity and the co-Jacobi identity, that is:
$$
(\io - \tau)\delta = 0 \; , \quad  (\io + \zeta + \zeta^2)(\io\otimes\delta)\delta=0,
$$
where
$$
\tau(u\otimes v)=v\otimes u, \quad \zeta(u\otimes v\otimes w)=v\otimes w\otimes u \quad
\ (\text{for all}\quad u,v,w\in \lieg)
$$
are the flip maps, and these Lie bracket and cobracket are compatible in the following
sense:
$$
\delta[u,v]=[u,v_{[1]}]\otimes v_{[2]}+v_{[1]}\otimes [u,v_{[2]}]+[u_{[1]},v]\otimes u_{[2]}
+u_{[1]}\otimes [u_{[2]},v].
$$
Any Lie algebra (respectively, Lie coalgebra, i.e., vector space dual to a Lie algebra)
is a Lie bialgebra with zero Lie cobracket (respectively, zero Lie bracket). The
definition of a morphism of Lie bialgebras is obvious.

Given a pair of Lie algebras $(\lieg_1,\lieg_2)$, let us ask if there exists a Lie
bialgebra $\lieg$ such that
$$
\lieg_2^* \longrightarrow \lieg
\longrightarrow \lieg_1
$$
is a short exact sequence in the category of Lie bialgebras. This means precisely
that $\lieg$ has a sub-bialgebra with trivial bracket, which is an ideal and such that
the quotient is a Lie bialgebra with trivial cobracket.

The theory of extensions in this framework has been developed in \cite{Mas2} and is quite
similar to the theory of extensions of l.c.\ groups that we have recalled above. Namely, for the
the existence of an extension $\lieg$ it is necessary and sufficient that $(\lieg_1,
\lieg_2)$
form a matched pair, and all extensions are
bicrossed products with cocycles. We consider this theory as an
infinitesimal version of the theory of extensions of Lie groups.

As we remember, for any matched pair of Lie algebras $(\lieg_1, \lieg_2)$,
there are mutual actions $\triangleright:
\lieg_2\otimes \lieg_1\to \lieg_1$ and $\triangleleft: \lieg_2\otimes \lieg_1\to \lieg_2$, compatible in a way
explained in Section 3 and such that
for all $a,b\in \lieg_1,\ x,y\in \lieg_2$ we have
$$
[a\oplus x,b\oplus y]=([a,b]+x\triangleright b - y\triangleright a)\oplus
([x,y]+x\triangleleft b-y\triangleleft a).
$$

For the general definition of a pair of 2-cocycles on such a matched
pair, we refer to \cite{Majbook}, \cite{Mas2}. For our needs, it suffices
to understand that these 2-cocycles are linear maps
$$
\cU : \lieg_1 \wedge \lieg_1 \recht
\lieg_2^* \; , \qquad \cV : \lieg_2 \wedge \lieg_2 \recht \lieg_1^* \;
$$
verifying certain 2-cocycle equations and compatibility equations
that are infinitesimal forms of Equations (\ref{cocU}). For the case of
dimension $n+1$, we give these equations explicitly below.
Let us formulate the link
between 2-cocycles on matched pairs of Lie algebras and those of Lie groups as
a proposition whose proof is straightforward.

\begin{proposition} \label{cocycles}
Let $(G_1,G_2)$ be a matched pair of Lie groups equipped with
cocycles $\cU$ and $\cV$, which are differentiable around the unit
elements, and let $(\lieg_1, \lieg_2)$ be the corresponding matched
pair of Lie algebras. Defining
\begin{align*}
\langle \cU(X,Y), A \rangle &= -i (X_e \ot Y_e \ot A_e - Y_e \ot
X_e \ot A_e)(\cU) \quad\text{and}\\ \langle \cV(A,B), X \rangle &=
-i (A_e \ot B_e \ot X_e - B_e \ot A_e \ot X_e)(\cV) \; ,
\end{align*}
for $X,Y \in \lieg_1$ and $A,B \in \lieg_2$, we get a
pair of cocycles on $(\lieg_1, \lieg_2)$.
\end{proposition}
Here $\langle \cdot , \cdot \rangle$ denotes the duality between $\lieg_i$
and $\lieg_i^*$ and $X_e,Y_e,A_e,B_e$ denote the partial derivatives at $e$ in
the direction of the corresponding generator. The factor $-i$ appears because
for Lie groups $\cU$ and $\cV$ take values in $\T$, and for real Lie algebras
we consider 2-cocycles as real linear maps.

For the dimension $n+1$, $\cV$ is necessarily trivial (if also
$n=1$, then also $\cU$ is trivial, so there are no non-trivial cocycles in the
dimension $1+1$). Returning to arbitrary $n$, we choose a generator $A$
for $\lieg_2$ and define maps $\be$ and $\chi$ as above. Then $\cU$
can be regarded as an antisymmetric, bilinear form on $\lieg_1$, and
the 2-cocycle equations of \cite{Majbook},\cite{Mas2} reduce to the equation
$$
\cU([X,Y],Z) + \chi(X) \cU(Y,Z) + \; \text{cyclic permutation} \; =
0 \quad\text{for all}\quad X,Y,Z \in \lieg_1 \; .
$$
It is clear that these 2-cocycles $\cU$ form a real vector space.

In Section 2, we defined the notion of the group of extensions for a
matched pair of Lie groups $(G_1,G_2)$ using the notion of cohomologous
2-cocycles. The same can be done for a matched pair of Lie algebras \cite{Mas2}.
In particular, for the dimension $n+1$, a 2-cocycle $\cU$ is called
\emph{cohomologous to trivial}, if there exists a
linear form $\rho$ in $\lieg_1^*$ such that
$$\cU(X,Y) = \rho([X,Y]) + \chi(X) \rho(Y) - \chi(Y) \rho(X) \; .$$
Two cocycles $\cU_1$ and $\cU_2$ are called \emph{cohomologous} if
$\cU_1-\cU_2$ is cohomologous to trivial.

The quotient space of 2-cocycles modulo 2-cocycles cohomologous to trivial,
with addition as the group operation, is called the \emph{group of extensions} of
the matched pair $(\lieg_1, \lieg_2)$.

Now let us describe all 2-cocycles on the matched pairs of real Lie
algebras of dimension $2+1$.

\begin{proposition} \label{prop4.1}
Referring to the classification of matched pairs of Lie algebras of
dimension $2+1$ given in Theorem~\ref{prop3.2}, the following
holds: the group of extensions is $\R$ in the cases 1.1, 2.1, 2.2,
2.3, 3 ($a=-1$), 4.1 ($d=-1$), 4.2 ($d=-1$) and 4.3. The cocycles are
defined by $\cU(X,Y)= \lambda$, for $\lambda \in \R$. In the other
cases, the group of extensions is trivial.
\end{proposition}
\begin{proof}
Since $\dim\ \lieg_1=2$, any
antisymmetric bilinear form $\cU$ on $\lieg_1$ is a cocycle. Suppose
that $X,Y$ are generators of $\lieg_1$ with $[X,Y]=dY$. A cocycle
$\cU$ is entirely determined by $\cU(X,Y)=\lambda$ for $\lambda \in
\R$. If $\chi=0$ and $d \neq 0$, we take $\rho(Y) = \lambda/d$ and get
that $\cU$ is cohomologous to trivial. If $\chi=0$ and $d=0$, it is clear that
$\cU$ is not cohomologous to trivial if $\lambda \neq 0$.

If $\chi \neq 0$, we may suppose that $\chi(X)=1$ and $\chi(Y)=0$. If
$d \neq -1$, we take $\rho(Y) = \lambda / (1+d)$ and get that $\cU$ is
cohomologous to trivial. If $d = -1$, it is clear that
$\cU$ is not cohomologous to trivial if $\lambda \neq 0$.
\end{proof}

Next, we want to exponentiate a cocycle on a matched pair of real Lie algebras,
i.e., to construct a measurable map $\cU: G_1 \times G_1
\times G_2 \recht U(1)$, with values in the unit circle of $\C$, satisfying
\begin{align*}
\cU(g,h,\al_k(s)) \; \cU(gh,k,s) &= \cU(h,k,s) \; \cU(g,hk,s), \\
\cU(g,h,s) \; \cU(\be_{\al_h(s)}(g),\be_s(h),t) &= \cU(g,h,ts)
\end{align*}
almost everywhere. Let us define a function $A(\cdot)$ by
$$\cU(g,h,s) = \exp(i A(g,h,s)) \; .$$
So $A(\cdot)$ should satisfy
\begin{align*}
A(g,h,\al_k(s)) + A(gh,k,s) &= A(h,k,s) + A(g,hk,s) \quad\text{mod}\; 2 \pi, \\
A(g,h,s) + A(\be_{\al_h(s)}(g),\be_s(h),t) &= A(g,h,ts) \quad\text{mod}\; 2 \pi
\end{align*}
almost everywhere.

\begin{proposition} \label{prop4.2}
If the group of extensions of a matched pair of Lie algebras of
dimension $2+1$ is non-trivial, there exists an exponentiation of this
matched pair with cocycles. These cocycles are labeled by $\R$ in all
the cases, except case 4.3 ($a=0,1$), where they are labeled by $\Z$.
\end{proposition}
\begin{proof}
Following \cite{VV}, Section~5.5, we look for the above function $A$ in the form
$$A(g,h,s) = P \int_0^s f(\phi_r(g,h)) \; dr \; ,$$
where $\phi_r(g,h):=(\be_{\al_h(r)}(g),\be_r(h))$ and where the function $f$ on $G_1 \times G_1$
is such that for almost all $g,h \in G_1$ the function $r \mapsto f(\phi_r(g,h))$
has a principal value integral over any interval in $\R$ ($dr$ is the Haar measure
on the 1-dimensional Lie group $(\R,+)$ or on $\R \backslash\{0\}$, in
which case we integrate from $1$ to $s$).
A necessary condition to be satisfied by $f$ is
$$\frac{d}{dt} \al_k(t) \big|_{t=0} \; f(g,h)+f(gh,k)=f(h,k) + f(g,hk)
\; .$$
Finally, having found such an $f$, we have to check if it really gives rise to a 2-cocycle.

In those cases where the actions $\al$ and $\be$ are everywhere
defined and smooth, one can check
that any smooth solution of this equation gives indeed
rise to a 2-cocycle (for the details see \cite{VV}, Section~5.5). In this way,
it is easy to find 2-cocycles in the
cases 1.1, 2.1, 2.2, 2.3, 3 ($a=-1$), and 4.2 ($d=-1$), namely:
in the cases 1.1, 2.1, 2.2 and 2.3, the action $\al$ is trivial, and
$G_1 = \R^2$ with addition. So, we can take $f(x_1,x_2;y_1,y_2)=
\lambda (x_1
y_2 - x_2 y_1)$, for any $\lambda \in \R$. In the cases 3 ($a=-1$) and
4.2 ($d=-1$), we observe that $G_1 = \{(a,x) \mid a > 0, x \in \R\}$
with $(a,x)(b,y) = (ab, x+ y/a)$. Because $\chi(X)=1$, the character
$g \mapsto \frac{d}{dt} \al_g(t) \big|_{t=0}$ is given by $(a,x)
\mapsto a$, and we can take $f(a,x;b,y)= \lambda abx \log b$, for any
$\lambda \in \R$.

The case 4.3 ($a=0$) has been studied in \cite{VV}, Section~5.5:
$f(a,x;b,y) = \lambda \frac{x \log b}{a b^2}$. Checking if we really get
2-cocycles, observe that
$$
f(\phi_r(a,x;b,y)) = \frac{\lambda x}{(b+ry) (ab + r(ay +
\frac{x}{b}))} \log |c + dr| \; ,
$$
then
$$
P \int_{-\infty}^\infty f(\phi_r(a,x;b,y)) \; dr =
\frac{\lambda}{2} \pi^2 \operatorname{Sgn} \bigl( \frac{y}{x}
(ay+\frac{x}{b}) \bigr) \; .
$$
From this, it follows that we do
get 2-cocycles if and only if $\lambda = \frac{4n}{\pi}$, with $n
\in \Z$ .

The same phenomenon happens in the cases 4.1 ($d=-1$) and 4.3 ($a>0$):
although the above principal
value integral is well defined, we do
not always get a 2-cocycle $\cU$, as explained in \cite{VV}, after Proposition~5.6.
In case 4.1 ($d=-1$), we use the matched pair
explicitly described in Equation~\eqref{exp41} with $d=-1$ and $b=0$. We take again
$f(a,x;b,y)= \lambda abx \log |b|$, for any
$\lambda \in \R$, and we can explicitly perform the integration, to
obtain 2-cocycles: $$A(a,x;b,y;s) = \lambda a x \bigl( - (b(s-1)+1) \log
|b(s-1)+1| + bs \log |bs| - b \log |b| \bigr) \; .$$

On the contrary, the situation of case 4.3 ($a>4$) is more
delicate. We use the matched pair explicitly described in
Equation~\eqref{exp43a4}. We can take $f(a,x;b,y)=\lambda \frac{y}{b} \log
|a|$ with $\lambda \in \R$, and our candidate for $A(\cdot)$ becomes:
\begin{align*}
& A(a,x;b,y;s)= \lambda \; P\!\!\int_1^s \frac{y}{yr + b - y}
\\ &\qquad \log \Bigl| \frac{ \bigl( (x+ay+1)r + ab - x - ay -1 \bigr) \;
\bigl( (x+ay)r + ab -x - ay \bigr)}{a \bigl( (y+1)r + b - y
-1\bigr) \; \bigl(yr + b - y\bigr)} \Bigr| \; dr \; .
\end{align*}
Because
$$P \!\!\int_{-\infty}^{+\infty} \frac{c}{cr + d} \log | ar + b| \; dr = \frac{\pi^2}{2}
\Sgn(\frac{b}{a} - \frac{d}{c}) \; ,$$ the same reasoning as in
Section~5.5 of \cite{VV} implies that we do get a 2-cocycle if $\lambda
= \frac{4n}{\pi}$ for $n \in \Z$.

Finally, in case 4.3 ($a>0$), with the explicit exponentiation given in
Equation~\eqref{exp43a-4}, the mutual actions are
defined everywhere and are smooth, but $G_1 = \T$. Taking
$f(a,x;b,y) = \lambda \frac{y}{b} \log a$ with $\lambda \in \R$, it is natural to use
$$A(a,x;b,y;\cos t, \sin t) = \lambda \int_0^t f \bigl(
\phi_{(\cos s,\sin s)}(a,x,b,y) \bigr) \; ds \; .$$
To have a 2-cocycle, we need that
$$\lambda \int_0^{2\pi} f \bigl( \phi_{(\cos s,\sin s)}(a,x,b,y) \bigr) \; ds = 0 \mod 2\pi \; .$$
Denote the left-hand side of this expression by $I_\lambda(a,x,b,y)$. Then, one can compute
that $$I_\lambda(a,x,b,y) = \lambda \bigl( H(a,x) + H(b,y) - H(ab,x+ay) \bigr) \; ,$$
where $H(a,x) := - 4 \pi \arctan \bigl( \frac{x}{1+a} \bigr)$. Hence,
there is no $\lambda$ which gives us a 2-cocycle.

We can, however, find 2-cocycles, using the other exponentiation of
the same matched pair of Lie algebras, as explained in the proof of Proposition~\ref{exponentiation_casethree}.
We obtain a matched pair $(G_1,\R)$, in which both
actions are everywhere defined and smooth. Hence, we obtain cocycles
labeled by $\R$, following the procedure described in the beginning of
the proof.
\end{proof}

\section{Infinitesimal objects for low-dimensional l.c.\ quantum groups}

\subsubsection*{The case of cocycle bicrossed product l.c.\ quantum groups}

Given a cocycle matched pair of Lie groups $(G_1,G_2)$ whose
2-cocycles $\cU$ and $\cV$ are differentiable around the unit
elements, we can construct the corresponding l.c.\ quantum group
using the cocycle bicrossed product construction.

But, in this situation,  we can also construct  two other intimately
related algebraic structures which can be viewed as infinitesimal
objects of this l.c.\ quantum group: a Lie bialgebra and a Hopf
$^*$-algebra, as in \cite{VV}, Section~5.2. The precise
mathematical link between these three  structures  is
not completely clear at the moment (it is tempting to consider it as a kind of a
Lie theory for our cocycle bicrossed product l.c. quantum groups).
We will discuss it mainly on the level of examples.

Let us recall the construction of infinitesimal Lie bialgebras and Hopf
$^*$-algebras in the special case of dimension $n+1$.

Let $(G_1,G_2)$ be a cocycle matched pair of Lie groups with $G_2 = \R$
(the case $\R \setminus \{0\}$ is completely
analogous, replacing differentials in $0$ by differentials in
$1$) and let $\cU$ be a 2-cocycle differentiable around the unit elements. Denote by
$\al_g(s)$ and $\be_s(g)$ the corresponding mutual actions. Then the
cocycle matched pair of Lie algebras is determined (see Section 4 and
Proposition \ref{cocycles}) by
\begin{align*}
\chi(X) &= X_e[ g \mapsto \frac{d}{ds} \al_g(s) |_{s=0} ] \; , \quad X \in \lieg_1 \; , \\
\be(X) &= \frac{d}{ds}( (d\be_s)(X) ) |_{s=0}\; , \quad X \in \lieg_1 \; , \\
\cU(X,Y) &= -i \frac{d}{ds} ( (X_e \ot Y_e - Y_e \ot
X_e)(\cU(\cdot,\cdot,s))) |_{s=0} ] \; \Atil \; , \quad X,Y \in \lieg_1 \;
.
\end{align*}
The infinitesimal Lie bialgebra is precisely the corresponding cocycle
bicrossed product Lie bialgebra and has generators $\Atil\in \lieg_2$ and $X \in
\lieg_1$, subject to the relations
\begin{align*}
[\Atil,X] &= \chi(X) \Atil \; , \\
[X,Y] &= [X,Y]_1 + \cU(X,Y) \Atil \; , \\
\sde(\Atil) &= 0 \; , \\
\sde(X) &= \be(X) \wedge \Atil \; .
\end{align*}

The dual infinitesimal Lie bialgebra has generators $A$ and
$\Xtil_i$, subject to the relations
\begin{align*}
[\Xtil_i,A] &= \sum_j \be(X_i)_j \Xtil_j \; , \\
[\Xtil_i,\Xtil_j] &= 0 \; , \\
\langle \sde(\Xtil_i), X \ot Y \rangle &= \langle \Xtil_i,[X,Y]_1 \rangle \; , \\
\sde(A) &= A \wedge \bigl( \sum_i \chi(X_i) \Xtil_i \bigr) +
\sum_{i < j} \cU(X_i,X_j) \Xtil_i \wedge \Xtil_j \; .
\end{align*}

Following \cite{VV}, Section~5.2 and in order to construct the infinitesimal Hopf
$^*$-algebra, which is an algebraic cocycle bicrossed product in the sense of
\cite{Majbook}, we denote by $\Atil$ the function $\Atil(x) = x$ on
$G_2$. Then, our Hopf $^*$-algebra has generators $\Atil$ and $\{X
\mid X \in \lieg_1\}$, where $\Atil^* = \Atil$ and $X^* = -X$,
with relations
\begin{align*}
[\Atil,X] & = X_e [g \mapsto \al_g(s) ] \; , \quad X \in \lieg_1
\; , \\
[X,Y] & = [X,Y]_1 + (X_e \ot Y_e - Y_e \ot
X_e)(\cU(\cdot,\cdot,s)) \;, \quad X,Y \in \lieg_1 \; , \\
\de(\Atil) & = \Atil \ot 1 + 1 \ot \Atil \; , \quad (\Atil \ot
\Atil \; \text{when} \; G_2 = \R \setminus \{0\}) \; , \\
\de(X) &= 1 \ot X + \sum_i X_i \ot (d\be_s)(X)_i \; , \quad X \in
\lieg_1 \; .
\end{align*}
This needs some explanation: several of the used expressions are
functions of $s \in G_2$, and it is understood that these
functions belong to some algebra of functions in $\Atil$ (in the simplest case -
to the algebra of polynomials in $\Atil$). Next,
$(X_i)$ is a basis for $\lieg_1$ and $(d\be_s)(X)_i$ is the $i$-th
component of $(d\be_s)(X)$ in this basis, and is, therefore,
again a function of $s \in G_2$. Finally, $[X,Y]_1$ denotes the
Lie bracket in $\lieg_1$. When $G_2 = (\R,+)$, co-unit and antipode are given by
$\vep(\Atil) = \vep(X)=0$ ($X \in \lieg_1$), $S(\Atil)=-\Atil$ and
$$S(X) = - \sum_i X_i \; (d \be_{-s})(X)_i \; , \quad X \in \lieg_1 \;
.$$
When $G_2=(\R \setminus \{0\},\cdot)$, we rather have $\vep(\Atil)=1, \vep(X)=0
$ ($X \in \lieg_1$), $S(\Atil)= \Atil^{-1}$ and $$S(X) = - \sum_i X_i
\; (d \be_{1/s})(X)_i \; , \quad  X \in \lieg_1 \; .$$
To describe the dual infinitesimal Hopf $^*$-algebra, suppose that we have
coordinate functions $(\Xtil_i)$ on $G_1$, dual to the basis
$(X_i)$ of $\lieg_1$. Then we can write down, with the same kind of
conventions, the dual Hopf $^*$-algebra, with generators
$(\Xtil_i)$ and $A$, such that $\Xtil_i^* = \Xtil_i$ and $A^* =
-A$, and with relations:
\begin{align*}
[\Xtil_i,A] &= \frac{d}{ds} (\Xtil_i(\be_s(g))) |_{s=0} \; , \\
[\Xtil_i,\Xtil_j] &= 0 \; , \\
\de(\Xtil_i)(g,h) &= \Xtil_i(gh) \; , \\
\de(A) &= 1 \ot A + A \ot \frac{d}{ds} (\al_g(s)) |_{s=0} +
\frac{d}{ds} \cU(g,h,s) |_{s=0} \; .
\end{align*}
Again, we write functions of $g,h \in G_1$ and they are tacitly
assumed to belong to some algebra of functions in $\Xtil_i$. Co-unit
and antipode are given by $\vep(\Xtil_i)=\Xtil_i(e)$, $\vep(A)=0$,
$S(\Xtil_i)(g)=\Xtil_i(g^{-1})$ and
$$S(A) = - A \; \frac{d}{ds} (\al_{g^{-1}}(s)) |_{s=0} - \frac{d}{ds}
\cU(g,g^{-1},s) |_{s=0} \; .$$

\begin{remark}~\label{Hopf-bialg} Observe that in order to pass
from the infinitesimal
Hopf $^*$-algebra to the infinitesimal Lie bialgebra, one
replaces functions on $G_2$ by their first-order approximations,
which are multiples of $\Atil$, one takes anti-auto-adjoint
generators $i \Atil$ and $X \in \lieg_1$ and takes as $\sde$ the
linear part of $i (\de - \de^{\text{op}})$, where
$\de^{\text{op}}$ denotes the opposite comultiplication.
\end{remark}

Let us turn now to concrete examples. For the matched pair of dimension
1+1 (see Remark~\ref{remnieuw}) we do not have cocycles. The infinitesimal Hopf
$^*$-algebra is determined by
\begin{align*}
& \Atil = \Atil^* \; , \quad X^* = - X \; , \quad [\Atil,X] =
\Atil - 1  \; , \\ &\de(\Atil) = \Atil \ot \Atil \; , \quad \de(X)
= X \ot \Atil^{-1} + 1 \ot X \; . \end{align*} Co-unit and
antipode are $\vep(\Atil)=1, \vep(X)=0, S(\Atil)=\Atil^{-1}, S(X)
= -X \Atil$. The infinitesimal Lie bialgebra is given by
$$[\Atil,X] = \Atil \; , \quad \sde(\Atil)=0 \; , \quad \sde(X) =
\Atil \wedge X \; .$$ We will now make the link between the Hopf
$^*$-algebra and the Lie bialgebra slightly more explicit. For a
smooth function $f$ in $\Atil$ we get $[f(\Atil),X] =
f'(\Atil)(\Atil - 1)$. If we take $A_1:=i \log \Atil$ as a new
generator, we observe that
\begin{align*}
& [A_1,X] = A_1 + \frac{i}{2} A_1^2 + \cdots \; , \quad \de(A_1) =
A_1 \ot 1 - 1 \ot A_1 \; , \\ & \de(X) = X \ot (1 + i A_1 -
\frac{1}{2} A_1^2 + \cdots) + 1 \ot X \; .\end{align*} Linearizing
now the above bracket and $i (\de - \de^{\text{op}})$, we
precisely get our Lie bialgebra.

Since this example is self-dual, the dual infinitesimal Hopf $^*$-algebra and Lie
bialgebra are the same as the original ones.

\begin{remark}~\label{ax+b} It is easy to show that the only
non-trivial Lie bialgebras of dimension 2 non-isomorphic to the above mentioned
Lie bialgebra, are defined by the relations
$$
[A,X] = X \; , \quad \sde(A) = 0 \; , \quad \sde(X) = q X \wedge
A \; .
$$
On the other hand, the Hopf $^*$-algebras corresponding to the l.c.\ quantum
"ax+b"-group considered by S.L.~Woronowicz and S. Zakrzevski \cite{WZ} and by
A.~Van Daele \cite{VanDaele2} are defined by the relations
$$
a^*=a, \; x^* = x \; , \quad ax = \exp(iq) xa \; , \quad \de(a) = a
\ot a \; , \quad \de(x) = x \ot a + 1 \ot x \; .
$$
Although this quantum group cannot be obtained by the bicrossed product
construction (see \cite{VV}, 5.4.d), we can apply to it the same formal procedure
to associate with it a Lie bialgebra. Namely, observing that
$[\log a, x] = i q x$, we put $A = \frac{-i}{q} \log a$ and
$X=ix$. Linearizing as above, we find precisely the previous Lie
bialgebra which then can be considered as an infinitesimal Lie bialgebra
of the l.c.\ quantum "ax+b"-group (we mentioned already that the exact relation
between these structures is not completely clear).
\end{remark}

Let us now list the infinitesimal Hopf $^*$-algebras and Lie
bialgebras of all the cocycle matched pairs of Lie algebras of
dimension 2+1, with two non-trivial actions (referring to Theorem~\ref{prop3.2},
case 4).

Case 4.1 has been exponentiated concretely in
Equation~\eqref{exp41}. Taking the obvious generators $X$ and $Y$
for $\lieg_1$ (differentiating to $a$ and $x$, respectively),
satisfying $[X,Y]=dY$, and $A$ for $\lieg_2$, we get the
infinitesimal form $\chi(X) = 1$, $\chi(Y) = 0$, $\be(X) = -X-bdY$
and $\be(Y) = -dY$. When $d \neq -1$, we do not have non-trivial
cocycles and we get the Hopf $^*$-algebra, with $\Atil^* = \Atil$,
$X^*=-X$, $Y^* = -Y$ and
\begin{align*}
[\Atil,X] &= \Atil - 1 \; , \quad [\Atil,Y]=0 \; , \quad [X,Y] =
dY \; , \\
\de(\Atil) &= \Atil \ot \Atil \; , \\
\de(X) &= X \ot \Atil^{-1} + 1 \ot X + b Y \ot
\Atil^{-d}(1-u'_d(\Atil)) \; , \\
\de(Y) &= Y \ot \Atil^{-d} + 1 \ot Y \; .
\end{align*}
Co-unit and antipode are $\vep(\Atil)=1, \vep(X)=\vep(Y)=0,
S(\Atil)=\Atil^{-1}, S(X) = -X \Atil -bd Y u_d(\Atil), S(Y)=-Y \Atil^d$.
The corresponding Lie bialgebra is
\begin{align*}
[\Atil,X] &= \Atil \; , \quad [\Atil,Y] = 0 \; , \quad [X,Y] = dY
\; , \\
\sde(\Atil) &= 0 \; , \quad \sde(X) = \Atil \wedge (X + bdY) \; ,
\quad \sde(Y) = d \Atil \wedge Y \; .
\end{align*}
When $d=-1$, we have found non-trivial cocycles, which are
infinitesimally given by $\cU(X,Y) = - \lambda \Atil$. For the
Hopf $^*$-algebra, there changes $[X,Y] = -Y -i \lambda \log
\Atil$, and on the Lie bialgebra level, this becomes $[X,Y]=-Y -
\lambda \Atil$. Co-unit and antipode remain unchanged.

The dual  Hopf $^*$-algebra for $d \neq -1$ has anti-self-adjoint generators
$\Xtil,\Ytil$, and a self-adjoint generator $A$, subject to the relations
\begin{align*}
[\Xtil,A] &= \Xtil (1-\Xtil) \; , \quad [\Ytil,A] = b \Xtil
(1-u'_d(\Xtil)) - d \Xtil \Ytil \; , \quad [\Xtil,\Ytil]=0 \; , \\
\de(\Xtil) & = \Xtil \ot \Xtil \;  , \\
\de(\Ytil) &= \Xtil^d \ot \Ytil + \Ytil \ot 1 \; , \\
\de(A) &= A \ot \Xtil + 1 \ot A \; .
\end{align*}
Co-unit and antipode are $\vep(\Xtil)=1, \vep(\Ytil)=\vep(A)=0,
S(\Xtil)=\Xtil^{-1}, S(\Ytil)=-\Xtil^{-d}\Ytil, S(A)=-A \Xtil^{-1}$.
The corresponding infinitesimal Lie bialgebra is determined by
\begin{align*}
[\Xtil,A] &= - \Xtil \; , \quad [\Ytil, A] = -bd \Xtil -d  \Ytil
\; , \quad [\Xtil,\Ytil] = 0 \; , \\
\sde(\Xtil) &= 0 \; , \quad \sde(\Ytil) = d \Xtil \wedge \Ytil \;
, \quad \sde(A) = A \wedge \Xtil \; .
\end{align*}
In the case $d=-1$ and with the same cocycle as above, we should
change $\de(A) = A \ot \Xtil + 1 \ot A + i \lambda \Xtil \Ytil \ot
\Xtil \log \Xtil$ in the definition of the Hopf $^*$-algebra and $\sde(A) = A \wedge
\Xtil - \lambda \Xtil \wedge \Ytil$ in the definition of the Lie
bialgebra. For the antipode, we change $S(A) = -A \Xtil + i \lambda
\Ytil \log \Xtil$.

Case 4.2 has been exponentiated in Equation~\eqref{exp42}.
Taking again the obvious generators, we
have the infinitesimal form $[X,Y]=dY$, $\be(X) = Y$, $\be(Y)=0$,
$\chi(X)=1$ and $\chi(Y)=0$. When $d \neq -1$, there are no
non-trivial cocycles. We find the Hopf $^*$-algebra with
generators $\Atil=\Atil^*$ and $X,Y$ anti-self-adjoint, satisfying
\begin{align*}
[\Atil,X] &= \Atil \; , \quad [\Atil,Y]=0 \; , \quad [X,Y] =
dY \; , \\
\de(\Atil) &= \Atil \ot \Atil \; , \\
\de(X) &= X \ot 1 + 1 \ot X + Y \ot
\Atil \; , \\
\de(Y) &= Y \ot 1 + 1 \ot Y \; .
\end{align*}
Co-unit and antipode are $\vep(\Atil)=1, \vep(X)=\vep(Y)=0,
S(\Atil)=\Atil^{-1}, S(X) = -X  -Y \Atil^{-1}, S(Y)=-Y$.
The corresponding Lie bialgebra is
\begin{align*}
[\Atil,X] &= \Atil \; , \quad [\Atil,Y] = 0 \; , \quad [X,Y] = dY
\; , \\
\sde(\Atil) &= 0 \; , \quad \sde(X) = Y \wedge \Atil \; ,
\quad \sde(Y) = 0 \; .
\end{align*}
When $d=-1$, we have found non-trivial cocycles, which are
infinitesimally given by $\cU(X,Y) = - \lambda \Atil$. For the
Hopf $^*$-algebra, there changes $[X,Y] = -Y -i \lambda
\Atil$, and on the Lie bialgebra level, this becomes $[X,Y]=-Y -
\lambda \Atil$. Co-unit and antipode remain unchanged.

The dual Hopf $^*$-algebra for $d \neq -1$ has anti-self-adjoint
generators $\Xtil,\Ytil$ and a self-adjoint generator $A$, subject to the
relations
\begin{align*}
[\Xtil,A] &= 0 \; , \quad [\Ytil,A] = u_d(\Xtil) \; , \quad [\Xtil,\Ytil]=0 \; , \\
\de(\Xtil) & = \Xtil \ot \Xtil \;  , \\
\de(\Ytil) &= \Xtil^d \ot \Ytil + \Ytil \ot 1 \; , \\
\de(A) &= A \ot \Xtil + 1 \ot A \; .
\end{align*}
Co-unit and antipode are $\vep(\Xtil)=1, \vep(\Ytil)=\vep(A)=0,
S(\Xtil)=\Xtil^{-1}, S(\Ytil)=-\Xtil^{-d}\Ytil, S(A)=-A \Xtil^{-1}$.
The corresponding infinitesimal Lie bialgebra is determined by
\begin{align*}
[\Xtil,A] &= 0 \; , \quad [\Ytil, A] = \Xtil
\; , \quad [\Xtil,\Ytil] = 0 \; , \\
\sde(\Xtil) &= 0 \; , \quad \sde(\Ytil) = d \Xtil \wedge \Ytil \;
, \quad \sde(A) = A \wedge \Xtil \; .
\end{align*}
In the case $d=-1$ and with the same cocycle as above, we should
change $\de(A) = A \ot \Xtil + 1 \ot A + i \lambda \Xtil \Ytil \ot
\Xtil \log \Xtil$ in the definition of the Hopf $^*$-algebra and $\sde(A) = A \wedge
\Xtil - \lambda \Xtil \wedge \Ytil$ in the definition of the Lie
bialgebra. For the antipode, we change $S(A)=-A \Xtil^{-1} + i \lambda
\Ytil \log \Xtil$.

Case 4.3 $(a<0)$. For the exponentiation
given in Equation~\eqref{exp43a4}, we have natural generators $X,Y$ and $A$, giving the
infinitesimal form $[X,Y]=Y$, $\chi(X)=-1$, $\chi(Y)=0$, $\be(X)=-X$
and $\be(Y)=2X+Y$. There are cocycles, infinitesimally given by
$\cU(X,Y)= - \lambda \Atil$. The Hopf $^*$-algebra has generators
$\Atil=\Atil^*$ and $X,Y$, anti-self-adjoint, satisfying
\begin{align*}
[\Atil,X] &= 1-\Atil \; , \quad [\Atil,Y]= -(1-\Atil)^2 \; , \quad [X,Y] =
Y -i \lambda \log \Atil \; , \\
\de(\Atil) &= \Atil \ot \Atil \; , \\
\de(X) &= X \ot \Atil^{-1} + 1 \ot X \; , \\
\de(Y) &= Y \ot \Atil + 1 \ot Y + X \ot (\Atil - \Atil^{-1}) \; .
\end{align*}
Co-unit and antipode are $\vep(\Atil)=1, \vep(X)=\vep(Y)=0,
S(\Atil)=\Atil^{-1}, S(X) = -X\Atil, S(Y)=-Y\Atil^{-1} + X(\Atil - \Atil^{-1})$.
The corresponding Lie bialgebra is
\begin{align*}
[\Atil,X] &= -\Atil \; , \quad [\Atil,Y] = 0 \; , \quad [X,Y] = Y-
\lambda \Atil
\; , \\
\sde(\Atil) &= 0 \; , \quad \sde(X) = \Atil \wedge X  \; ,
\quad \sde(Y) = (2X + Y) \wedge \Atil  \; .
\end{align*}

The dual Hopf $^*$-algebra has self-adjoint generators
$\Xtil,\Ytil$ and an anti-self-adjoint generator $A$, subject to the
relations
\begin{align*}
[\Xtil,A] &= 1-\Xtil+2\Ytil \; , \quad [\Ytil,A] =
\Xtil^{-1}\Ytil(1+\Ytil) \; , \quad [\Xtil,\Ytil]=0 \; , \\
\de(\Xtil) & = \Xtil \ot \Xtil \;  , \\
\de(\Ytil) &= \Xtil \ot \Ytil + \Ytil \ot 1 \; , \\
\de(A) &= A \ot \Xtil^{-1} + 1 \ot A + i \lambda \log \Xtil \ot
\Xtil^{-1} \Ytil \; .
\end{align*}
Co-unit and antipode are $\vep(\Xtil)=1, \vep(\Ytil)=\vep(A)=0,
S(\Xtil)=\Xtil^{-1}, S(\Ytil)=-\Xtil^{-1}\Ytil, S(A)=-A \Xtil + i
\lambda \Ytil \log \Xtil$.
The corresponding Lie bialgebra is determined by the equations
\begin{align*}
[\Xtil,A] &= - \Xtil +2\Ytil \; , \quad [\Ytil, A] = \Ytil
\; , \quad [\Xtil,\Ytil] = 0 \; , \\
\sde(\Xtil) &= 0 \; , \quad \sde(\Ytil) = \Xtil \wedge \Ytil \;
, \quad \sde(A) = \Xtil \wedge A - \lambda \Xtil \wedge \Ytil \; .
\end{align*}

For the exponentiation of the case 4.3 $(a=0)$ described
in Equation~\eqref{exp43a0} with
infinitesimal form $[X,Y]=2Y$, $\chi(X)= -2$, $\chi(Y)=0$, $\be(X)=0$
and $\be(Y) = X$, there are non-trivial cocycles determined by
$\cU(X,Y)= - \lambda \Atil$, giving rise to the following Hopf
$^*$-algebra, with a self-adjoint generator $\Atil$ and
anti-self-adjoint generators $X,Y$, subject to the relations
\begin{align*}
[\Atil,X] &= -2\Atil \; , \quad [\Atil,Y]= - \Atil^2 \; , \quad [X,Y] =
2Y - i \lambda \Atil \; , \\
\de(\Atil) &= \Atil \ot 1 + 1 \ot \Atil \; , \\
\de(X) &= X \ot 1 + 1 \ot X  \; , \\
\de(Y) &= Y \ot 1 + X \ot \Atil + 1 \ot Y \; .
\end{align*}
Co-unit and antipode are $\vep(\Atil)=\vep(X)=\vep(Y)=0,
S(\Atil)=-\Atil, S(X) = -X, S(Y)=-Y + X \Atil$.
The corresponding Lie bialgebra is determined by
\begin{align*}
[\Atil,X] &= -2\Atil \; , \quad [\Atil,Y] = 0 \; , \quad [X,Y] = 2Y-
\lambda \Atil
\; , \\
\sde(\Atil) &= 0 \; , \quad \sde(X) = 0  \; ,
\quad \sde(Y) = X \wedge \Atil \; .
\end{align*}

The dual Hopf $^*$-algebra has self-adjoint generators
$\Xtil,\Ytil$ and an anti-self-adjoint generator $A$, subject to the relations
\begin{align*}
[\Xtil,A] &= \Ytil \; , \quad [\Ytil,A] = 0 \; , \quad [\Xtil,\Ytil]=0 \; , \\
\de(\Xtil) & = \Xtil \ot \Xtil \;  , \\
\de(\Ytil) &= \Xtil \ot \Ytil + \Ytil \ot \Xtil^{-1} \; , \\
\de(A) &= A \ot \Xtil^{-2} + 1 \ot A +i\lambda \Xtil^{-1} \Ytil \ot
\Xtil^{-2} \log \Xtil\; .
\end{align*}
Co-unit and antipode are $\vep(\Xtil)=1, \vep(\Ytil)=\vep(A)=0,
S(\Xtil)=\Xtil^{-1}, S(\Ytil)=-\Ytil, S(A)=-A \Xtil^2 + i
\lambda \Xtil \Ytil \log \Xtil$.
The corresponding Lie bialgebra is determined by
\begin{align*}
[\Xtil,A] &= \Ytil \; , \quad [\Ytil, A] = 0
\; , \quad [\Xtil,\Ytil] = 0 \; , \\
\sde(\Xtil) &= 0 \; , \quad \sde(\Ytil) = 2 \Xtil \wedge \Ytil \;
, \quad \sde(A) = -2 A \wedge \Xtil +\lambda \Ytil \wedge \Xtil \; .
\end{align*}

Finally, to treat case 4.3 $(a<0)$, we use the modification of the
exponentiation \eqref{exp43a-4} with the usage of the universal covering
Lie group (see the very end of Section 5). This gives $[X,Y]=Y$,
$\chi(X)=-1$, $\chi(Y)=0$, $\be(X)=-Y$ and $\be(Y)=X$. There are
non-trivial cocycles, given by $\cU(X,Y)=-\lambda \Atil$. Hence, we
get a Hopf $^*$-algebra with generators $\Atil=\Atil^*$, $X = -
X^*$ and $Y = -Y^*$, subject to the relations
\begin{align*}
[\Atil,X] &= - \sin \Atil \; , \quad [\Atil,Y]= -1 + \cos \Atil \; , \quad [X,Y] =
Y -i \lambda \Atil \; , \\
\de(\Atil) &= \Atil \ot 1 + 1 \ot \Atil \; , \\
\de(X) &= X \ot \cos \Atil - Y \ot \sin \Atil + 1 \ot X \; , \\
\de(Y) &= X \ot \sin \Atil + Y \ot \cos \Atil + 1 \ot Y \; .
\end{align*}
Co-unit and antipode are $\vep(\Atil)=\vep(X)=\vep(Y)=0,
S(\Atil)=-\Atil, S(X) = -X \cos \Atil - Y \sin \Atil, S(Y)=X \sin
\Atil - Y \cos \Atil$.
The corresponding Lie bialgebra is
\begin{align*}
[\Atil,X] &= -\Atil \; , \quad [\Atil,Y] = 0 \; , \quad [X,Y] = Y-
\lambda \Atil
\; , \\
\sde(\Atil) &= 0 \; , \quad \sde(X) = \Atil \wedge Y  \; ,
\quad \sde(Y) = X \wedge \Atil  \; .
\end{align*}

The dual Hopf
$^*$-algebra has self-adjoint generators $\Xtil, \Ytil$ and the
generator $A=-A^*$, subject to the relations
\begin{align*}
[\Xtil,A] &= \Ytil \; , \quad [\Ytil,A] = \frac{1}{2} (\Xtil^{-1}
\Ytil^2 - \Xtil + \Xtil^{-1}) \; , \quad [\Xtil,\Ytil]=0 \; , \\
\de(\Xtil) & = \Xtil \ot \Xtil \;  , \\
\de(\Ytil) &= \Xtil \ot \Ytil + \Ytil \ot 1 \; , \\
\de(A) &= A \ot \Xtil^{-1} + 1 \ot A -i\lambda \log \Xtil \ot
\Xtil^{-1} \Ytil \; .
\end{align*}
Co-unit and antipode are $\vep(\Xtil)=1, \vep(\Ytil)=\vep(A)=0,
S(\Xtil)=\Xtil^{-1}, S(\Ytil)=-\Xtil^{-1}\Ytil, S(A)=-A \Xtil - i
\lambda \Ytil \log \Xtil$.
The corresponding Lie bialgebra is determined by
\begin{align*}
[\Xtil,A] &= \Ytil \; , \quad [\Ytil, A] = -\Xtil
\; , \quad [\Xtil,\Ytil] = 0 \; , \\
\sde(\Xtil) &= 0 \; , \quad \sde(\Ytil) = \Xtil \wedge \Ytil \;
, \quad \sde(A) = \Xtil \wedge A - \lambda \Xtil \wedge \Ytil \; .
\end{align*}

The above list of Lie bialgebras $\lieg$ covers all decomposable Lie bialgebras of
dimension 3, i.e., those that are extensions of the form
$$
\lieg_2^* \longrightarrow \lieg
\longrightarrow \lieg_1.
$$
Naturally, their duals are extensions of the form
$$
\lieg_1^* \longrightarrow \lieg
\longrightarrow \lieg_2.
$$
Here $\lieg_1$ and $\lieg_2$ are the Lie algebras forming the cocycle
matched pair.

\subsubsection*{The case of indecomposable low-dimensional l.c.\ quantum groups}

Comparing the above list with the full list of all mutually
non-isomorphic Lie bialgebras $\lieg$ of dimension 3 obtained by
X.~Gomez \cite{G}, one can observe that there are also
indecomposable Lie bialgebras of dimension 3. Below we list them
and describe the corresponding Hopf $^*$-algebras and, when they
are known, the corresponding locally compact quantum groups.

We start with the Hopf $^*$-algebras $U_q(\su_2)$,
$U_q(\su_{1,1})$, $U_q(\sll_2(\R))$, which are real forms of the
Hopf algebra $U_q(\sll_2(\C))$, and their duals $\SU_q(2)$,
$\SU_q(1,1)$ and $\SL_q(2,\R)$. Hopf algebraically, they are
studied by several authors, see e.g.\ \cite{FRT}.
Operator algebraically, $SU_q(2)$ and its dual are studied in
\cite{VS,W2,W3} as an example of a compact (dually, discrete) quantum
group. $\SL_q(2,\R)$ is studied as a l.c.\ quantum group recently
by W.~Pusz and S.L.~Woronowicz, see \cite{PuszWor} for a related example.
$SU_q(1,1)$ and its dual is treated as a l.c.\ quantum group
in \cite{KK,KK2,W6}.

The Hopf algebra $U_q(\sll_2(\C))$ is defined by 3 generators $a,x,y$
and the following relations ($q$ is a complex parameter):
\begin{align*}
ax &= qxa \; , \quad ay = q^{-1}ya \; , \quad
[x,y]=\frac{1}{q-q^{-1}}(a^2 - a^{-2}) \; , \\
\de(a) &= a \ot a \; , \\
\de(x) &= x \ot a + a^{-1} \ot x \; , \\
\de(y) &= y \ot a + a^{-1} \ot y \; .
\end{align*}
It is well known that there exist three different Hopf $^*$-algebra
structures on $U_q(\sll_2(\C))$: we get $U_q(\su_2)$ by putting $q > 0$, $q \neq 1$, $a=a^*$
and $y=x^*$, we get $U_q(\su_{1,1})$ for the same values of $q$, $a= a^*$ and $x =
-y^*$ and we finally get $U_q(\sll_2(\R))$ for $|q|=1$, $q \neq \pm
1$, $a=a^*, x^*=-x,y^*=-y$.

One can construct the corresponding Lie bialgebras using the above mentioned formal
procedure of linearization. For $U_q(\su_2)$, we
put $H = -i
\frac{1}{\log q} \log a$, $X = i(x+y)$ and $Y=x-y$, and we arrive at the
Lie bialgebra
\begin{align*}
[H,X] &= Y \; , \quad [H,Y]=-X \; , \quad [X,Y] = \frac{8 \log q}{q -
  q^{-1}} H \; , \\
\sde(H) &= 0 \; , \quad \sde(X) = 2 \log q \; H \wedge X \; , \quad
  \sde(Y) = 2 \log q \; H \wedge Y \; .
\end{align*}
For $U_q(\su_{1,1})$, we put $H=-i \frac{1}{\log q} \log a$, $X = x +
y$ and $Y=i(x-y)$ to obtain
\begin{align*}
[H,X] &= -Y \; , \quad [H,Y]=X \; , \quad [X,Y] = \frac{8 \log q}{q -
  q^{-1}} H \; , \\
\sde(H) &= 0 \; , \quad \sde(X) = 2 \log q \; H \wedge X \; , \quad
  \sde(Y) = 2 \log q \; H \wedge Y \; .
\end{align*}
Finally, for $U_q(\sll_2)$, we put $q = \exp(ir)$, $H = - i
\frac{2}{r} \log a$, $X=x$ and $Y=y$ to arrive at
\begin{align*}
[H,X] &= 2X \; , \quad [H,Y]=-2Y \; , \quad [X,Y] = \frac{r}{\sin r} H \; , \\
\sde(H) &= 0 \; , \quad \sde(X) = r H \wedge X \; , \quad
  \sde(Y) = r H \wedge Y \; .
\end{align*}

For the dual Hopf $^*$-algebras, we take the following
versions. $\SU_q(2)$ has generators $a,a^*,b,b^*$, parameter $q > 0$
and relations
\begin{align*}
ab &= qba \; , \quad ab^* = q b^* a \; , \quad bb^* = b^* b \; , \\
aa^*- a^* a &= (q^{-2}-1) b b^* \; , \quad a a^* + b b^* = 1 \; , \\
\de(a) &= a \ot a - q^{-1} b \ot b^* \; , \\
\de(b) &= a \ot b + b \ot a^* \; .
\end{align*}
Suppose $a = \exp(A)$. We first formally calculate that
$$[A,b] = \log q \; b \; , \quad [A^*,b] = - \log q \; b \; .$$
So, we put $H = A - A^*$, $X=i(b+b^*)$ and $Y = b-b^*$ to obtain the
linearization
\begin{align*}
[H,X] &= 2 \log q X \; , \quad [H,Y]= 2 \log q Y \; , \quad [X,Y] = 0 \; , \\
\sde(H) &= q^{-1} X \wedge Y \; , \quad \sde(X) = Y \wedge H \; , \quad
  \sde(Y) = H \wedge X \; .
\end{align*}

Next, we write $\SL_q(2,\R)$ with $|q|=1$ and self-adjoint generators
$a,b,c,d$, subject to the relations
\begin{align*}
ab &= qba \; , \quad ac = q c a \; , \quad bd = q db \; , \quad cd = q
dc \; ,\\
bc &= cb \; , \quad [a,d]=(q-q^{-1})bc \; , \quad ad-qbc=1 \; , \\
\de(a) &= a \ot a + b \ot c \; , \quad
\de(b) = a \ot b + b \ot d \; , \\
\de(c) &= c \ot a + d \ot c \; , \quad
\de(d) = c \ot b + d \ot d \; .
\end{align*}
Again, writing $a = \exp(A)$, $q
= \exp(ir)$, $H=iA$, $X=ib$ and $Y=ic$, we get the Lie bialgebra
\begin{align*}
[H,X] &= -r X \; , \quad [H,Y]= -r Y \; , \quad [X,Y] = 0 \; , \\
\sde(H) &= X \wedge Y \; , \quad \sde(X) = 2 H \wedge X \; , \quad
  \sde(Y) = 2 Y \wedge H \; .
\end{align*}

Finally, we present $\SU_q(1,1)$ with $q > 0$, generators $a,b$ and relations
\begin{align*}
ab &= qba \; , \quad ab^* = q b^* a \; , \quad bb^* = b^* b \; , \\
aa^*- a^* a &= (1-q^{-2}) b b^* \; , \quad a a^* - b b^* = 1 \; , \\
\de(a) &= a \ot a + q^{-1} b \ot b^* \; , \\
\de(b) &= a \ot b + b \ot a^* \; .
\end{align*}
Following the same road as for $\SU_q(2)$, we get the Lie bialgebra
\begin{align*}
[H,X] &= 2 \log q \; X \; , \quad [H,Y]= 2 \log q \; Y \; , \quad [X,Y] = 0 \; , \\
\sde(H) &= q^{-1} Y \wedge X \; , \quad \sde(X) = Y \wedge H \; , \quad
  \sde(Y) = H \wedge X \; .
\end{align*}

The Hopf $^*$-algebras corresponding to the l.c.\ quantum group of
motions of the plane and its dual was considered in e.g.\ \cite{Koe}.
They are treated as l.c.\ quantum groups in
\cite{B}, \cite{VDW}, \cite{VK} and \cite{W4}.

Take $\mu > 0$ and consider the Hopf algebra defined by
$$ax = \mu xa \; , \quad \de(a) = a \ot a \; , \quad \de(x) = a \ot x
+ x \ot a^{-1} \; .$$ We can put two different Hopf $^*$-algebra
structures. First, we get $U_\mu(\mathfrak{e}_2)$ by taking $a$
self-adjoint and $x$ normal. Next, we get $E_\mu(2)$ by supposing that
$a$ is unitary and $x$ is normal. We linearize $U_\mu(\mathfrak{e}_2)$
by writing $H = i \log a$, $X = i(x + x^*)$ and $Y = x-x^*$. This
gives us the Lie bialgebra
\begin{align*}
[H,X] &= -\log \mu \; Y \; , \quad [H,Y]= \log \mu \; X \; , \quad [X,Y] = 0 \; , \\
\sde(H) &= 0 \; , \quad \sde(X) = 2 H \wedge X \; , \quad
  \sde(Y) = 2 H \wedge Y \; .
\end{align*}
For $E_\mu(2)$, we write $H = \log a$ (which is indeed
anti-self-adjoint), $X=i(x+x^*)$ and $Y = x - x^*$, to arrive at the
Lie bialgebra
\begin{align*}
[H,X] &= \log \mu \; X \; , \quad [H,Y]= \log \mu \; Y \; , \quad [X,Y] = 0 \; , \\
\sde(H) &= 0 \; , \quad \sde(X) = 2 Y \wedge H \; , \quad
  \sde(Y) = 2 H \wedge X \; .
\end{align*}

Observe that the list of 3-dimensional Lie bialgebras \cite{G}
contains some more objects, and we now want to present the
corresponding Hopf $^*$-algebras, which are less known. As far as we
know, they have not yet been considered on the level of l.c.\
quantum groups.

Let $\mu \in \C$, $\mu \neq 0$ and let $\rho > 0$. Put $\sla = -
\frac{\log \rho}{\overline{\mu}}$. Then, we can define a Hopf
$^*$-algebra with relations
\begin{align*}
[a,x] &= \mu x \; , \quad xx^* = \rho x^* x \; , \quad
a= a^* \; , \\
\de(a) &= a \ot 1 + 1 \ot a \; , \\
\de(x) &= x \ot \exp(\lambda a) + 1 \ot x \; .
\end{align*}
Co-unit and antipode are given by $\vep(a)=\vep(x)=0, S(a)=-a, S(x)=-x
\exp(-\lambda a)$.
The specific form of $\sla$ is needed to ensure that $\de$ respects
the relation $x x^* = \rho x^* x$. Then, putting $H = i a$, $X = i(x+x^*)$ and $Y =
x-x^*$, and observing that $X$ and $Y$ commute in a first order approximation, we
get the corresponding Lie bialgebra
\begin{align*}
& [H,X] = - \Im \mu \; X - \Re \mu \; Y \; , \quad [H,Y]= \Re \mu
\; X - \Im \mu
\; Y \; , \quad [X,Y] = 0 \; , \\
& \sde(H) = 0 \; , \; \sde(X) = (\Re \sla \; X - \Im \sla \; Y)
\wedge H \; , \;
  \sde(Y) = (\Im \sla \; X + \Re \sla \; Y) \wedge H \; .
\end{align*}
One can check that $\sde$ respects the relation $[X,Y]=0$, because
$\Im (\sla \overline{\mu})=0$. Also, one can check that this
family of Lie bialgebras is self-dual, i.e., the dual of any Lie
bialgebra with specific values of $\mu$ and $\rho$ belongs again
to this family (but with different values of $\mu$ and $\rho$).
So, the dual Hopf $^*$-algebras are of the same form as above.

Next, we take real numbers $\al$ and $\be$, and we write the Hopf
$^*$-algebra with self-adjoint generators $a,x,y$ and relations:
\begin{align*}
[a,x] &= -i x \; , \quad [a,y]=-i \al y \; , \quad
xy = \exp(-i\al\be)yx \; , \\
\de(a) &= a \ot 1 + 1 \ot a \; , \\
\de(x) &= x \ot \exp(\be a) + 1 \ot x \; , \\
\de(y) &= y \ot \exp(- \al \be a) + 1 \ot y \; .
\end{align*}
Co-unit and antipode are given by $S(x)=-x \exp(-\beta a), S(y) =
-y \exp(\alpha \beta a)$, $\vep(a)=\vep(x)=\vep(y)=0$. To
linearize, we write $H = ia$, $X=ix$ and $Y=iy$. Observing again
that $X$ and $Y$ commute in a first order approximation, we
obtain the corresponding Lie bialgebra
\begin{align*}
[H,X] &= X \; , \quad [H,Y]= \al Y \; , \quad [X,Y] = 0 \; , \\
\sde(H) &= 0 \; , \quad \sde(X) = \be X \wedge H \; , \quad
  \sde(Y) =  \al \be H \wedge Y \; .
\end{align*}
In the same sense as in the previous paragraph, this family of Lie
bialgebras is self-dual, so the dual Hopf $^*$-algebras are of the same form.

Finally, there is one isolated Lie bialgebra, which is defined by
\begin{align*}
[H,X] &= 2X \; , \quad [H,Y]= -2 Y \; , \quad [X,Y] = H \; , \\
\sde(H) &= H \wedge Y  \; , \quad \sde(X) =  X \wedge Y \; , \quad
  \sde(Y) =  0 \; .
\end{align*}
We can write the following Hopf $^*$-algebra, which appears
in \cite{CP}, Section~6.4.F  and which has generators $h^* = -h$,
$x=x^*,y=y^*$ and relations
\begin{align*}
[h,x] &= 2x - \frac{1}{2}h^2 \; , \quad [h,y]=2(1-\exp(y)) \; , \quad
[x,y] = h \; , \\
\de(h) &= h \ot \exp(y) + 1 \ot h \; , \\
\de(x) &= x \ot \exp(y) + 1 \ot x \; , \\
\de(y) &= y \ot 1 + 1 \ot y \; .
\end{align*}
One can check that $[x,\exp(y)]=\exp(y) h + \exp(y)(1-\exp(y))$,
so taking $H=h$, $X=-ix$ and $Y = i y$, and linearizing we get the
above Lie bialgebra. For the dual Lie bialgebra, we cannot
construct at the moment a corresponding Hopf $^*$-algebra. The
main problem to construct this exponentiation is the fact that the
dual Lie bialgebra has no non-trivial Lie sub-bialgebra.

\end{document}